\documentclass[a4paper,10pt]{article}

\usepackage{lmodern}

\usepackage{amsfonts}
\usepackage{amssymb}
\usepackage{graphicx}
\usepackage{mathtools}
\usepackage{skak}
\usepackage{MnSymbol}
\usepackage{bbding}
\usepackage{amscd}
\usepackage{MnSymbol}

\usepackage{wasysym}
\usepackage{stmaryrd}
\usepackage{ifsym} 
\usepackage{tgadventor}

\makeatletter

\newcommand{\Rmnum}[1]{\expandafter\@slowromancap\romannumeral #1@}
\makeatother

\title{ Morse Spectra, Homology  Measures, Spaces of Cycles  and Parametric Packing Problems.}

\author{Misha Gromov}

\begin{document}

\maketitle \tableofcontents

\begin{abstract} 

 An "ensemble" ${\Psi}={\Psi}(X)$ of  (finitely or infinitely many) particles in a  space $X$, e.g. in the Euclidean $3$-space, is customary characterised by  the set function 
 $$U\mapsto ent_U({\Psi})=ent({\Psi}_{|U}), \mbox { }  U\subset X,$$
   that assigns  the {\it  entropies} of {\it the $U$-reductions  ${\Psi}_{|U}$ of  ${\Psi}$},  to   all  bounded  open subsets
$U\subset X$. 
In  the physicists' parlance, this entropy is \vspace {1mm}

\hspace {13mm} "{\sl  the  logarithm of the number of  the states of $\cal E$ 
 
 \hspace {19mm}  that are 
 effectively  observable from $U$"},\vspace {1mm}

 \hspace {-6 mm} This     "definition",  in the   context of  mathematical statistical mechanics,  is    translated  to the language  of the measure/probability    theory.\footnote{See:  Lanford's {\sl Entropy and equilibrium states in classical statistical
mechanics, Lecture Notes in Physics, Volume 20, pp. 1-113, 1973}    and Ruelle's   {\sl Thermodynamic formalism : the mathematical structures of classical equilibrium statistical mechanics}, 2nd Edition,
 Cambridge Mathematical Library  2004, where the emphasis is laid upon   (discrete) {\it lattice}    systems.  Also a {\sl categorical rendition} of Boltzmann-Shannon   entropy is suggested in  {\sl "In a Search for a Structure, Part 1: On Entropy"},  www.ihes.fr/$\sim$gromov/PDF/structre-serch-entropy-july5-2012.pdf}

But what happens if     {\sl "effectively observable number of states"}  is replaced  by \vspace {1mm}

 \hspace {10mm}  {\sl "the number of effective/persistent degrees of freedom 
 
 \hspace {23mm}  of ensembles of moving particles"?}\vspace {1mm}

We suggest in this paper     several    mathematical counterparts to the idea of   {\sl " persistent degrees of freedom "} and  formulate specific questions, many of which  are   inspired by Larry Guth's  results  and ideas on   
the Hermann Weyl kind of  asymptotics of {\it the Morse (co)homology spectra}  of the {\it volume energy function} on the spaces of   {\it  cycles}  in  balls.\footnote {{\sl Minimax problems related to cup powers and Steenrod squares},
Geometric and Functional Analysis, 18 (6), 1917-1987 (2009).}
And often we present variable aspects of the same idea in     different sections of this paper.

Hardly anything that can be called  "new theorem"  can be found in our paper but we reshuffle   many known results and expose them from a  particular angle. This article is meant as   an introductory chapter to something    yet to be written with    
 much  of what we present here  
extracted from my yet unfinished manuscript {\sl Number of 
Questions}.

\end{abstract}

\section  {Overview  of Concepts and Examples.}

We    introduce below  the idea of {\sl "parametric packing"} and of related concepts which are  expanded  in detail in the rest of the paper.\vspace {1mm}

\textbf {A.} Let $X$ be a topological space, e.g. a manifold,  and $I$ is a countable index set  that may be  finite, especially if  $X$ is compact.

A  collection of $I$-tuples of  non-empty open (sometimes closed) subsets  $U_i\subset X$, $i\in I$,  is called {\it a packing} or an {\it $I$-packing} of $X$ if these subsets {\it do not intersect}. 

Denote by  ${\Psi}(X;I)$ the space of these packings with some natural topology, where, observe  there are several candidates for such a topology  if $X$ is {\it non-compact}.

\textbf {B.} {\it Homotopies  Constrained by  Inradii  and Waists.}  We are interested in the homotopy, especially (co)homology, properties of  subspaces  ${\cal P}\subset {\Psi}(X;I)$ defined by imposing  {\it lower bounds}   on   the sizes of $U_i$.
where the  two such  invariants of $U\subset X$  we shall be often (albeit sometimes implicitly) use in this paper are the  {\it the inradius   of $U$} and the {\it the $k$-waist of $U$, $k=1,2,...,dim(X)-1$},
 defined later on in the metric and the symplectic categories.
\vspace {1mm}

 \textbf {C.} {\it Metric Category and  Packings by Balls.} Let $X$ be a metric space, let   $r_i\geq 0$ be non-negative numbers and take metric balls in $X$ of radii $R_i\geq r_i$ for $U_i$. Packings by such balls are traditionally called {\it sphere packings} where one is especially concerned with  packing homogeneous spaces (e.g. spheres and Euclidean spaces)   by {\it equal} balls. 

 The corresponding  space ${\cal P}= {\cal P}(X;\{\geq r_i\}_{i\in I}) $ naturally embeds into  the  {\it Cartesian power  space}  
 $$X^I=  \underbrace { X\times X \times ...\times X}_I$$  
 
 of $I$-tuples of points   $X$  where it is distinguished  by the inequalities 
 $$dist(x_i,x_j)\geq d_{ij}=r_i+r_j, \mbox { } i,j\in I, i\neq j,$$
and where, observe all these spaces with $d_{ij}>0$ lie in the Cartesian power space   $X^I$ minus diagonals,
 $$X^I\setminus \bigcup_{i,j\in I}Diag_{ij}\mbox { for $ Diag_{ij} \subset X^I$ defined by the equations $dist(x_i,x_j)=0$}.$$
 
And if $X$ is Riemannian manifold with {\it the convexity  (injectivity?) radius} $\geq R$, then, clearly, the 
inclusion  $${\cal P}(X;\{\geq r_i\}_{i\in I}) \hookrightarrow X^I\setminus \bigcup_{i,j\in I}Diag_{ij}$$
is {\it a  homotopy equivalence} for $\sum_{i\in I}r_i < R/2$;
moreover, if all $r_i$ are mutually equal, this 
homotopy equivalence is {\it equivariant}   for  the permutation group that acts on $I$ and thus on    $X^I$ and on $X^I\setminus \bigcup_{i,j\in I}Diag_{ij}\subset X^I$ and  ${\cal P}(X;\{\geq r_i=r\})\subset X^I\setminus \bigcup_{i,j\in I}Diag_{ij} $.\vspace {2mm}

The packings spaces  {\it covarinatly functorially} behave under {\it expanding} maps 
between metric spaces $X\to Y$.  

 But   {\it contravariant}  functoriality under  {\it  contracting}, i.e. {\it distance decreasing}, maps $f: X\to Y$   needs the following  extension  of the concept of ball packings to   $I$-tuples of {\it subsets} (rather than points)  $V_i\subset X$  instead of points $x_i\in X$.
 \vspace {1mm}
 
 \textbf {D.}  {\it Packings by Tubes. }  These are    $I$-tuples of closed  subsets $V_i\in X$, such that   mutual distances\footnote{Recall that $dist (V_1,V_2)$  between two subsets in a metric space $X$  is defined as the {\it infimum} of the distances between points $x_1\in V_1$ and $x_2\in V_2$.}  between them satisfy $dist(V_i,V_j)\geq d_{ij}$.

  \textbf {E.} {\it Packing by Cycles.}  The above  becomes  interesting if all
$V_i\subset X$  support   given {\it nonzero}   homology classes $h_i$ of dimension $k$, $k=0,1,...,$ in $X$, or if they support $k$-cycles some of which are linked in $X$,  either individually or "parmetrically" (compare  \textbf {P} below and  section 18)

\vspace {1mm}
 
 \textbf {F.} {\it Packings by Maps.} This is yet  another variation of the concept of "packing".
Here  {\sl subsets}  $U_i\subset X$  are replaced by   {\it maps} $\psi_i:U\to X$, where in general, the domain $U$ of $\psi_i$  may depend on $i\in I$.

Now, "packing by $\psi_i$" means packing by the images of these maps, i.e. these images should not intersect  with  specific "packing  conditions" expressed in terms of geometry of $U$ and of these maps. 

For instance, 
if $X$ and $U$ are equidimensional Riemannian  manifolds on may  require the maps $\psi_i$ to be expanding.  
Or $\psi_i$  may  belong to  a  particular  category   (pseudogroup) of  maps.

\vspace {1mm}

 \textbf {G.} {\it Symplectic  Packings by Balls. }   Here  $X=(X,\omega) $ is a  $2m$-dimensional symplectic manifold,  $U_i$ are balls of radii $R_i$ in the standard  symplectic space $\mathbb R^{2n}=(\mathbb R^2)^n$ and symplectic  packings  are given by $I$-tuples of symplectic embeddings  $U_i\to X$ with disjoint images.

Another attractive   class of symplectic packings is that by polydiscs and by  {\it $R$-tubes around  Lagrangian 
submanifolds} in $X$.

\vspace {1mm}

 \textbf {H.}  {\it Holomorphic Packings}.  These make sense and look interesting  for packing by holomorphic maps
$U_i\to \mathbb C^n$, $n\geq 2$,  with {\it Jacobians one}, in particular, by {\it symplectic holomorphic} maps for $n$ even,
but I have not thought about these.  


\vspace {2mm}

\textbf {I}. {\it Essential Homotopy.} This is  the part of the  homotopy, e.g. cohomology, structure of a  {\it geometrically} defined  subspace  ${\cal P}\subset {\Psi}$  that comes from the  
ambient space ${\Psi}$, where $\Psi$ itself is  defined in a purely {\it topological} terms.

One think of  ${\Psi}$  as the background  that supports the {\it geometric information} on $\cal P$ written  in the  {\it homotopy theoretic} language of $\Psi$.

This information concerns  the {\it relative homotopy size} of $\cal P$ in ${\Psi}$, often  expressed by particular (quasi)numerical  invariants, such, for instance,  as {\it homotopy hight, cell numbers, cohomology valued measures.}

\vspace {1mm}
\textbf {J.} {\it Example: Packings by Two Balls.} The  space ${\cal P}(R) \subset X\times X$, $R\geq 0$,  of  packings of a metric space $X$ by two $R$-balls is  defined by the inequality
 $$ X\times X\supset {\cal P}(R)
=_{def} \{x_1,x_2\}_{dist(x_1,x_2)\geq 2R},   $$   
where, observe  the distance function $d$  on $X\times X$ is related to the distance in $X\times X$   to the diagonal  $X_{dia}=Diag_{12}\subset X\times X$ by
 $$d=dist_X(x_1,x_2)=\sqrt 2 \cdot dist_{X\times X} ((x_1,x_2), X_{dia}).$$

The {\it algebraic   topology} of the distance function  $d$
on $X\times X$, more specifically the (co)homologies  of the {\it inter-levels} 
$$d^{-1}[R_1,R_2]\subset X\times X,$$
(that carries, in general,  more {\it geometric} information about $X$  than
what we call "essential homotopy" of these subsets) can be 
thought of  (be it essential or non-essential)  as {\it  a cohomology valued "measure-like"  set function} on the real line,  namely
 $$[a,b]\mapsto H^\ast(d^{-1}[a,b]) \mbox { for all segments  $[a, b]\subset \mathbb R$.}$$

Exceptionally, e.g. if $X$ is a symmetric space,  the distance function is "{\sl Morse perfect}"\footnote{"Morse perfect" functions do not have to be "textbooks Morse": they    may be {\it non-smooth} and have {\it positive dimensional sets of    critical points}.} :  all of homotopy topology (e.g. homology)  of $f$ is "essential",  quite transparent instances of which are  {\it   projective spaces over $\mathbb R,\mathbb C$ and $\mathbb H$} as well  {\it flat tori}, where these functions obviously induced from  (Morse perfect distance) functions $d_0:X\to \mathbb R$
for $d_0(x)=dist(x,x_0)$.

On the other hands  much of geometrically significant {\it geometric information} carried  by {\it the topology} of the distance functions on manifolds of {\it negative curvature,} does not  quite conform  to such  concept of "essential". 
(The simplest part of the information encoded by  $d$, that  is represented by the set of   the lengths of  {\it undistorted}\footnote{A submanifold $Y$ in a Riemannian manifold $X$ is called   {\it undistorted} if  the distance in $Y$ associated with the {\it induced Riemannian} structure coincide with the {\it restriction of the  distance function}  from $X\supset Y$. For instance the shortest non-contractible curve in $X$ is undistorted.}
  closed geodesics in $X$, is homotopically, but not necessarily homologically, essential.)

\vspace {2mm}

\textbf {K.} {\it Permutation Symmetries.}  The Cartesian power space $X^I$ is    acted upon by 
the symmetric group 
$\mathbb S_N=Sym(I)$, $N=card(I)$  and that this action is {\it free}  in the complement to the diagonals.

Let $G $ be a subgroup in 
$\mathbb S_N=Sym(I)$, $N=card(I)$, e.g. $G=\mathbb S_N$,   and let $\cal P\subset {\Psi}$ be a subspace  invariant under this action. 
We are especially interested in those  homotopy   characteristics   of $\cal P$  that are  encoded by the kernel of the  cohomology homomorphism $\kappa^\ast: H^\ast( BG)\to H^\ast( {\cal P}$$\slash$$G )$ for  the classifying map $\kappa$ from  the {\it  quotient space}\footnote{ Our  $\cal P$ is contained in the complement to the  diagonals in $X$ and, hence, the action of $G$  is free on  ${\cal P}$; otherwise, we  would replace the quotient  space ${\cal P}\slash$$G$ by  the {\it homotopy quotient} ${\cal P}$$\sslash$$G$.}   ${\cal P}$$\slash$$ G$ to the  classifying space $BG$ of the group $G$,
$$\kappa:  \mbox {${\cal P}$$\slash$$G \to BG.$}$$

\vspace {1mm}

Some of the  questions we want to know the answer to are  as follows:\footnote{The cohomology of $\mathbb S_N$  are well understood, (e.g. see \cite {adem})  but I shamefully  failed to extract a rough estimate of  the ranks $rank(H^i(\mathbb S_N,\mathbb F_p)$ from what I read. But even if these ranks are bound by $const^N$, the minimal number of cells in the optimal cell decorposition of the classifying space $B\mathbb S_N$ must(?) grow  roughly as $N!$  for $N\to \infty$.}\vspace {1mm}

  \textbf {L.}  Let $F_0$ be a  $G$-equivariant map from a topological space $S$ with an action of a group $G=G(I)\subset \mathbb S_N$, $N=card(I)$, on it, to the Cartesian power  space $X^I$, 
$$F_0:S\to X^I, $$
and let $[F_0]_{G}$ denote the equivarinat homotopy class of this map.

What is the maximal radius
 $R=R(S,[F_0]),$
 such that $F_0$ admits an equivariant homotopy to an equivariant map $F$ from $S$ to the space ${\cal P}(X;I,R)\subset X^I$, of $I$-packings of $X$ by balls of radii $R$,
 $$F: S\to {\cal P}(X;I,R)\subset X^I,  \mbox { } F\in [F_0]_{\mathbb S_N}?$$

What is the supremum $R_{max}(k)$ of these $R= R(S,[F_0])$, over all $(S,[F_0])$ with $dim(S)=k$?

In other words, what is the maximal $R=R_{max}(k)=R(X,N,k)$, such that every $k$-dimensional  $G$-invariant subset in the complement of the diagonals in $X^I$ (where the action of  $\mathbb S_N\supset G$ is free)  admits an equivariant homotopy to    ${\cal P}(X;I,R)\subset X^I$? 

\vspace {2mm}

 \textbf {M.}  Let $$K^\ast =K^\ast (N,\varepsilon)\subset H^\ast( BG; \mathbb F_p),\mbox { }  \mathbb F_p = Z/p\mathbb Z,\mbox {  for $G=G(N), N=card(I)$}$$
 be the kernels   of the homomrphisms   $$\kappa\ast= \kappa_ {\mathbb F_p}^\ast: H^\ast( BG;  \mathbb F_p )\to H^\ast( {\cal P}/G; \mathbb F_p ),$$
where 
${\cal P}={\cal P}(X; I, \varepsilon)$ is  the space of packings of a given Riemannian manifold $X$, e.g. of the  unit sphere $S^l$ or the torus $\mathbb T^l$, by $N$ balls of radius $\varepsilon.$
  
  \vspace {1mm}
  
   (\textbf {M$_1$}) {\sl What is the  behaviour  of the (graded) ranks of these  kernels 
  $K^\ast (N,\varepsilon)$ as functions of $\varepsilon$?}
    \vspace {1mm}
  
   (\textbf {M$_2$})  {\sl What is the asymptotic   behaviour  of these  $rank(K^\ast (N,\varepsilon))$ for $N\to \infty$ and $\varepsilon \to 0$,
   where a particular case of interest is that of  $$\varepsilon=\varepsilon_N=const\cdot N^\alpha\mbox { for some $\alpha<0?$}$$}

  \vspace {1mm}

\textbf {N.} {\it Packing Energies and  Morse Packing Spectra.} The space 
${\cal P}={\cal P}_\varepsilon ={\cal P}(X; I, \varepsilon)$ can be seen  as {\it the sublevel} of a suitable "{\sl energy function}" $E$ on the ambient space ${\Psi}=X^I\supset \cal P_\varepsilon$,
where any  
   monotone decreasing function in $$\rho(\psi)\mbox { for 
   $\psi=\{x_i\}$ and 
$\rho(\psi)=\min_{x_i\neq x_j}dist (x_i,x_j)$}$$ 
will do\footnote { The  role of real numbers $\mathbb R$  here reduces  to  indexing the subsets  $\Psi_r\subset \Psi$, $r\in \mathbb R$, according to their order by inclusion:  $\Psi_{r_1}\subset \Psi_{r_2}$ for  $r_1 \leq r_2$. 

In fact,  our   "spectra" make sense for functions with values  in an arbitrary {\it lattice}  (that is a partially ordered set  that  admits   $inf$ and $sup$), while 
 {\sl additivity}, that is the most essential feature of the physical energy,   becomes visible only for spaces  $\Psi$   that   split as  $\Psi=\Psi_1\times \Psi_2$.}
  for    $E(\psi_1,\psi_2)=E(\psi_1)+E(\psi_2)$  and  where  simple candidates for such functions are
$$E(\psi)=\frac{1} {\rho(\psi)}\mbox { or  $E(\psi)=- \rho(\psi),$}$$
or, that seems  most appropriate from  a certain perspective,
 $$E(\psi)=-\log  {\rho(\psi)}.$$ 

 Notice that $\rho(\psi)$ equals $\sqrt 2$ times the distance  from $\{x_i\}\in X^I$ to the union of the  diagonals $Diag_{ij}\subset X^I$ that are defined by the equations $x_i=x_j$,
 $$\rho(\psi)=\sqrt 2\cdot \min_{ij}dist_{X^I}(\psi, Diag_{ij}), \mbox { } \psi=\{x_i\}.$$

\vspace {2mm}
 \textbf {N$^{\star}$.}  We are  predominantly interested in {\it the  homotopy significant  (Morse) spectra} of such  energy functions $E:{\Psi}\to \mathbb R$,  on  topological spaces $\Psi$, where such a spectrum  is the set of  those values $e\in \mathbb R$ where the homotopy type of the  sublevel $E^{-1}(-\infty, e]$
undergoes an {\it irreversible change} (precise  definitions are given in section 4)  
 and the above  (\textbf {M$_1$}) concerns such changes that are recored by the variation of the kernels  $K^\ast (N,\varepsilon)$.

 \vspace {1mm}
 
 \textbf {O.} {\it "Duality" between Homology Spectra of Packings and of Cycles.}
 Evaluation of the  homotopy (or homology)  spectrum  of packings, in terms  of the above \textbf {A}, needs establishing two opposite geometric inequalities, similarly  how it goes  for the spectra of {\it Laplacians} associated to {\it Dirichlet's energies.}
 \vspace {1mm}
 
 \textbf {O$_{\mathbf I}$}:{\it  Upper Bounds on Packing Spectra.} Such a  bound for  packings a  manifold $X$ by $R$-balls, means  an  inequality  $R\leq \rho_{up}$ for some $\rho_{up}=\rho_{up}(S, F_0)$ (or several such inequalities for various $S$ and $F_0$), that  would guarantee that a  map $F_0:S\to X^I$ {\it is homotopic (or at least  homologous)}  to a map with image in the packing space ${\cal P}(X;I,R)\subset X^I$.

This, in all(?)  known examples,  is achieved by {\it explicit constructions}  of specific "homotopically  (or homologically) significant"  packings families $P_s\in {\cal P}(X;I,R)$, $s\in S$, for  $R\geq \rho_{up}$.
  
   \vspace {0.8mm}
  
 \textbf { O$_{\mathbf {II}}$} :{\it  Lower Spectral  Bounds.} Such a  bound $R\geq \rho_{low}= \rho_{low}(S, F_0)$ is supposes to signify  that $F_0:S\to X^I$ is {\it not homotopic (or not even homologous)} to a map  $S\to {\cal P}(X;I,R)\subset X^I$.
  
 All (?) known bounds of this kind are  obtained by {\it  (parametric) homological localisation}  that is by confronting such  maps $F_0$, think of them as   {\it families of   $N=card(I)$ moving balls     in $X$ parametrised by $S$}, with {\it families of cycles} in $X$ where the two families have a  nontrivial (co)homology pairing between them. 
   
   A  simple, yet instructive,  instance of this is where: \vspace {1mm}

there is  
{\it a  (necessarily non-zero) homology class}
$h_\circ\in H_k (X)$ for some 

$k=1,2,...,n=dim (X)$,  
such that \vspace {0.5mm}

 the image   $h_S\in H_\ast(X^I)$ of some homology class from $H_\ast(S)$ under the induced homomorphism $$(F_0)_\ast :  H_\ast(S) \to H_\ast(X^I)$$ 
 have {\it non-zero homology intersection   with the power class
$h_\circ^{\otimes I}\in  H_\ast(X^I)$},
$$    h_S\frown h_\circ^{\otimes I}\neq 0.$$
 
Obviously, in this case, \vspace {0.7mm}

 {\it if  a map $F_0:S\to X^I$ is homotopic to a map into the packing  space  ${\cal P}(X; I,R)\subset X^I$, then every closed subset $Y\subset X$ that supports the class $h_\circ$ admits a packing by $N$-balls  $U_i\subset Y$, $i\in I$,   $N=card (I)$,  of radii $R$  for the restriction of the  distance function from $X$ to $Y$.}  
 {\it Consequently, $$N\cdot  R^{k} \leq  const _X<\infty.$$}
   
\vspace {1mm}   

 {\it Example.}  Let $X$ be a Riemannian product of two closed connected  Riemannian manifolds, $X=Y\times Z$, let  $S=Z^I$ and $F_0:Z^I=(Z\times y_0)^I \subset  X^I$, $ y_0\in Y$, be the tautological   embedding.   
 If this $S\subset X^I$ can be moved by a homotopy in $X^I$  to   ${\cal P}(X; I,R)\subset X^I$, then $Y$ can be packed by $N$-balls of radius $R$. 
 
 Notice that the converse is also true in this case. In fact,    if   balls $U_{y_i}(R) \subset Y$   
 pack $Y$, then  the Cartesian  product 
 $$S=\bigtimes _{i\in I}(Z\times y_i) \subset X^I$$
 is contained in ${\cal P}(X; I,R)$.

    \vspace {1mm}
      
   In general, when cycle moves, this kind of argument,    besides suitable  
   nontrivial (co)homology pairing, needs  lower bounds on spectra of the {\it volume-energies} in  spaces of $k$-cycles, in particular lower bounds on {\it $k$-waists} of our manifolds that correspond to the bottoms of such spectra.
   
In particular, such  bounds on {\it symplectic waists}, are used in the symplectic geometry for proving \vspace {0.8mm} 

{\it non-existence  of individual   packings, as well as of  multi parametric 
families 

of  certain symplectic packings}.     
   
       \vspace {2 mm}
 \textbf { P.}  {\it Packings by Tubes around Cycles.} The concepts of   spaces of packings and those  of cycles can be brought to  the common ground  by introducing    the space  of $I$-tuples of {\it disjoint $k$-cycles} $V_i$ in $X$ and  of (the  homotopy spectrum of) the  function $E\{V_i\}$ that would  somehow incorporate   $vol_k(V_i)$ along with $\log dist_X(V_i,V_j)$.\footnote  {One may think of these $V_i$ as images of   $k$-manifolds mapped to  $X$  that faithfully corresponds to cycles with $\mathbb Z_2$-coefficients.} 

(One may replace  the distances $dist_X(V_i,V_j)$ in $X$  by distances in {\it the flat metric} in the space of cycles, where $dist_{flat}(V_i,V_j)$  is defined as  the $(k+1)$-volume of {\it the minimal  $(k+1)$-chain }  between $V_i$ and $V_j$ in $X$, but this would lead to a quite different picture.)     
     
      \vspace {1mm}
  
  \vspace {1mm}
 
  \textbf { Q.} {\it  Spaces of Infinite Packings.} If $X$ is a non-compact manifold, e.g. the Euclidean space $\mathbb R^n$, then there are many  candidates for {\sc the space of packings}, all of which are infinite dimensional spaces with infinite dimensional (co)homologies, where this infinities  may be  (partly) offset by  actions of  infinite groups  
 on these spaces.

 For instance, spaces   $\cal P$ of packings of $\mathbb R^n$ by countable sets
 of $R$-balls $U_i(R)\subset \mathbb R^n$ are acted upon by the isometry group $iso(\mathbb R^n)$ that, observe , commute with  the action of the group $Sym(I)$ of bijective transformations of the (infinite countable) set $I$.

The simplest(?) instance of an interesting infinite packing space ${\cal P}={\cal P}(\mathbb R^n;I,R)$ is  where
$I=\mathbb Z^n\subset \mathbb R^n$ is the integer lattice, where the ambient space 
$\Psi$ equals the space of {\it bounded displacements} of $\mathbb Z^n\subset \mathbb R^n$ that are  maps $\psi:i\mapsto x_i\in \mathbb R^n$, such that
$$dist(i,x_i)\leq C=C(\psi)<\infty \mbox { for all } i\in \mathbb Z^n$$  and where ${\cal P}\subset \Psi$
is distinguished by the inequalities $dist(x_i,x_j)\geq 2R.$

 The essential part of the infinite dimensionality of this  $\Psi$ comes from the infinite  group $\Upsilon=\mathbb Z^n\subset iso(\mathbb R^n)$, that acts on it. 
  In fact,  $\Psi$ naturally (and $\Upsilon$-equivarinatly) imbeds into the union of compact infinite product spaces 
 $${\Psi}\subset \bigcup_{C>0} B(C)^{\Upsilon},$$
 where $ B(C)\subset \mathbb R^n$ is the Euclidean ball of radius $C$ with the centre at the origin, and where  several  entropy-like topological invariants, such as the {\it mean dimension} $dim({\cal P}): \Upsilon$ and {\it polynomial  entropy} $H^\ast({\cal P}):\Upsilon$, are available.\vspace {1mm}

The quotient space of the above  space $\Psi$, or rather of this $\Psi$ minus the diagonals, by the infinite "permutation group"\footnote{One only needs here the subgroup of  $Sym(I=\mathbb Z^n)$ that consists  of {\it bounded} bijective displacements    $\mathbb Z^n \to \mathbb Z^n$, i.e. where these "displacements"  $i\mapsto j$ satisfy $dist(i,j)\leq C<\infty$. } $Sym(I)$  consists of the set of certain discrete  subsets  
$\Delta \subset \mathbb R^n$. These $\Delta$ has the property that the intersections of it with all bounded open subsets  $V\subset \mathbb R^n$ satisfy {\it the uniform density condition}, 
    $$card(V\cap \Delta)-card (V\cap \mathbb Z^n)\leq const_\Delta\cdot vol (U_1(\partial V)) \leqno{\mbox {\tiny \AsteriskRoundedEnds}}$$
where $U_1(\partial V)\subset \mathbb R^n$ denotes the union of the unit balls with their centres in the boundary of $V$.

But there are by far more uniformly dense (i.e. that satisfy {\tiny \AsteriskRoundedEnds}) subsets  than the images of $\mathbb Z^n$ under bounded displacement.

In fact, it is far from clear  what   topology (or rather homotopy) structure should be used  in the  space of uniformly dense subsets, that, for instance, would render this space {\it (path)connected.} 

  \textbf {R.} {\it On Stochasticity.} The traditional probability may be brought back  to this picture, if, for instance, 
packing spaces are defined by the inequalities
      $$dist(x_i, x_j)\geq \rho(i,j),$$
where   $\rho(i,j)$   assumes two values, $R_1>0$ and $R_2>R_1$, taken independently with given  probabilities $p(i,j)=p(i-j)$, $i,j\in I=\Upsilon=\mathbb Z^n$.

\vspace {2mm}
  \textbf { S.} {\it From Packings to Partitions and Back.}  An  $I$-packing $P$ of a metric space $X$ by subsets $U_i$ can be canonically  extended to the corresponding {\it Dirichlet-Voronoi partition}\footnote{Here, "{\sl $I$-Partition}"  means  a covering of $X$ by closed subsets $V_i\subset X$, $i\in I$,   with non-empty non-intersecting interiors, where we often tacitly assume certain regularity of  the boundaries of these  $U_i$.   } 
by subsets $U_i^+\supset U_i$, where each  $U_i^+\subset X$, $i\in I$, consists of the points $x\in X$ nearest to $U_i$, i.e. such that
$$dist(x,U_i)\leq dist(X,U_j), \mbox { } j\neq i.$$
If, for instance, $X$ a  convex\footnote {{"\sl Convex"} means that every two points are joint by a unique geodesic.}  Riemannian space with constant curvature and if all $U_i$ are convex then  $U_i^+$ are  {\it convex polyhedral sets.}  

Conversely, convex subsets $U\subset X$  can be often {\it canonically shrunk} to single points $u\in U$ by families of {\it convex} subsets, $U_t\subset U$, $0\leq t \leq 1$,
$$\mbox {  where $U_0=U$, $U_1=u$  and  $U_{t_2}\subset U_{t_1}$ for $t_2\geq t_1$ }.$$  For instance, if $X$ has non-positive curvature, such a shrinking can be accomplished  with the  {\it inward equidistante deformation} of the boundary $\partial U$.

This shows, in particular,  that  the space of convex  $I$-partitions (as well as of convex $I$-packings)  of  a convex space $X$ of constant curvature is {\it $\mathbb S_N$-equivarinatly, $N=card(I)$,   homotopy equivalent} to the space of $I$-tuples of distinct points $x_i\in X$. (The case of negative curvature reduces to that of the positive one via projective isomorphisms  between bounded  convex  spaces of constant curvatures.)

 Partitions of metric spaces, especially convex ones whenever these are available, reflects  finer aspects the geometry of $X$ than sphere packings. For instance,  families  of convex partitions  obtained by consecutive division of convex sets by hyperplanes are used for sharp evaluation of {\it waists} of spheres as we shall explain later on. 

On the other hand, a typical Riemannian manifold  $X$ of dimension $n\geq 3$  admits  only approximately convex  partitions (along with convex packings), where  the geometric significance  of these remains problematic.
\vspace {2mm}

  \textbf { T.} {\it  Composition of Packings, (Multi)Categories and Operads.}  If  $U_i$ pack $X$ then packings of $U_i$, $i\in I$,  by $U_{ij}$, $j\in J_i$, define a packing of $X$ by all these $U_{ij}$. 

Thus, for instance,  in the case of  $X$ and all $U$ being Euclidean balls, this composition defines a {\it  topological/homotopy operad} structure in the space of packings of balls by balls. 

The significance of such a  structure is questionable for round ball packings, especially for those of {\it high density},  since composed packings constitute only a small part of the space of all packings.

 But   packing  {\it of cubes by smaller  cubes} and    {\it symplectic} packings of balls by    smaller balls seem more promising in this respect, since  even quite  dense packings in these cases, even partitions,  may  have a significant amount of (persistent) degrees of freedom.\footnote{Besides composition, there are other operations on (nested) packings.  For instance, (close  to each other) large  balls (cubes) may  "exchange" small balls  (cubes) in them.}

\vspace {2mm}

 \textbf { U.} {\it  On   Faithfulness of the "Infinitesimal Packings" Functor.} The  (multi)category  structures in spaces of packings define, in the limit,  similar structures in spaces of packings of spaces $X$  by "infinitely many infinitesimally small" subsets $U_i\subset X$.  

\smallskip

{\it Question \textbf {U1}}. How much of the geometry of a (compact) space  $X$, say with a metric or symplectic geometry, can be seen in the homotopies  of spaces of packings of $X$ by such $U_i$?  

{\it Question \textbf {U2}.}  Is there a good   category of  "abstract packing-like objects",  that are not, a prori, associated to actual packings of geometric spaces?

\vspace {1mm}

Concerning  Question \textbf {U1}, notice that the above mentioned pairing between "cycles" and packings, shows that the volumes of certain {\it minimal subvarieties} in a Riemannian manifold $X$, can be  reconstructed from the homotopies of packings of $X$ by arbitrarily  small balls.\vspace {1mm}

For instance, \vspace {1mm}

if $X$ is a complete Riemannian manifold with {\it non-negative sectional  curvatures}, then {\it the lengths of its closed geodesics} are (easily) seen in the homotopy spectra  of these packing spaces. (see section 9)

And  if  $X$ is an {\it orientable surface}, then this remains true with {\it no  assumption on the curvatures} for the  geodesics that are {\it length minimising in their respective homotopy classes}. 

Similarly, much of the  geometry of {\it waists} of a   convex set $X$, say in the sphere $S^n$,   may be seen in the homotopies of spaces of partition of $X$ into convex subsets, see \cite {waists}.

\vspace {1mm} 

{\it Question \textbf {U3}.} Would it be useful to enhance the homotopy structure of a  packing  space  of an $X$, say  by (infinitesimally) small balls, by keeping track of  (infinitesimal) geometric sizes of the homotopies in such a space?

\vspace {2mm}

 \textbf { V. }  {\it Limited Intersections.}   A similar to packings  (somewhat  less interesting?)  space is that   of $N$-tuples of balls with  {\it no  $k$-multiple intersections} between them,  (this space  contains the space of $(k-1)$-tuples  of packings)  can be seen with the distance function  to the union of the   {\it $k$-diagonals} --  there are  $N\choose k$ of them -- that are the pullbacks of the principal diagonal $\{x_1=x_2=...=x_k\}$ in $X^J$, where $card(J)=k$,   under  maps $X^I\to X^J$ corresponding to   $N!/(N-k)! $ imbedding  $J\to I$.

\textbf { W. } {\it Spaces of Coverings.} Individual packing often go together with coverings, say, with  minimal covering of metric spaces by $r$-balls. Possibly, this companionship extends to that between {\it spaces} of packings and {\it spaces} of coverings.

\section {A few Words on Non-Parametric Packings.}

Classically, one is concerned with  {\it maximally dense} packings of  spaces $X$ by {\it disjoint  balls}, rather than with the homotopy properties of families of moving  balls  in $X$.

Recall, that {\it a sphere packing} or, more precisely,   {\it a  packing of a metric space $X$ by  balls   of radii $r_i$, $i\in I$, $r_i>0$}, for a given {\it indexing set} $I$  of finite or countable cardinality $N=card(I)$  is, by definition,  a collection of  (closed or open) balls  $U_{x_i}(r_i)\subset X$, $x_i\in X$,  with mutually non-intersecting interiors. 

Obvioulsy, points $x_i\in X$ serve as centres of such balls if and only if 
$$dist(x_i,x_j)\geq d_{ij}=r_i+r_j.$$

{\it Basic Problem.}  What is the  {\it maximal radius}  $r=r_{max}(X;N)$ such   that $X$ admits a packing by $N$   balls of radius $r$?
 
 In particular, \vspace {0.7mm}
 
\hspace {10mm}  {\it what is the asymptotics of  $r_{max}(X;N)$  for $N\to 
 \infty$?}  \vspace {0.7mm}
 
If $X$ is a compact $n$-dimensional Riemannian manifold (possibly with boundary), then the principal term of this asymptotics depends  only on the volume of $X$, namely, one has the following (nearly obvious) \vspace {1mm}
      
      \hspace {25mm}   {\sc  Asymptotic Packing Equality.}
$$\lim_{N\to \infty}\frac{N\cdot r_{max}(X;N)^n }{vol_n(X)} =
 \circledcirc_n,$$
where  $\circledcirc_n>0$ is a  universal (i.e. independent of $X$) {\it Euclidean packing constant}   that corresponds in an obvious way to the {\it optimal density} of the    sphere packings  of the  Euclidean space $\mathbb R^n$.
  
(Probably, the full asymptotic expansion of $r_{max}(X;N)_{N\to \infty}$ is  expressible in terms of the  derivatives of the curvature of $X$ and derivatives the curvature  similarly to Minakshisundaram-Pleijel   formulae for spectral asymptotics.)

The explicit value of $\circledcirc_m$ is known only for   $n=1,2,3$. In fact,  
the  optimal, i.e. maximal,  packing  density  of $\mathbb R^n$ for $n\leq 3$  can be implemented by a {\it $\mathbb Z^m$-periodic}  (i.e. invariant under some  discrete action of $\mathbb Z^n$ on $\mathbb R^n$)  packing, where the case $n=1$ is obvious, the case $n=2$ is due to Lagrange (who proved that  the optimal packing is the hexagonal one)  and 
the case of $n=3$, conjectured by Kepler, 
was resolved by Thomas  Hales.

(Notice  that   $\mathbb R^3$, unlike $\mathbb R^2$ where {\it the  only}  densest packing is the hexagonal one,  admits  {\it infinitely many}  different  packings;   most of these are {\it not $\mathbb Z^3$-periodic}, albeit they are $\mathbb Z_2$-periodic.

Probably,   none of    densest   packings of $\mathbb R^n$   is  
 $\mathbb Z^n$-periodic  for large $m$, possibly  for  $n\geq  4$. Moreover,  {\it the topological  entropy} of the action of $\mathbb R^n$ on the space of optimal packings may be non-zero. 

 Also,  there may be infinitely many 
 algebraically independent numbers  among 
 $\circledcirc_1,\circledcirc_  2, ...$; moreover,  the  number of algebraically independent among $\circledcirc_1,\circledcirc_  2, ..., \circledcirc_n$  may grow  as
 $const\cdot n$, $const>0$.)

.

\section {Homological Interpretation of     the  Dirichlet-Laplace  Spectrum.}  

Let $\Psi$ be  a topological space and $E:\Psi\to\mathbb R$ a continuous  real valued function, that is  thought of as an energy $E(\psi)$ of states $\psi\in \Psi$ or as a Morse-like function on $\Psi$.

The subsets $$\Psi_e=\Psi_{\leq e}=E^{-1}(\-\infty, e]\subset \Psi,\mbox {  } r\in \mathbb R,$$ are called the (closed) {\it $e$-sublevels of} $E$.

A number $e_\circ\in \mathbb R$ is said {\it to lie in the homotopy significant spectrum of} $E$ if the homotopy type of  $\Psi_r$ undergoes a  {\it significant}, that is  {\it irreversible}, change as $e$ passes through  the value $e=e_\circ$, that may be understood as {\it non-existence of a homotopy} of the subset  $\Psi_{e_\circ}$  in $\Psi$ that would bring it to the sublevel $\Psi_{e<e_\circ}\subset \Psi_{e_\circ} $.

\vspace {1mm}

 {\it Basic  Quadratic Example.} Let  $\Psi$ be an infinite dimensional projective space and $E$ equal the ratio of two quadratic functionals. More specifically,  let $E_{Dir}$ be
the  Dirichlet function(al) on differentiable functions $\psi=a(x)$ normalised by the $L_2$-norm  on a  compact Riemannian manifold $X$,
$$E_{Dir}(\psi)=\frac{||d\psi||^2_{L_2}}{||\psi||^2_{L_2}}=\frac {\int_X ||d\psi(x)||^2dx}{\int_X \psi^2(x)dx}.$$

The eigenvalues  $e_0, e_1, e_2,...,e_N,...$ of   $E_{Dir}$ (i.e. of the  corresponding Laplace operator)  are  {\sl homotopy significat} since 
the rank of the inclusion homology homomorphism $H_\ast(\Psi_{r}; \mathbb Z_2)\to H_\ast(\Psi; \mathbb Z_2)$ {\it  strictly} increases (for $\ast=N$) as $e$ passes through $e_N$.\vspace {1mm}

 An essential feature of Dirichlet energy that, as we shall see, is shared by many other examples is {\it homological localisation}.
 
 Let  $X$ be partitioned by closed  subsets $U_i$, $i\in I$, $card (I)=n$, with piecewise smooth boundaries. Then
 
{\it  the $N$-th  eigenvalue $e_N=e_N(X)$  is bounded  from below by
 the minimum of

  the first 
 eigenvalues of $U_i$, 
 $$e_N(X)\geq \min _{i\in I} e_1(U_i).$$}

 Indeed, by linear algebra,  every $N$-dimensional projective space of functions on $X$, say  $S=P^N\subset \Psi=P^\infty$,
 contains a (necessarily non-zero) function $\psi_\star(x)$ such that 
 $$\int _{U_i}\psi_\star(x)dx=0\mbox   { for all $i\in I$}.\leqno {\star_i}$$ 
 Therefore, $$\sup_{a\in P^n}E_{Dir}(\psi) \geq  \max _{i\in I} E_{Dir}(\psi_{\star |U_i})\geq e_i,$$
where the key feature of this argument. \smallskip

\hspace {10mm} -- {\it simultaneous solvability of the equations} $\star_i$ -- 
 
\hspace {0mm}   does not truly  need the linear structure in $P^n$.
   but only the fact that  \smallskip
   
 \hspace {6mm}   {\it  $S\subset \Psi$  supports a  cohomology class $h\in H^1 (\Psi;\mathbb Z_2)=\mathbb Z_2=\mathbb Z/2\mathbb Z$, 
   
   \hspace {19mm}with non vanishing $\smile$power $h^{\smile N}\in H^N(\Psi;\mathbb Z_2)$.}
\smallskip
 
 If $U_i$ are   small  approximately round subsets, then  
 $$e_1(U_i)\geq \varepsilon_n \left({\frac {1}{vol(U_i)}}\right)^{\frac {2}{n}},   \mbox {} n=dim X,$$
and, with suitable partitions into such subsets,  one bounds $e_N(X)$  
from below\footnote{It is    obvious that  $e_N(X)\leq C_X\cdot  \left (\frac {N}{vol(X)}\right)^{\frac {2}{n}}.$  In fact, the numbers $N_{sp\leq e}(U)$ of the eigenvalues $e_i(U)\leq e$ of  open subsets $U\subset X$ satisfy {\it Hermann Weyl's asymptotic formula}  
$N_{sp\leq e}(U)\asymp D_n vol(U)^\frac{n}{2} $, $n=dim(X)$, where  the existence of the limits  $ N_\asymp(U)$ of $ N_{sp\leq e}(U)e^\frac{2}{n}$ for $e\to \infty$  and additivity of the set function $U\to N_\asymp(U)$ follows  from the locality of the $\smile$product,  while the  evaluation of  $D_n$,  that happens to be equal  $2\pi^{-n} vol(B^n(1))$ for $B^n(1)$ being the unit  Euclidean ball, depends on the (Riemannian)  geometry  of the Dirichlet energy.   }
 by
    $$e_N(X)\geq \varepsilon_X\cdot  \left (\frac {N}{vol(X)}\right)^{\frac {2}{n}}.$$

This  can  be quantified in terms of the geometry of $X$. 

 For instance, if   $Ricci(X)\geq -n\delta$, $\delta\geq 0$, and $diam(X)\leq D$, then    the {\it homology localisation} for {\it $(n-1)$-volume energy $$E_{vol}(\psi)=vol_{n-1}(\psi^{-1}(0))$$}  in conjunction with  Cheeger's  spectral  inequality implies that 
 $$ e_N(X)\geq \varepsilon_n^{1+D\sqrt{\delta}}D^{-2}N^{\frac {2}{n+1}}.$$
  (See section 7  and   "Paul Levy Appendix" in    \cite{metric}.)
  
  \section{Induced Energies on  (Co)Homotopies on  (Co)Homologies.}

{\it On  Stable and Unstable Critical points.} If $E$ is a   Morse function on a smooth manifold $\Psi$, then the homotopy type of  the energy  sublevels $\Psi_{e}=E^{-1}(-\infty, e]\subset \Psi$  does change at all critical values $e_{cri}$  of $E$. However, only exceptionally rarely, for the so called {\it perfect Morse functions},   such as for the above quadratic energies, these changes are irreversible. In fact, every value $r_0\in \mathbb R$   can be made critical by an {\it arbitrarily
small $C^0$-perturbation\footnote{"$C^0$" refers to the uniform topology in the space of continuos functions.} $E'$} of a smooth function $E(\varphi)$, such that $E'$ equals $E$ outside the subset  $E^{-1}[r_0-\varepsilon, r_0+\varepsilon]\subset \Psi$; thus, the topology change of the sublevels of $E'$  at $r_0$ is  insignificant.

But the   spectra of  Morse-like functions introduced below  have such  homotopy significance  built into their very  definitions.  \vspace {1mm}

\vspace {1mm}

{\it ${\cal H}_\circ(\Psi)$, $E_\circ$  and  the Homotopy Spectrum.}  Let $\cal S$ be a class of topologicl spaces $S$  and  let ${\cal H}_\circ(\Psi)={\cal H}_\circ (\Psi;{\cal S})$ be the category   where  the objects  are homotopy classes of continuous maps $ \phi: S\to \Psi$ and morphisms are homotopy classes of maps  $\varphi_{12}:S_1\to S_2$, such that the corresponding  triangular diagrams are  (homotopy) commutative, i.e.  the composed maps
$\phi_2\circ\varphi_{12}:S_1\to \Psi$  are homotopic to $\phi_1$.

 Extend functions  $E:\Psi\to \mathbb R$ from $\Psi$    to    ${\cal H}_\circ(\Psi)$  as follows.   
Given a continuous map   $ 
\phi: S\to \Psi$ let
  $$ E(\phi)=E_{max}(\phi)=\sup_{s\in S} E\circ \phi(s),$$
denote by  $ [\phi]= [\phi]_{hmt}$  the homotopy class of $\phi$.
and set
$$ E_\circ[\phi]=E_{mnmx}[\phi]= \inf_{\phi\in [\phi]}E(\phi).$$  
 In other words, \vspace {1mm}
 
\hspace {-5mm}{\sl $E_\circ[\phi]\leq e\in \mathbb R$ if and only if the map  $\phi=\phi_0$ admits a homotopy of maps

\hspace {-5mm}$\phi_t: S\to A,$ $0\leq t\leq 1,$ 
 such that $\phi_1$  sends $S$ to the sublevel $\Psi_e=E^{-1}(-\infty ,e]\subset \Psi$.} \vspace {1mm}

  {\it The   covariant (homotopy) ${\cal S}$-spectrum of} $E$ is the set of values  $E_\circ[\phi]$  for some class  $\cal S$  of (homotopy types of) topological spaces $S$ and (all)  continuous maps $\phi :S\to \Psi$. 

For instance, one may take for $\cal S$ the set of  homemorphism classes  of  countable (or just  finite)  cellular spaces. In fact,   the set of sublevels $\Psi_r$, $r\in \mathbb R,$  themselves is  sufficient for most  purposes. \vspace {0.7mm}

{\it Lower and Upper Bounds on Spectra.} Lower bounds on  homotopy spectra 
say, in effect, that "homotopically {\it complicated/large}" maps $\phi$ (that may need  complicated parameter spaces $S$ supporting them) necessarily  have {\it large energies   $E_\circ[\phi]$}. 

Conversely,  upper bounds depend on  construction of {\it complicated} $\phi :S\to \Psi$ with  {\it small energies.} \vspace {0.7mm}

 {\it On Topology,  Homotopy and on Semisimplicial Spaces.} The topology  of  a space  $\Psi$ per se is not required  for the definition of homotopy (and cohomotopy below)  spectra.  What is needed  is a {\sl "homotopy structure" in $\Psi$}  defined by distinguishing a  class of maps from
  "simple spaces" $S$ into $\Psi$. 
  
If such a structure in $\Psi$ is associated  with polyhedra taken for  "simple $S$", then      $\Psi$  is called a   
{\it semisimlicial (homotopy)  space}  with its "homotopyy structure" defined via   the contravariant  functor $S  \rightsquigarrow maps (S\to \Psi)$  from the category of simplicial complexes and simplicial maps to the category of sets. \vspace {1mm}

{\it ${\cal H}^\circ(\Psi)$,  $E^\circ$  and the  Cohomotopy $\cal S$-Spectra.} Now, instead of  ${\cal H}_\circ(\Psi)$ we extend $E$ to the category  ${\cal H}^\circ(\Psi)$ of homotopy classes of maps $ \varphi :\Psi\to T$, $T\in \cal S$, by defining $E^\circ[\varphi ]$ as the supremum of those $e\in \mathbb R$  for which the restriction map of $\varphi$ to the  energy sublevel $\Psi_e = E^{-1}(-\infty, e]\subset \Psi$, 
$$\varphi_{|\Psi_e}:\Psi_e\to T,$$
is  {\it contractible.}\footnote{In some cases,  e.g. for maps $\varphi$ into discrete spaces $T$ such as Eilenberg-MacLane's $K(\Pi;0)$, {\sl "contractible"}, must be replaced by {\sl"contractible to a marked point serving as zero"} in $T$ that is expressed in writing as $[\varphi]=0$. }
Then the set of the values $E^\circ[\varphi]$, is called {\it the contravariant homotopy (or cohomotopy) $\cal S$-spectrum of $E$.}

For instance, if $\cal S$ is comprised  of the {\it Eilenberg-MacLane   $K(\Pi,n)$-spaces}, $n=1,2,3,...$, then this is called {\it the  $\Pi$-cohomology spectrum of} $E$.\vspace {1mm}

{\it Relaxing Contractibility  via Cohomotopy Operations.} Let us express "contractible" in writing as $[\varphi]=0$, let $\sigma: 
T\to T'$  be a continuous map and  let us  regard  the (homotopy classes of the) compositions of $\sigma$  with  $\varphi:\Psi\to T$  as an   {\it operation} $[\varphi]\overset {\sigma}\mapsto [\sigma\circ\varphi]$.

Then  define  $E^\circ[\varphi]_\sigma\geq E^\circ[\varphi]$  by maximising 
 over those $e$ where   $[\sigma\circ\varphi_{|\Psi_e}]=0$ rather than $[\varphi_{\Psi_e}]=0$.

{\it Pairing between Homotopy and Cohomotopy.} Given a  pair of maps $(\phi, \varphi)$,  where   $\phi:S\to \Psi$ and $\varphi:\Psi\to T$, write\vspace {0.7mm}

{\it $[\varphi\circ \phi]=0$ if the composed map $S\to T$ is contractible,

$[\varphi\circ \phi]\neq0$ otherwise.}\vspace {0.7mm}

\hspace {-6mm}Think of this as a function with value $"0"$ and $"\neq 0"$ on these pairs.\footnote{If the space $T$ is {\it disconnected}, it should be better endowed with a marking $t_0\in T$ with "contractible"  understood as "contractible to $t_0$".}\vspace {1mm}

{\it  $E_\ast $,   $E^\ast $ and  the (Co)homology Spectra.} 
If $h$ is a homology class in the space $\Psi$ then $E_\ast (h) $ denotes {\it the infimum} of  $E_\circ [\phi] $ over all (homotopy classes) of maps $\phi:S\to \Psi$ such that $h$ is {\it contaned in the image} of the  homology homomorphism induced by $\phi$.

Dually,  the energy   $E^\ast (h) $ on a cohomology class $h\in H^\ast (\Psi;\Pi)$ for an Abelian group $\Pi$, is defined as  $E^\circ [\varphi_h] $ for the $h$-inducing map from $\Psi$ to the product of Eilenberg-MacLane spaces:
 $$\varphi_h:\Psi\to \bigtimes_nK(\Pi,n), \mbox { } n=0,1, 2,... \mbox { }  .$$ 
 In simple words,  $E^\ast (h) $ equals
the supremum of those $e$ for which $h$ vanishes on $\Psi_e=E^{-1}(\infty,e]\subset \Psi$.\footnote {The definitions of energy on homology and cohomology obviously extend to  generalised  homology and  and cohomology theories.}  

Then one defines the (co)homology spectra as the sets of values of these energies    $E_\ast $ and  $E^\ast $ on homology and on cohomology.

\vspace {1mm}

{\it Homotopy Dimension (Height) Growth. } The roughest invariant one wishes to extract from the (co)homotopy spectra of an energy  $E:\Psi\to\mathbb R$ is {\it the rate of the growth of the homotopy dimension} of the sublevels $\Psi_e=E^{-1}(-\infty,e]\subset \Psi $,   where the homotopy dimension a subset  $B\subset A$  is the minimal $d$ such that $B$, or at least every polyhedral   space $P$ mapped to  $B$, is  contractible in $A$ to a subset $Q\subset A$ of dimension $d$.\footnote  {This is called {\it essential dimension} in \cite {non-linear} and it equals the homotopy height of (the homotopy class   of) the inclusion $B\hookrightarrow A$.}

In many cases,   this dimension is known to satisfy a polynomial bound $homdim(\Psi_e)\leq c e^\delta$ for some constants $c=c(E)$ and $\delta=\delta(E)$,
where such an inequality  amounts to {\it a lower} bound on the spectrum of $E$.

In the simplest case  of  $\Psi$  homotopy equivalent to $P^\infty$, this dimension as function of $e$  carries {\it all} spectral information about $E$.   

 For instance if $E$ is the     Dirichlet energy of  function  on a Riemannian  $n$-manifold $X$ where the  eigenvalues are  bounded from below by
$e_N(X)\geq \varepsilon_X {N}^\frac {2}{n}$,
one has $homdim(\Psi_e)\leq c \cdot e^\frac {n}{2}$  for $c= \varepsilon_X^{-\frac{n}{2}}$.

\vspace {2mm}

{\it Multidimensional Spectra.} Let   ${\cal E}=\{E_j\}_{j\in J}:\Psi\to\mathbb R^J$ be a continuous map. Let $h$ be a cohomology class
of  $\Psi$ and   define the {\it spectral hypersurface} $\Sigma_h\subset\mathbb R^J$ in the Euclidean space $\mathbb R^J=\mathbb R^{l=card(J)}$ as the boundary  of the subset $\Omega_h\subset \mathbb R^J$ of the $J$-tuples of numbers 
$(e_j)$ such that the  class $h$ vanishes on the subset 
$\Psi_{<e_j}\subset \Psi$ defined by the inequalities 
$$E_j(\psi)<e_j,\mbox { } j\in J.$$
$$ \Sigma_h=\partial \Omega_h,\mbox { } \Omega_h=  \{e_j\}_{ h | \bigtimes_{j\in J} \Psi_{e_j}=0}.$$

(This also make sense for general cohomotopy classes $h$ on  $\Psi$ with $h=0$ understood as contractibility of the map $\psi:\Psi\to T$ that represent $h$  to a marked "zero"  point in $T$ where marking is unnecessary for connected spaces $T$.)

More generally, given a continuous map ${\cal E}: \Psi\to Z$,  one  "measures" open subsets $U\subset Z$ according to the  sizes of "the parts" of the  homology of $\Psi$ that are "contained" in ${\cal E}^{-1}(U)\subset \Psi$,  that are the images of the homology homomorphism $H_\ast (U)\to  H_\ast(\Psi)$;   similarly,  kernels  of the cohomology homorphisms $H^\ast(\Psi \setminus U)\to  H^\ast(\Psi)$ serve  as a measure-like function $U\mapsto \mu^\ast(U)$ on $Z$, (see section 11).

If $ \Psi\to Z$ are smooth manifolds, and $\cal E$ is a proper  smooth map, then the  the set $\Sigma_{\cal E}\subset Z$ of critical values of $\cal E$  "cuts" $Z$ into subsets where the measure $\mu^\ast$ is (nearly) constant. (It is truly  constant on   the connected components of  $Z\setminus \Sigma_{\cal E}$ but may vary at the boundaries of these subsets, since these boundaries are contained in $\Sigma_{\cal E}$.)
Here,  "the cohomology spectrum of  $\cal E$" should be  somehow defined via  "coarse-graining(s)" of the "partition"  of $Z$ into these subsets according to the values of $\mu^\ast$.

\vspace {1mm}

\vspace {1mm}

{\it Packing Example. } Take   $\Psi$ equal the $I$-Cartesian power of a Riemannian manifold, $\Psi=X^I$, let   $J$ consist of unordered pairs  $(i_1,i_2)$ $i_1\neq i_2$, thus 
$card (J)=l= N(N-1)/2$, $N=card (I)$,  and let 
${\cal E}$ be given by the reciprocals of  the $l$ functions $ E_j=dist_X(x_{i_1},x_{i_2})$. 
(This map is equivariant for the natural actions of the permutation group $Sym_N=aut(I)$ on $\Psi $ and on $\mathbb R^J$ and the most interesting aspects of the topology of this $\cal E$ that are  indicated below  become visible   only in the  equivariant setting of 
 section 12) 

\vspace {1mm}

\vspace {1mm}

{\it Spectral Families.}  A similar (dual?) picture arises when one has a family of   functions  
$E_z: \Psi_z\to \mathbb R$ parametrised by a topological space  $Z\ni z$,  where the family $homospec_z\subset \mathbb R$ of  homotopy spectra of $E_z$ is seen as  the  {\it spectral hypersurface} $\Sigma$ of $\{E_z\}_{z\in Z},$
$$\Sigma= \bigcup_{z\in Z}homospec_z \subset  Z\times \mathbb R.$$

\vspace {1mm}

{\it On Positive and Negative Spectra.} Our  definitions of homotopy and homology spectra are best  adapted to  functions $E(\psi)$ bounded from below but they can be adjusted  to more general   functions $E$  such  as  $E(x)=\sum_k a_kx_k^2$  where there may be infinitely many negative  as well as positive numbers among    $\psi_k$.

For instance, one may define the  spectrum  of a  $E$ unbounded from below  as the limit of the homotopy spectra of  $E_\sigma =E_\sigma(\psi)=\max(E(\psi), -\sigma) $ for $\sigma\to +\infty$.

But often, e.g. for the action-like functions in the symplectic geometry, one needs something more sophisticated than a simple minded cut-off of "undesirable infinities".\footnote{ It seems,  however, that neither a  general theory nor a comprehensive list of examples  exit for the moment.}

\vspace {1mm}

{\it On  Continuous Homotopy Spectra.} There also  is a  homotopy theoretic rendition/generalisation of {\it continuous spectra}  with  some Fredholm-like notion of homotopy,\footnote{ See {\sl On the uniqueness of degree in infinite dimension} by P. Benevieri and M. Furi,
http://sugarcane.icmc.usp.br/PDFs/icmc-giugno2013-short.pdf.} such that, for instance,  the natural inclusion of the projectivised Hilbert subspace  $PL_2[0, t]\subset  PL_2[0, 1]$, $0<t<1$, would  not contract to any  $PL_2[0, t-\varepsilon]$.

\section{  Families of $k$-Dimensional Entities,  Spaces of Cycles and  Spectra of   
$k$-Volume Energies.}

The  {\sl "$k$-dimensional size"} of  a metric space $X$, e.g. of a Riemannian manifold, may be  defined   in terms of 
 the (co)homotopy or  the (co)homology spectrum of the  {\it $k$-volume} or  a similar energy function on a {\sl space of "virtually $k$-dimensional entities"}  $Y$ in  $X$. These spectra and related invariants of $X$ can be defined by means  of families of such "entities" as in the following  examples. \vspace {1mm}
 
 
 (A) If $X$ and  $S$ are  topological spaces of dimensions $n=dim(X) $,  and $m=n-k=dim(S)$ 
 then continuous maps $\varsigma:X\to S$ define $S$-families of the fibres $Y_s=\varsigma^{-1}(s)\subset X$.   \vspace {1mm}

 (B) Given a pair of spaces $T$ and $T_0\subset T$, where $dim(T_0)=dim(T)-m$, and   an $S$-family of maps $\phi_s:X\to T$, $s\in S$, one arrives at an  $S$-family of  "virtually $m$-codimensional"  subsets in $X$ by taking the  pullbacks $Y_s=\phi_s^{-1}\subset X$.  \vspace {1mm}
 
 (C) In the case of  {\it smooth manifolds}  $X$ and $S$ and {\it generic smooth maps} $\varsigma :X\to S$, the   families  from (A) and similarly for (B) can be seen  geometrically  as  maps from $S$ to the  {\it space of $k$-manifolds} with the natural (homotopy) semisimplicial  structure defined by bordisms.
 
 And for general continuous maps  $\varsigma$ we think of $Y_s=\varsigma^{-1}(s)\subset X$ as an $S$-families of {\it virtual} $k$-(sub)manifolds in $X$.

\vspace {1mm}

(D) More generally, if $X$ and $S$ are {\it pseudomanifolds}  of dimensions $n$ and $m=n-k$ and  $\varsigma :X\to S$  is a {\it simplicial} map, then  the fibres $Y_s\subset X$ are  $k$-dimensional pseudomanifods for all $s$ in $S$ away from the $(m-1)$-skeleton of $S$  with the map $s\mapsto Y_s$ being  semisimlicial for 
 the natural semisimplicial structure in the space of pseudomanifolds. 
 
  \vspace {1mm}
  
 (E) In order to admit  families of {\it mutually intersecting} subsets $Y_s\subset X$ we need an auxiliary space $\Sigma$ mapped to $S$ by $\varsigma:\Sigma \to S$. Then we  let $\tilde Y_s=\varsigma^{-1}(s)\subset \Sigma$   and  define  $S$-families $Y_s\subset X$ via maps 
 $\chi:\Sigma\to X$ by taking images $Y_s=\chi(\tilde Y_s)\subset X$.
  
 ($\mathrm{\tilde E}$)  Such a $\Sigma$ may be constructed starting from a family of subsets $Y_s\subset X$  as the subset   $\Sigma=\tilde X\subset X\times S $ that consists of the pairs $(x,s)$ such that $x\in Y_s$. However, this  $\tilde X\to X$, unlike more general $\Sigma\to X$, does not account for possible self-intersections of  $\tilde Y_s$ mapped to $X$.)

 \vspace {1mm}

(F) If we want to keep track of multiplicities of maps $\tilde Y_s\to Y_s\subset X$ it is worthwhile to regard the maps 
$$\chi_{|\tilde Y_s}:\tilde Y_s\to X$$
themselves, rather tan their images, as our (virtual)  "$k$-dimensional entities in $X$".\vspace {1mm}

(G= C+F) The space of smooth $k$-submanifolds in an $n$-dimensional manifold  $X$ can be represented by the space of continuous maps $\sigma$ from  $X$ to the Thom space $T$  of the inversal $m$-dimensional vector bundle ofr $m=n-k$, where  "virtual $k$-submanifolds" in $X$ come as the $\sigma$-pullbacks of the zero section of this bundle.

Then {\it the space of bordisms} associated to   arbitrary topological space  $X$ may be defined with $n$-manifolds  $\Sigma$,   maps $\sigma: \Sigma\to T$ and  maps $\chi:\Sigma\to X$. \vspace {2mm}

                       \hspace {30mm} {\sc Space ${\cal C}_k(X;\Pi)$ of $k$-Cycles in $X$.}\vspace {1mm}

 There are several homotopy equivalent candidates for  {\it the space of $k$-cycles with $\Pi$-coefficients}  of a topological space $X$.
 
For instance, one may apply the   construction   of a {\it semisimplicial space of $k$-cycles} associated to a chain complex of Abelian groups (see      section 2.2 of \cite {singularities2})  to the complex  $\{C_i\overset {\partial_i}\to C_{i-1}\}_{i=0,1,...,k,...}$  of {\it singular $\Pi$-chains} that are  $\sum_\nu\pi_\nu\sigma_\nu$  for $ \pi_\nu\in \Pi$, where  $\Pi$ is an Abelian group  and where   $\sigma_\nu:\triangle^i\to X$ are continuous maps of the $i$-simplex to $X$.

For instance, $k$-dimensional psedomanifolds mapped to $X$ define singular $\mathbb Z_2$-cycles in $X$ ($\mathbb Z$-cycles if  $\Sigma$ and 
 $S$ are oriented) where the  above $S$-families agree with the    semisimplicial structure in ${\cal C}_k(X;\mathbb Z_2)$. 
 \vspace {1mm}
 
More generally,   $l$-chains in the space of $k$-cycles in $X$ can be represented by  
$l$-dimensional families of $k$-cycles in X that are    $(k + l)$-chains in $X$. This   defines a map between the corresponding  spaces of cycles $$\top :{\cal C}_l({\cal C}_k(X;\Pi);\Pi))\to {\cal C}_{k+l}(X;\Pi\otimes_\mathbb Z \Pi)$$
 for all Abelian (coefficient) groups $\Pi$,  as well as a   natural homomorphism of degree $-k$
from the cohomology of a space $X$ to the cohomology of the space of $k$-cycles in $X$ with $\Pi$-coefficients, denoted 
$$\top ^{[-k]}: H^{n}(X;\Pi)\to  H^{n-k}({\cal C}_k(X;\Pi);\Pi), \mbox{  for all  $ n\geq k$  },$$
provided  $\Pi=\mathbb Z$ or $\Pi=\mathbb Z_p=\mathbb Z/p\mathbb Z$.




 \vspace {2mm}

 {\it Almgren-Dold-Thom Theorem for the Spaces of Cycles.} Let $X$ be a   triangulated space and  $f:X\to \mathbb R^m$ a generic piece-wise linear  map.  Then the "slicings" of  
 generic $(m+k)$-cycles $V\subset X$ by the fibres of $f$, define  homomorphisms from 
the homology groups  $H_{m+k}(X;\Pi)$ to the homotopy groups $\pi_m({\cal C}_k(X;\Pi))$ of the spaces   ${\cal C}_k(X;\Pi)$ of  $k$-dimensional $\Pi$-cycles in    
   $X$ for all Abelian groups $\Pi$.

  In fact, the map $r\mapsto f^{-1}(r)\cap V$, $r\in \mathbb R^m$,  sends $\mathbb R^m$ to ${\cal C}_k(X;\Pi)$, where the complement to the (compact!)  image $f(V)\subset \mathbb R^m$ goes to the the zero cycle.

The   Almgren-Dold-Thom Theorem claims that these  homomorphisms 
     $$H_{m+k}(X;\Pi)\to\pi_m({\cal C}_k(X;\Pi))$$
 are isomorphisms.   This is easy for the above semisimlicial spaces  (see \cite {singularities2}) and, as it was shown by Almgren, this remains true for  
 spaces of {\it rectifiable cycles with  flat topology.}

(Recall, that the  {\it flat/filling distance}  between two homologous $k$-cycles $C_1,C_2$ in $X$  with $\mathbb Z$- or $\mathbb Z_p$-coefficients 
is defined as the  infimum of $(k+1)$-volumes (see below) of $(k+1)$-chains  $D$ such that $\partial D=C_1-C_2$.

Also notice, that If $\Pi$ is a field, then 
    ${\cal C}_k(X;\Pi)$  {\it splits into 
   the product} of {\it the Eilenberg-MacLane spaces} corresponding to the homology
   groups of $X$, 
$${\cal C}_k(X;\Pi)=\bigtimes_{n=0,1,2,...} K(\Pi_n, n),  \mbox { where }  \Pi_n = H_{k+n}(X;\Pi),$$
 where, observe,  $H^{n} (K(\Pi, n);\Pi)$ is {\it canonically} isomorphic to $\Pi$ for the cyclic groups $\Pi$.)

\hspace {12mm} {\sc $k$-Volume and Volume-like  Energies.}  \vspace {1mm}

$\bullet$ {\it Hausdorff Measures of Spaces and Maps.} The  $k$-dimensional Hausdorff measure of a {\it semimetric}\footnote{"Semi"   allows   a possibility  of  $dist(y_1,y_2)=0$ for $y_1\neq y_2$. }  space $Y$  is defined for all positive real numbers $k\geq 0$ as 
    $$Haumes_k(Y)=\beta_k \cdot \inf_{\{r_i\}} \sum_{i\in I}r_i^{-k},$$
  where $I$ is a countable set, 
  $\beta_k
     = \frac{\pi^{k/2}}{\Gamma(\frac{k}{2} + 1)} $ 
  is the normalising constant that, for integer k, equals    
      the volume of the unit   ball  $B^k$  and where the infimum is taken over all  $I$-tuples   $\{r_i\}$ of positive numbers, such 
that $Y$ admits an $I$-covering  by balls   of radii $r_i$.

The corresponding  {\it Hausdorff measure of a map $f:Y\to X$}, where 
$X$ is a metric space, is, by definition, 
 the Hausdorff measure of  $Y$ with the  
semimetric induced by  $f$ from $X$.

If  $f$ is one-to-one then 
$Haumes(f)= Haumes(f(Y))$
but, if  $f$ has a "significant multiplicity" then $Haumes(f)> Haumes f(Y)$.
 
  \vspace {1mm}
 
 {\it  $\delta$-Neighbourhoods $U_\delta (Y)\subset X$ and Minkowski Volume.}
  Minkowski $k$-volume of a subset $Y$ in an $n$-dimensional Riemannian manifold $X$ and/or in a similar space with a distinguished measure regarded as the $n$-volume,  is   defined as
 $$Mink_k(Y)= \liminf_{\delta\to 0}   \frac{vol_n (U_\delta(Y))}{\delta^{n-k}vol_{n-k}(B^{n-k})}.$$  \vspace {1mm}

In general, the the Minkowski  $k$-volume may be much smaller than the Hausdorff $k$-measure but the two are equal for  "regular" subsets $Y\subset X$ where "regular"   includes

$\bullet$ compact smooth and piecewise smooth submanifolds in smooth  manifolds;
 
 $\bullet$  compact real analytic and semianalitic subspaces in real analytic spaces; 
  
  $\bullet$ compact minimal subvarieties in Riemannian manifolds.

  \vspace{1mm}
 
Besides  Minkowski volumes themselves, 

\hspace {6mm}{\it the normalised $n$-volumes of the $\delta$-neughbourhoods  of    subsets $Y\subset X$},

$$\delta\mbox{-}Mink_k(Y)=_{def}   \frac{vol_n (U_\delta(Y))}{\delta^{n-k}vol_{n-k}(B^{n-k})},$$
also can be used as "volume-like energies" that have interesting  homotopy spectra. 

The pleasant, albeit obvious,  feature of the volume  $\delta$-$Mink_k(Y)$ for $\delta>0$   is its {\it continuity} with respect to {\it the Hausdorff metric} in the space of subsets $Y\subset X$.  \vspace{1mm}
   
 {\it On $\delta$-Covers and $\delta$-Packings of $Y$.}   The  {\it minimal} number of 
 $\delta$-balls  {\it needed to cover}  $Y$   provides an  alternative to
  $\delta$-$Mink_k(Y)$ and the {\it maximal} number of {\it disjoint}  $\delta$-balls in $X$ with {\it centers in} $Y$ plays a similar role.    
 
 The definitions of these  numbers  {\it does not depend} on  $vol_n$; yet, they are closely (and obviously) related  to  $\delta$-$Mink_k(Y)$.
   
\section {Minmax Volumes and Waists.}
   
   Granted a space $\Psi$  of "{\sl virtually $k$-dimensional entities}" in $X$ and a volume-like energy function $E=E_{vol_k}: \Psi\to \mathbb R_+$, 
   "the first eigenvalue" -- the  bottom of the homotopy/(co)homology spectrum of  this $E$   
is called the {\it $k$-waist of $X$}.


 \vspace {1mm}

\vspace {1mm}

To be concrete, we define waist(s) below via the  two basic operations  of  producing $S$-parametric families of subsets -- taking  pullbacks and images of  maps represented by the  following diagrams 
$${\cal D}_X=\{X\overset {\chi} \leftarrow \Sigma\overset {\varsigma}\to S\}$$
where $S$ and $\Sigma$ are  simplicial  (i.e. triangulated topological) spaces of dimensions $m=\dim(S)$ and $m+k=dim (\Sigma)$  and where  our "entities"  are the pullbacks 
$$\tilde Y_s=\varsigma^{-1}(s)\subset \Sigma,\mbox { } s\in S,$$ 
that are  mapped to $X$ by $\chi$. \vspace{1mm} 

  \hspace {40mm}    {\sc Definitions.} \vspace{1mm} 
 
[\textbf {A}] {\it The maximal  $k$-volume} -- be it  Hausdorff, Minkowski, $\delta$-Minkowski, etc,--  of such a family is defined  as the supremum of the  of the corresponding volumes  of the image restrictions of the map $\chi$ to $Y_s$, that is 
  $$ \sup_{s\in S} vol_k(\chi(\tilde Y_s)).$$
(It is more logical  to use the {\it volumes of the maps} $\chi_{|Y_s}:\tilde Y_s\to X$ rather than  of their images but this is not essential at this point.)

[\textbf{B}]  {\it The minmax  $k$-volume} of the pair of  the {\it homotopy classes} of maps $\varsigma $ and $\chi$ 
 denoted $vol_k [\varsigma, \chi]$
 is defined as
 $$\inf_{{\varsigma, \chi}}\sup_{s\in S} Haumes_k(\chi(Y_s)),$$
where the infimum is taken over all pairs of maps $(\varsigma, \chi)$
in  a given homotopy class  $[\varsigma, \chi]$ of (pairs of) maps.

 \vspace {1mm}
 
 [\textbf {C}] {\it The  $k$-Waist} of a  Riemannian Manifold $X$, possibly with a boundary, is    
 $$waist_k(X)=\inf_{{\cal D}_X}\sup_{s\in S} vol_k(\chi(Y_s)),$$
where the infimum  
    (that, probably,  leads to  the same outcome as taking {\sl maximum}) 
   is  taken over all  diagrams ${\cal D}_X=\{X\overset {\chi} \leftarrow \Sigma\overset {\varsigma}\to S\}$ that represent "{\sl homologicaly substantial}" families of  subsets $Y_s= \chi(\tilde Y_s=\varsigma^{-1}(s))\subset X$  that support   $k$-cycles in $X$. \vspace{1mm}
   
{\it What is Homologically Substantial?} A family of subsets $Y_s\subset X$ 
is regarded as homologically substantial if it satisfy some {\it (co)homology condition} that insures  that  the subsets 
$S_{\ni x}\subset S$, $x\in X$, that consist of $s\in S$ such that $Y_s\ni x$,
is non-empty for all (some?) $x\in X$. 

In the setting of smooth manifolds and smooth families, the simplest such condition  is {\it non-vanishing} of the "{\sl algebraic number of points}" in  $S_{\ni x}$ for generic $x\in X$ that make sense if $dim(Y_x)+dim S=dim X$ for these $x$.

More generally, if $dim(Y_s)+dim(S)\geq  dim X$,  then the corresponding condition in the  {\it bordisms homology theory}  asserts       {\it non-vanishing of the cobordism  class of submanifold}   $S_{\ni x}\subset X$ (for cobordisms regarded as homologies  of a point in the bordism homology theory).

\vspace{1mm}

 {\it $\mathbb Z_2$-waists}.  One arrives at a particular  definition, namely, of what we  call   {\it  $\mathbb Z_2$-waist}   if "homologically" refers to homologies with {\it $\mathbb Z_2$-coefficients} (that is the best understood case),  and, accordingly,   the above  $inf$ is taken over all diagrams  $D_X=\{X\overset {\chi} \leftarrow \Sigma\overset {\varsigma}\to S\}$ where   $S$ and  $\Sigma$ are {\it pseudumanifolds  of dimensions $n$ and $n-k$ with boundaries}\footnote{ An "$n$-psedumanifold with a boundary" is understood here as a simpicial  polyhedral space,  where 
all $m$-simplices  for $m<n$ are contained in the boundaries of $n$-simples   and where the boundary $\partial \Sigma\subset \Sigma$ is comprised of the $(n-1)$-simplifies that have  {\it odd} numbers of $n$-simplices adjacent to them.}  (probably,  using only {\it  smooth manifolds}  $S$ and $\Sigma$ in our diagrams  would lead to the same $\mathbb Z_2$-waist) and  where  
 $\varsigma$ and  $\chi$   are continuous maps, such that
 {\it $\chi: \Sigma\to X$ respects the boundaries,\footnote{ $\Sigma $ has non-empty boundary only if $X$ does and $\varsigma:\Sigma\to S$, unlike $\chi: \Sigma\to X$,  {\it does not have} to send $\partial \Sigma\to \partial S$.} i.e.  $\partial \Sigma\to \partial X$}  and   where    
 {\it $\chi$ has non-zero $\mathbb Z_2$} degree that  exemplifies the idea of   "homological  substantionality".\footnote{If one makes the definition of waists with  the volumes of maps $\chi_{|Y_s}:Y_s\to X$ instead of  the volumes  of their images $\chi(Y_s)\subset X$
that would  leads, a priori,  to {\it larger} waists. However, in the $\mathbb Z_2$-case,  the   waists defined with the volumes of images,  probably, equal  the ones, defined via the volumes of maps,    even for non-manifold targets $X$. This is   obvious under mild regularity/genericity  assumptions on the maps   $\varsigma$ and $\chi$, but needs verification in our setting of general continuous maps. }
 
  \vspace {1mm}

One may render this definition more  algebraic by

$\bullet$ admitting an {\it arbitrary} (decent)  topological  space for the role  $\Sigma$ (that is continuousy mapped by $\varsigma$ to an  $m$-dimensional {\it pseudomanifold} $S$);

and 

$\bullet$ replacing the  {\sl "non-zero degree condition"}  by requiring that the fundamental $\mathbb Z_2$-homology  class  $[X]$ of $X$  should  lie in the image  of homology homomorphism $\chi_\ast: H_\ast(\Sigma;\mathbb Z_2)\to H_\ast(X;\mathbb Z_2)$, or equivalently, that  the cohomology homomorphism 
$\chi^\ast:  H^\ast(X;\mathbb Z_2)\to H^\ast(\Sigma;\mathbb Z_2)$, does not vanish on the fundamental {\it cohomology}  class of $X$.

(This    naturally leads  to a definition of the  {\it  the waists of an $n$-dimensional 
homology and cohomology  class $h$}  of   dimension $N\geq n$, where one may 
generalize/refine  further by requiring   non-vanishing  of some cohomology operation applied on $\chi^\ast(h)$, as in \cite {guth-steenrod}.)

{\it  Examples of Homologically Substantial Families.} (i)  
The simplest, yet significant, instances   of such  families of (virtually  $k$-dimensional) subsets in $n$-manifolds $X$ are  the pullbacks $Y_s=\varsigma^{-1}(s)\subset X$, $s\in \mathbb R^{n-k}$, for continuous maps $\varsigma:X\to  \mathbb R^{n-k}$. 
(The actual dimension of some among  these $Y_s$ may be strictly greater than $k$.)

\vspace {1mm}

(ii) Let $S$ be a subset in  the  projectivized  space\footnote{ {\it Projectivization $P(L)$ of a  linear space $L$}  (e.g. of  maps $X\to \mathbb R^m$)
is obtained by  removing zero and dividing $L\setminus 0$   by the action of $\mathbb R^\times$.} 
$P^\infty$ of continuous maps $X\to \mathbb R^m$  and $Y_s=s^{-1}(0)\subset X$, $s\in S$, be the  zero sets of these maps.\footnote {This  $P^\infty$  {\it is} 
an infinite dimensional projective space, unless $X$ is a finite set.} 
Let $P^\infty_{reg}\subset P^\infty$ consist of the maps $X\to \mathbb R^m$ {\it the images of which linearly span all of} $\mathbb R^m$, where, observe, the inclusion  $P^\infty_{reg}\subset P^\infty$ is a {\it homotopy equivalence} in the present  (infinite dimensional) case.

Observe that  there is a natural map, say $G$,   from  $P_{reg}^\infty$  to the Grassmnian $Gr^{\infty}_m$ of $m$-planes in the linear space of functions $X\to \mathbb R$, where  the linear subspace  $G(p)$ in the space of functions $X\to \mathbb R$  for a  (projectifized)   map $p:X\to \mathbb R^m$  
 consists the compositions  $l\circ p:X\to \mathbb R$ for all linear functions $l:\mathbb R^m\to \mathbb R$.
 
\vspace {1mm}

{\it If the cohomology  homomorphism $\mathbb Z_2=H^m(P^\infty;\mathbb Z_2)\to H^m(S; \mathbb Z_2)$ does not vanish, then the $S$-family  $Y_s\subset X$  is $\mathbb Z_2$-homologically substantial, provided $S\subset P^\infty_{reg}$ i.e. if  {\it the images} of $s:X\to \mathbb R^m$, $s\in S$, {\it linearly spans all of $\mathbb R^m$}  for all $s\in S$.}\vspace {1mm}

Indeed,  if $S$ is an $m$-dimensional psedomanifold, for $m=n-k$, then  the number of points in the subset $S_{\ni x}\subset$
for which $Y_s\ni x$, (that is defined under standard genericity conditions)   does not vanish $\mod  2$. In fact,   it is easy to identify  this number with the value of the  Stiefel-Whitney class  of the complementary bundle to the canonical line bundle over the projective space $P^\infty \eqsim P^\infty_{reg})$, where this "complementary bundle"  is unduced from the the canonical  $(\infty-m)$-bundle over the Grassmnian $Gr^{\infty}_m$  by the above map $G:P_{reg}^\infty\to Gr_m^\infty$.
 
\vspace {1mm}

(iii) In the smooth situation,  the above $Y_s=s^{-1}(0)\subset X$ are, generically,  submanifolds of codimension $m$  in $X$ with {\it trivial normal bundles.}

General submanifolds  and families of these  are obtained by mappings $\phi$ from $X$ to the Thom space $T_m$  of a  universal $\mathbb R^m$
bundle $V$ by taking the pullbacks of the zero $\mathbf {0}\subset V\subset T_m$. 

 And if    $\{\phi_s\}:X\to T_m$ is a family of maps  parametrised by the $m$-sphere  $S\ni s$ that
   equals the closure of a  fibre of $V$ in $T_m\supset V$, then the family 
  $ Y_s=\phi_s^{-1}(\mathbf {0})\subset X$ is homologically substantial, if   the  map $S\to T_m$ defined by  $s\mapsto \phi_s(x_0)\in T_m$, for some point $ x_0\in X$, 
  has non-zero intersection index with the zero   $\mathbf {0}\subset V \subset T_m.$ \vspace{1mm}
    
 (iv)   The above, when   applied to maps between the suspensions of our spaces, $\phi^{\wedge_k}: X\wedge S^k\to T_m\wedge S^k$, delivers families of (virtual) submanifods  of dimensions $n-m$, $n=dim (X),$ mapped to $X$ via the projection $X\wedge S^k\to X$, where these families are nomologically substantial under the condition similar to that in   (iii).   
\enskip

(v) There are more general (non-Thom) spaces, $T$, where   ($(n-m)$-volumes of) pullbacks of $m$-codimensional subspaces $T_0\subset T $  are of  some interest.

For instance, since the space of $m$-dimensional $\Pi$-cocycles of  (the singular chain complex of) $X$ (see  \cite {singularities2}) is homotopy equivalent to space of continuous  maps  $X\to K(\Pi,m)$, it  may be worthwhile to look from this perspective at 
 (e.g.  cellular)  spaces  $T$  that represent/approximate   {\it Eilenberg Maclane's}   $K(\Pi,m)$  with  the codimension $m$ skeletons of such $T$ taken of $T_0\subset T$.







\vspace {2mm}

\textbf{  Positivity of Waists.} {\it The $k$-dimensional $\mathbb Z_2$-waists of all Remanning $n$-manifolds $X$
 are strictly positive  for all $k=0,1,2,...,n$:}
 $$\mathbb Z_2\mbox{-}waist_k(X)\geq w_{\mathbb Z_2}= w_{\mathbb Z_2}(X)>0.\leqno  [\geq_{nonsharp}]_{\mathbb Z_2} $$
 
{\it About the Proof.}  The waists defined with all of the above "$k$-volumes" are {\it monotone under inclusions},   $$ waist_k (U_1)\leq waist_k(U_2)\mbox { for all open subsets $U_1\subset U_2\subset X$, }$$
and they also   properly behave  under $\lambda$-Lipschitz maps   $f:X_1\to X_2$ of non-zero degrees,
 $$waist_k(X_2=f(X_1))\leq \lambda^k waist_k(X_1)\mbox { for maps  $f:X_1\to X_2$ with    $deg_{\mathbb Z_2}(f)\neq 0$}. $$

 Thus, $[\geq_{nonsharp}]_{\mathbb Z_2}$   for all $X$   follows from  that  for   the  unit (Euclidean) $n$-ball $B^n$. 
  
 Then the case $X=B^N$  reduces  to  that  of the unit $n$-sphere $S^n$,\footnote {In fact, the sphere   $S^n$   is bi-Lipschitz equivalent to the double of the ball $B^n$. 
 
 Also the ball $B^n=B^n(1)$ admits {\it a radial $k$-volume contracting map  onto} the sphere $ S^n(R)$ of radius $R$, such that $vol_k (S^k(R))=vol_k(B^k(1)$; this allows    {\it a  sharp evaluation}  of certain "regular $k$-waists" of $B^n$, see \cite {waists}, \cite {singularities2}, \cite {guth-steenrod}.  }
while the lower bound   $[\geq_{nonsharp}]_{\mathbb Z_2}$ for the $k$-waist  of $S^n$  defined with  the {\it Hausdorff measure} of  $Y\subset X=S^n$ is proven in   \cite {singularities2} by a reduction to a {\it combinatorial filling  inequality.}  \vspace {1mm}

On the other hands  the {\it sharp values}  of  $\mathbb Z_2$-waists of spheres
are known for the  {\it Minkowski volumes} and, more generally, for all    $\delta$-$Mink_k$, $\delta>0$.  Namely 
$$ \delta\mbox{-} Mink_k\mbox {-}waist_{\mathbb Z_2}(S^n)= \delta\mbox{-} Mink_k(S^k),\leqno {[waist]_{sharp}}, $$
where  $S^k\subset S^n$ is an equatorial $k$-sphere and \vspace {0.5mm}
 
 \hspace {20mm} $  \delta${-}$Mink_k(S^k\subset S^n) =_{def}\frac{vol_n(U_\delta(S^k))}{\delta^{n-k} vol_{n-k}(B^{n-k})}$\vspace {0.5mm}

\hspace {-6mm}  for  the $\delta$-neighbourhood $U_\delta(S^k)\subset S^n$ of this 
sphere  in $ S^n$  and   the unit ball  $B^{n-k}\subset \mathbb R^{n-k}$.

Every   $
\mathbb Z_2$-homologically substantial $S$-family of "$k$-cycles"  $Y_s\subset S^n$   has a member  say $Y_{s_\circ}$ such that \vspace {1mm}
 
  {\it \hspace {10mm} $   \delta${-}$Mink_k(Y_{s_\circ})\geq   \delta${-}$Mink_k(S^k)$ for {\it all} $\delta>0$  simultaneously.}\vspace {1mm}

 In particular, \vspace {1mm}

{\it given an arbitrary continuous map $\varsigma :S^n\to \mathbb R^{n-k}$, there exist a value $s_\circ\in \mathbb R^{n-k}$, such that  the $\delta$-neighbourhoods of the $s_\circ$-fiber  $Y_{s_\circ}=\varsigma^{-1}(s_\circ)\subset S^n$ are  bounded from below by the volumes of such neighbourhoods of the equatorial $k$-subspheres $S^k\subset S^n$,
$$vol_n(U_\delta(Y_{s_\circ}))\geq vol_n(U_\delta(S^k))\mbox { for all } \delta>0.$$}
Consequently the Minkowski volumes of this $Y_{s_\circ}$ is greater than the volume of the equator, 
$$Mink_k(Y_{s_\circ})\geq Mink_k(S^k)\mbox { for this very  } s_\circ\in \mathbb R^{n-k}. \leqno {[\circ]_k}.$$

This is shown in \cite {waists} by a {\it parametric homological localisation argument}. (See section  7 below; also  see section 19 for further remarks, examples  and  conjectures.)

 \vspace {2mm}

\section{Low Bounds on Volume Spectra via Homological Localization.}

 Start with a simple  instance of homological localisation for {\it the $(n-1)$-volumes of zeros} of families {\it of real functions} on Riemannian manifolds. \vspace {1mm}

[\textbf {I}] {\it Spectrum of $E_{vol_{n-1}}$}.  Let $L$ be an $(N+1)$-dimensional linear  space of  functions $s:X\to \mathbb R$ on a  compact $N$-dimensional  Riemannian manifold $X$ and let $U_1,U_2,...U_N\subset X$  be disjoint balls of radii $\rho_N$, such that   
$$\rho^n_N \sim  \circledcirc_n \frac{vol_n(X)}{N}\mbox { for large $N\to \infty$ },$$ 
where   $\circledcirc_n>0$  is a universal    constant   (as in the "packing  section" 2).       
By the {\it Borsuk-Ulam theorem}, there exists a non-zero function $l\in L$ such that the zero set $Y_l =s^{-1}(0) \subset X$ cuts all $U_i$ into equal halves and  the Hausdorff measures (volumes) of the intersections of $Y_l$ with $U_i$  and {\it the isoperimetric inequality} implies the lower bound 
$$vol_{n-1}(Y_l\cap U_i) \geq \beta_{n-1} \rho_N^{n-1}-o(\rho_N^{n-1}),\mbox { } N\to \infty,$$for $$\beta_{n-1}
      =vol_{n-1}(  B^{n-1}(1))= \frac{\pi^{n-1/2}}{\Gamma(\frac{n-1}{2} + 1)}.$$

Therefore,  

{\it the supremum of the the Hausdorff measures  of the zeros $Y_l\subset X$ of non-identically  zero functions $l:X\to \mathbb R$ from an arbitrary $(N+1)$-dimensioanl linear space   $ L$ of functions on $X$  is bounded from below  for large $N\geq N_0=N_0(X)$ by
$$\sup_{l\in L\setminus 0}vol_{n-1}(Y_l)\geq \delta_n N^\frac {1}{n} vol^\frac {n-1}{n}_n(X)$$
for a universal constant $\delta_n>0$.}\vspace {1mm}

[\textbf {I*}] {\it Homological   Generalisation.} This  inequality  remains valid  for all  {\it non-linear} spaces $L$ (of functions on compact Romanian manifolds $X$) that are  invariant under scaling $l\mapsto \lambda l$, $\lambda\in \mathbb R^\times$, provided the projectivisations $(L\setminus 0)/\mathbb R^\times \subset  P^\infty$  of  $L$ in the projective space $P^\infty$ of all continuous functions on $X$ {\it support the (only)  nonzero cohomology class }  from  $H^N(P^\infty;\mathbb Z_2)=\mathbb Z_2.$

\vspace {1mm}

[\textbf {II}] {\it $E_{vol_{n-m}}$-Spectra  of Riemannian manifolds $X$}. Let $L$ be an $(mN+1)$-dimensional  space of continuous maps $l:X\to \mathbb R^m$  that is invariant under scaling $l\to \lambda l$, $\lambda\in \mathbb R^\times$,
and let the projectivized space $S=(L\setminus 0)/\mathbb R^\times \subset P^\infty$ in the projective space $P^\infty$ of  continuous maps $X\to \mathbb R^m$  modulo   scaling support non-zero cohomology class from $H^{mN}(P^\infty;\mathbb Z_2)=\mathbb Z_2$, e.g. $L$ is a liner pace of maps $X\to \mathbb R^m$ of dimension $mN+1$.

{\it Then the zero sets $Y_s=Y_l=l^{-1}(0)\subset X$,  for $s=s(l)\in S\subset P^\infty$  being non-zero maps $l:X\to \mathbb R^m$ mod $\mathbb R^\times $-scaling,
$$\sup_{s\in S}vol_{n-m}(Y_s)\geq \delta_n N^\frac {m}{n} vol^\frac {n-m}{n}_n(X)\leqno {[\ast]_{n-m}}$$
for   large $N\geq N_0=N_0(X)$ and a universal constant $\delta_n>0$, and where  
 $vol_{n-m}$ stands for the Minkowski $(n-m)$-volume. }

\vspace {1mm} $$vol_{n-m}(Y_s\cup U)\leq (n-m)\mbox {-}waist(U)-\varepsilon$$
Then, by the definition of waist, the cohomology restriction homomorphism 
$H^m(P^\infty;\mathbb Z_2)\to  H^m(S^\setminus (U ) ;\mathbb Z_2)$ vanishes for all $\varepsilon>0$

Apply this to    $N$   open subsets $U_i\subset X$, $i=1,...,N,$ and   observe  that the {\it  non-zero} cohomology  class  that comes to   our $S\subset P^\infty$   from  $H^{Nm}(P^\infty;\mathbb Z_2)(=\mathbb Z_2)$  equals  the $\smile$-product of {\it necessarily non-zero} $m$-dimensional classes coming from $H^m(P^\infty;\mathbb Z_2)$.

This, interpreted as the "simultaneous  $\mathbb Z_2$-homological  substantionality" of the families  $ Y_s(i)=Y_s\cup U_i\subset U_i$, for all $i=1,...,N$, shows that

{\it there exists an $s\in S$, such that 
$$vol_{n-m}(Y_s\cup U_i)\geq (n-m)\mbox {-}waist(U)-\varepsilon,$$}
for all $i=1,...,N$,  in agreement with the definition of waists in the previous section.

Finally,  take an efficient packing of $X$ by $N$ balls $U_i$ as in the above [\textbf {I}]  and 
[\textbf {I*}]  and derive the   above  ${[\ast]_{n-m}}$ from the lower bound on waists (see    $[\geq_{nonsharp}]_{\mathbb Z_2}$ in the previous section)  of $U_i$. QED.\vspace {0.7mm}

(The lower bound   $[\geq_{nonsharp}]_{\mathbb Z_2}$ on waists was formulated under   technical, probably redundant,   assumption on $L$ saying that 
the restrictions of spaces  $L$ to $U_i\subset X$ lie in the subspaces $P^\infty_{reg}(U_i) $ (such that  the  images of the "regular"   $l:U_i\to \mathbb R^m$ span $\mathbb R^m$)  of the corresponding projective spaces  $P^\infty= P^\infty(U_i)$ of maps $U_i\to \mathbb R^m$. 
 This assumption can be easily   removed for $vol_{n-m}=Mink_{n-m}$ by a simple  approximation argument applied to $\delta$-$Mink_{n-m}$-waists
and, probably, it also seems   not hard(?) to remove  for $vol_{n-m}=Haumes_{n-m}$ as well.\footnote{   This was stated as a problem   in section 4.2 of \cite {non-linear} but I do not recall if the major source of   complication  was the issue of regularity.}\vspace {2mm}

\hspace {17mm}{\it  \sc Further Results,  Remarks  and Questions. } \vspace {1mm}

\vspace {1mm}

(a) The (projective) space of the above    $Y_s$ can be seen, at least if  $Y_s$ are "regular",      as a  (tiny for $n-k> 1$) part of the space ${\cal C}_k (X;\mathbb Z_2)$ of all $k$-dimensional  $\mathbb Z_2$-cycles in $X$, say for the $n$-ball $X=B^n$,  where the full space  ${\cal C}_k (B^n;\mathbb Z_2)$ of the relative  $k$-cycles  mod 2 is   Eilenberg-MacLane's $K(\mathbb Z_2, n-k)$ by the Dold-Thom-Almgren theorem, see \cite {almgren-homotopy} and  section 2.2 in \cite {singularities2}.

 If $n-k-1$,  then  ${\cal C}_k (B^n;\mathbb Z_2)=P^\infty $ and  
 $H^\ast({\cal C}_k (B^n;\mathbb Z_2);\mathbb Z_2)$  is the polynomial algebra in a single variable of degree one, but if $n-m\geq 2$, then  the  cohomology algebra of  ${\cal C}_k (X;\mathbb Z_2)$  is freely generated by infinitely many  monomials in  {\it Steenrod  squares}   of the generator of $H^{n-k}({\cal C}_k (B^n;\mathbb Z_2))=\mathbb Z_2$.

 Thus, if $n-k>2$, the  cohomology spectrum  of $ E_{vol_k}$  is indexed not by integers (that,  if $n-k=1$, correspond to graded ideals of the polynomial algebra in a single variable), but  by the graded ideals in a more complicated   algebra $H^{n-k}({\cal C}_k (B^n;\mathbb Z_2))=\mathbb Z_2$ with the Steenrod algebra acting on it. \vspace {1mm} 
 
(b) {\it Guth' Theorem.} The asymptotic of this "{\sl Morse-Steenrod spectra"} of the spaces ${\cal C}_k (B^n;\mathbb Z_2)$ of $k$-cykles in the $n$-balls   were  evaluated, up to  a, probably redundant,  lower order term,  by Larry Guth (see  \cite {guth-steenrod}  where a deceptively simple looking corollary of his results is the following 
 
\vspace {2mm}
 
\hspace{5mm}{\sc Polynomial   Bound on the  Spectral Homotopy Dimension 

\hspace{33mm} for the Volume Energy.}\vspace {1mm}
 
  {\sl The homotopy dimensions (heights) of the sublevels $\Psi_e=E^{-1}(-\infty.e] \subset \Psi$  of the $vol_k$-energy $E=E_{vol_k}:\Psi\to \mathbb R_+$ 
 on the space $\Psi=   {\cal C}_k (X;\mathbb Z_2)$ of $k$-dimensional rectifiable  $\mathbb Z_2$-cycles in a compact  Riemannian manifold $X$
 satisfies
$$homdim(\Psi_e)\leq ce^\delta$$
where the constant  $c=c(X)$ depends on the geometry of $X$ while $\delta=\delta(n)$ depends only on the dimension $n$ of $X$.}\vspace {1mm}

This, reformulated as a lower spectral bound on $E_{vol_k}$, reads.\vspace{1mm}

 {\sl Let $X$ be a compact Riemannian manifold and   let $ Y_s\subset X$, $s\in S\subset   {\cal C}_k (X;\mathbb Z_2)$,
be an $S$-family of $k$-cycles (one may think of these as $k$-pseudo-submanifolds in $X$) that is {\it \textbf {not} contractible to any $N$-dimensional subset} in  $\Psi={\cal C}_k (X;\mathbb Z_2)$.
Then 
$$\sup_{s\in S} vol_k(Y_s)\geq \varepsilon\cdot  N^\alpha,$$
where $\varepsilon =\varepsilon(X)>0$ depends on the geometry of $X$, while  $\alpha=\alpha(n) (=\delta^{-1})$, $n=dim(X)$, is a universal positive  constant.} \vspace {1mm}

In fact, Guth's   results yield a   nearly sharp bound 
$$\sup_{s\in S} vol_k(Y_s)\geq  \varepsilon(X,\alpha)  \cdot  N^\alpha,$$
for all $\alpha<\frac {1}{k+1}$, where, conjecturally, this must be also true for $\alpha=\frac {1}{k+1}$.


  \vspace{2mm}

(c) {\sl  Is there a direct simple proof of the above inequaliy  with some, let it be non-sharp,
$\alpha$  that would bypass  fine   analysis (due to Guth) of the Morse-Steenrod cohomology spectrum of $E_{vol_k}$ on  the space of cycles?

 Does a  polynomial lower bound  hold for  ($k$-volumes of)  $\mathbb Z_p$- and for $\mathbb Z$-cycles?} 

Apparently, compactness  of spaces of  (quasi)minimal subvarieties in $X$   implies      discreteness of  homological volumes spectra  via the  Almgren-Morse theory,   but  this does not seem(?) to deliver   even logarithmic   lower spectral bounds due to the  absence(?)  of corresponding  bounds   on Algren-Morse indices of
 minimal subvarities  in terms of their volumes.\vspace {1mm}

\vspace {1mm}

(d) The     lower bounds  for the $k$-volumes of zeros  of families  maps $X\to \mathbb R^m$  (see [\textbf {II}] above) can be, probably,  generalized in the spirit of Guth' results, at least  in  the $\mathbb Z_2$-setting,    to spaces of  maps  $\psi $ from $X$ to the total spaces   of $\mathbb R^m$-bundles $V$ 
and to the Thom spaces of such bundles where $E_{vol_k}(\psi)=_{def}vol_k(\psi^{-1}(\mathbf {0}))$  for the zero sections $\mathbf {0}\subset V$ of these bundles\footnote {A  similar effect can be achieved by replacing   $\mathbb R^\times$-scaling   by the action of the full linear group $GL_m(\mathbb R)$ on $\mathbb R^m$ and working with the  equivaruant  cohomology of the space of maps $X\to  \mathbb R^m$  with  the corresponding action of  $GL_m(\mathbb R)$ on it. }  where the Steenrod squares should be replaced by the bordism  cohomology operations   that are in the $\mathbb Z_2$-case amount  to  taking Stefel-Whitney classes.

 Here, as well as  for spaces of maps $\psi$ from $X$ to more  general spaces $T$ where $E_{vol_k}(\psi )=vol_k(\psi^{-1}(T_0))$ for a given $T_0\subset V$, one needs somehow to factor away the homotopy classes of maps $\psi$  with  $E_{vol_k}(\psi )=0$ 
 (compare \cite {manin}). 
 
Also it may be interesting to augment the $k$-volume by other (integral) invariants of $Y=\psi^{-1}(T_0)$, where the natural candidates in the case of $k$-dimensional (mildly singular)  $Y=\psi^{-1}(\mathbf {0})$ would be curvature integrals expressing  the $k$-volumes of the tangential lifts of  these $Y\subset X$ to the Grassmann spaces  $Gr_k(X)$, $k=dim Y$, of target $k$-planes in $X$  (compare section 3 in\cite {proteins}).
  
  \vspace {1mm}

(e)   One can define    spaces of subsets $Y\subset X$ that support   " $k$-cycles (or rather {\it cocycles}), where these $Y$ do not have to be regular in any way, e.g. rectifiable as in Guth' theorem, or even geometrically $k$-dimensional.  But Guth's parametric  homological localization along with the bounds on waists  from  the previous section yield the same lower bounds on the volume spectra on these space as in the rectifiable case.

\vspace {1mm}

\vspace {1mm}

\section {Variable spaces, Homotopy Spectra  in Families and   Parametric  Homological  Localisation. 
  }   \vspace {1mm}

Topological spaces  $\Psi$  with (energy) functions  $E$ on them often come in families. 

In fact, the proofs  of the sharp lower bound on the Minkowski waists of spheres (see  section 6)  and  of Guth' lower bounds on  the full  $vol_k$-spectra (see the previous section)  depend on localizing {\it not} to (smallish) {\it fixed} disjoint subsets  $U_i \subset X$ but to variable or  "{\it parametric}" ones  that may change/move along with  the subsets  $Y_s$, $s\in S$, in them. 

 In general, families $\{\Psi_q\}$ are  constituted by  "fibres" of  continuous maps 
$F$ from a space $\Psi=\Psi_Q$ to  $Q$ where the "fibers" $\Psi_q=F^{-1}(q)\subset \Psi$, $q\in Q$, serve as  the members of these families and  where the energies $E_q$ on $\Psi_q$ are obtained by restricting functions  $E$ from $\Psi$ to $\Psi_q\subset \Psi$.\footnote{ Also one may have functions with the range also depending on $q$, say  $ a_q:  \Psi_q\to R_q$  and  one may generalise further by defining families as some  (topological)  sheaves over Grothendieck sites.}\vspace {1mm}

Homotopy spectra in this situation are defined with continuous families of spaces $S_q$  that are "fiberes"  of continuous maps $S\to Q$ and  where
the relevant maps $\phi:S\to \Psi$ send   $ S_q\to \Psi_q$  for all $q\in Q$ with these  maps denoted $\phi_q=\phi_{|S_q}$.

Then the energy of  the fibered  homotopy class $[\phi]_Q$ of such a fiber preserving map $\phi$  is defined as earlier as 
   $$E[\phi]_Q= \inf_{\phi\in [\phi]_Q}\sup_{s\in S} E\circ \phi(s)\leq \sup_{q\in Q} E_q[\phi_q],$$
where the latter inequality is, in fact,  an equality in many cases.\vspace {1mm}

{\it Example 1:  $k$-Cycles in Moving  Subsets}.  Let $U_q$ be  a $Q$-family of open subsets in a Riemannian manifold $X$.  
An  instance of this is the family of the   $\rho$-balls $U_x=U_x(\rho)\in X$ for a given  $\rho\geq 0$ where $X$ itself plays the role of $Q$.

Define  $\Psi= \Psi_Q$   as     the space    $C_k\{U_q;\Pi\}_{q\in Q}$  of $k$-dimensional  
$\Pi$-cycles{\footnote{ "$k$-Cycle" in $U\subset X$    means  a {\it relative} $k$-cycle in $(U,\partial U)$, that is  a $k$-chain with boundary contained in the boundary of $U$.  Alternatively, if $U$ is an  open non compact subset,  "$k$-cycle" means   an {\it infinite} $k$-cycle, i.e. with (a priori) {\it non-compact support}.} 
 $c=c_q$ in      $U_q$ for  all $q\in Q$, that is $\Psi=\Psi_Q$ equals the space of pairs ($q\in Q, c_q\in C_k(U_q;\Pi) )$, where, as earlier,  $\Pi$ is an Abelian (coefficient) group with a norm-like function; then we  take     $E(c)=E_q(c_q)=vol_k(c)$ for the energy.\vspace {1mm}

{\it Example 2: Cycles in Spaces  Mapped to an $X$.} Here, instead of subsets in $X$ we take locally diffeomorphic maps  $y$ from a fixed Riemannian  manifold $U$ into $X$  and take the Cartesian product $C_k(U;\Pi)\times Q$ for $\Psi=\Psi_Q$. \vspace {1mm}

 {\it Example  1+2: Maps with variable domains.}  One may   deal  families of spaces $U_q$  (e.g. "fibers" $U_q=\psi^{-1}(q)$ of a  map between smooth manifolds $\psi:Z\to Q$) along with maps $y_q:U_q\to X$.
\vspace {1mm}

{\it On   Reduction 1 $\Rightarrow$ 2.} There  are cases,  where the spaces $\Psi_Q=C_k\{U_q;\Pi\}_{q\in Q}$ of cycles in moving  subsets $U_q\subset X$ topologically
split: $$\Psi_Q=C_k(B;\Pi)\times Q, \mbox { for a fixed manifold $U$}.$$


A simple, yet representative, example is where  $Q=X$ for  the $m$-torus, $X=\mathbb T^m=\mathbb R^m/\mathbb Z^m$, where   $B=U_0$ is an open  subset in $\mathbb T^m $ and where $Y=X=\mathbb T^m$ equals the space of translates  $U_0\mapsto U_0+x$, $x\in \mathbb T^m$.

For instance, if $U_0$ is a ball of radius $\varepsilon\leq 1/2$, then it can  be identified with the Euclidean $\varepsilon$-ball 
$B=B(\varepsilon)\subset \mathbb R^m$.

Similar splitting is also possible  for {\it parallelizable} manifolds $X$ with injectivity radii $>\varepsilon$  where  moving $\varepsilon$-balls $U_x\subset X$ are obtained via the exponential maps $exp_q:  T_q=\mathbb R^m\to X$ from a fixed ball $B=B(\varepsilon)\subset \mathbb R^m$.

In general, if $X$ is {\it non-parallelizable},  one may take  the space  of the tangent orthonormal frames in $X$ for $Q$, where,  the  product space   $C_k(B;\Pi)\times Q$,  where $B=B(\varepsilon)\subset \mathbb R^m$, makes a principle 
$O(m)$- fibration, $m=dim(X)$, over the space  $C_k\{U_x(\varepsilon);\Pi\}_{x\in X}$ of cycles of moving $\varepsilon$-balls $U_x(\varepsilon)\subset X$.\footnote{ Vanishing of Stiefel-Whitney classes seems to  suffice for   (homological) splitting of this vibration  in the case  $\Pi=Z_2$  as  in section 6.3 of my article {\sl  "Isoperimetry of Waists...")} in GAFA }



\vspace {2mm}

 Waists of   "variable metric spaces"  needed for  homological localization  of  the  spectra of the volume energies  on {\it fixed spaces}   can be defined as follows.\vspace {1mm}


 Let ${\cal X}=\{X_q\}_{q\in Q }$ be a family of metric spaces   seen as the   fibres, i.e. the pullbacks of  points,  of a continuous map $ \varpi: {\cal X}\to Q$ and consider   
  $S$-families  of subsets in $X_q $  that are    $Y_s\subset X_{q(s)}$ defined with some   maps $S\to Q $  for $s\mapsto 
q=q(s)$.  
  
  The $k$-waist of such a family is defined as 
  $$waist_k\{X_q\}=\inf \sup_{s\in S}vol_k(Y_s)$$
 where $vol_k$ is one of the "volumes"  from  section 6, e.g.  Hausdorff's $k$-measure and where "inf" is taken over all {\sl "homologically substantial"} families $Y_s$.
 
 In order to define the latter and to keep the geometric picture in mind, we 
 
 $\bullet $ fix a section $Q_0\subset \cal X$,
 that is a continuous  family\footnote{ In fact, one only needs a distinguished "horizontal (co)homology class" in $\cal X$.} of points $x(q)\in X_q$; 
 
  $\bullet$  assume that 
 $X$ is an $n$-manifold and that the family $Y_s$ is given by fibres of a map $\Sigma=\cup_{s\in S} Y_s\to S$  for $\Sigma$ being an $n$-pseudomanifold   of dimensions $n$.  

Then the family $Y_s\subset X_{q(s)}$ comes via a map $\Sigma\to \cal  X$ and "homologically substantial" is  understood as {\it non-vanishing} of  the intersection index of  $Q_0 \subset \cal X$ with $\Sigma$ mapped to $\cal X$.

For instance, if $ \varpi: {\cal X}\to Q$ is a fibration with {\it contractible} fibres $X_s$, the section $Q_0\subset \cal X$ exits and homotopically unique that makes  our   index non-ambiguously  defined.

Notice that $waist_k\{X_q\}$ may be strictly smaller than $\inf_{q\in Q} waist_k(X_q)$; yet,  the argumnt(s) used for individual $X$ show that the waists of {\it compact}   families of (connected) Riemannian manifolds are {\it strictly positive}.\vspace {1mm}

{\it Example: "Ameba" Penetrating a  
 Membrane.}  Let   a domain  $U_t\subset \mathbb R^3$, $0\leq t \leq 1$, be composed   of  a pair of disjoint  balls of radii $t$ and $1-t$ joint by a $\delta$-thin tube. The $2$-waist of $U_t$ is at least  $\pi/4$, that is the waist of the ball of radius $t/2$ for all $t\in [0,1]$.  But the waist of the "variable domain" $U_t$  equals 
 the area  of the section of the tube that is   $\pi \delta^{2}$.

\section {Restriction and Stabilisation  of Packing Spaces.} There is no significant relations between {\it individual} packings of  manifolds  and their submanifolds, but such 
relations do exist for  {\it spaces of packings}.

For instance, let $X_0\subset X$ be a closed  $n_0$-dimensional  submanifold in an $n$-dimensional, Riemannian manifold $X$ and  let
 $$\frown_{ I}:  H_\ast  (X^I)\to   H_{\ast-\ast'}  (X_0^I), \mbox {  $ \ast'=   N(n-n_0)$,  { }$N=card (I),$}$$
 be the homomorphism corresponding to the  (generic) intersections of cycles in the Cartesian power  $X^I$ with the submanifold 
 $X_0^I\subset X^I$, where the  homology groups are understood  with $\mathbb Z_2=\mathbb Z/2\mathbb Z$-coefficients, since  we do not assume that $X_0$ is orientable.
Then \vspace {1mm}

{\it the homological packing energies $E_\ast$   of $X$ and $X_0$} (obviously)  {\it satisfy
$$E_\ast(\frown_I(h)) \leq E_\ast(h),$$}
for all homology classes   $h\in H_\ast(X^I)=H_\ast(X)^{\otimes I}$,  where packings of $X_0$ are understood with respect to the metric, i.e. distance function,  induced 
from $X\supset X_0$.

{\it Corollary}. Let $P(X;I,r) \subset X^I$ be the space of $I$-packings of $X$ by balls of radii $r$ and let $S\subset P(X;I,r)$ be a $K$-cycle, $K=N(n-n_0)$ that has a non-zero intersection index with 
$X_0^I\subset X^I$.

{\it Then $X_0$ admits a packing by $N$-balls of radius $r$.}\vspace {1mm}

Next, let us invert the intersection homomorphism $\frown_I$ in a presence of a {\it projection} also called {\it retraction}  $p: X\to X_0$ of   $X$ to $X_0$,  i.e. where  $p$ fixes  $X_0$.

  If $p$ is a fibration or, more generally it is a  generic smooth $p$, then the 
pullbacks $Y_i= p_{-1}(x_i) X$ are $k$-cycles, $k-n-n_0$, that transversally meet $X_0$ at the points $x_i\in X_0$.   It follows, that  the Cartesian product $N(n-n_0)$-cycle 
$$S=\bigtimes_{i\in I}Y_i \subset X^I$$
has non-zero intersection index with $X_0^I\subset X_I$.

 And -- now geometry enters -- if $p$ is a  (non-strictly) {\it distance decreasing} map, then, obviously, this $S$ is positioned in the space  $P(X;I,r) \subset X^I$ of $I$-packing of 
  $X$ by  $r$-balls  $U_r(x_i)$ and  multiplication of  cycles in $ C\subset P(X_0;I,r)$ with $S$, that is $C\mapsto C\times S$, 
    {\it imbeds} $$H_\ast(P(X_0;I,r))  \underset {\times S}\to H_{\ast+Nk}(P(X;I,r)),$$
    such that the composed map 
    $$H_\ast(P(X_0;I,r ))  \underset {\times S}\to H_{\ast+Nk} ((P(X;I,r)) \underset  {\frown_I}\to H_\ast(P(X_0;I,r))$$ equals the identity.

 Thus, \vspace {1mm}
 
\hspace {-4mm} {\it the homology packing spectrum of $X$ fully determines such a 
 spectrum of $X_0$.}\footnote{If $p$ is {\it a homotopy retraction}, e.g. $p:X\to X_0$     
 is a  vector bundle, then the  above homomorphisms come from the the Thom {\it isomorphisms} between corresponding spaces. Yet, there is  more to packings of $X$ than what comes from $X_0$, it   already seen in the case where $X$ is a ball and  $X_0=\{0\}$.  
 
 However,  ball packings of $X_0$ in this case  properly reflect  properties of packing of $X$ by $r$-neighbourhoods of $k$-cycles  $Y_i\subset X$ the  intersection indices of which with $X_0$ equal one. }
 
\vspace {1mm}
 
{\it Example.}  Let $\underline X$ be a compact  manifold of negative curvature and $X_0\subset \underline X$ be a closed geodesic. Then the above applies to the covering
 $X =X(X_0)$  of $ \underline X$ with the cyclic  fundamental  group generated by the homotopy class of $X_0$.
 
 Therefore, \vspace{1mm}
 
 {\it the lengths of all the closed geodesics in $\underline X$  are determined by  
 
 the homotopy packing spectrum of $\underline X$.} \vspace{1mm}

{\it Questions.}  How much of the geometry of  minimal varieties $V$  in $X$, that are critical points of the volume energies, can be seen in terms of the above   families of  packings of $X$ by balls (or by non-round subsets as in  \cite {waists}) "moving transversally to" $V$?

\vspace{1mm}

(Minimal subvarieties  $V$ can be approximated by  sets of centres  of  small $\delta$-balls densely packing these $V$; this suggests looking at spaces of packings of 
$c\delta$-neihgbourhoods  $U_{c\delta}(V)\subset X$ by $(1-\varepsilon)\delta$-balls   for  some  $c>1$.\vspace{1mm}

 {\it Symplectic  Remark.} The above  relation between, say,  individual  packings of an $X_0\subset X$ by $N$ balls and 
$N(n-n_0)$-dimensional families  of packings of $X$ by balls "moving   transversally to $X_0$" is reminiscent of {\it hyperbolic stabilisation of Morse functions} as it used in the study of   {\it generating functions} in the symplectic geometry, see \cite \cite {viterbo},  {eliashberg-gromov} and references therein.  

Is there something more profound here than just a  simple minded similarity?
 
 

\vspace {1mm}

\section {Homotopy Height, Cell Numbers and Homology. }

The homotopy    spectral values $r\in \mathbb R$  of $E(\psi)$  are  "named" after (indexed by) the homotopy classes $[\phi]$ of maps $\phi: S\to \Psi$,  where $r=r_{[\phi]}$ is, by definition,  the minimal $r$   such that $[\phi]$ comes from a map $S\to \Psi_r\subset \Psi$ for 
$\Psi_r= E^{-1}(-\infty,r]$.     In fact, such a "name" depends only on the partially ordered set, cal it ${\cal H}_\gtrless(\Psi)$, that is {\it the maximal partially ordered reduction} of   ${\cal H}_\circ(\Psi)$ defined as follows.

Write $[\phi_1] \prec [\phi_2]$  if there  is a morphism $\psi_{12}:  [\phi_1] \to  [\phi_2]$ in  ${\cal H}_\circ(\Psi)$
and turn this  into a partial order by identifying objects, say      $[\phi]$  and $[\phi']$, whenever  $[\phi] \prec [\phi']$ as well as  $[\phi'] \prec [\phi]$. \vspace {1mm}

{\it Perfect Example.} 
If  $X$ is (homotopy equivalent to) the real  projective space $P^\infty$ then the partially ordered set  ${\cal H}_\gtrless(\Psi)$ is isomorphic to the set of nonnegative  integers $\mathbb Z_+=\{0, 1,2, 3,....\}$. This is why spectral  (eigen) values  are indexed by integers in the classical case.

\vspace {1mm}
In general the set ${\cal H}_\gtrless(\Psi)$ may have undesirable(?) "twists". For instance, if $\Psi$ is homotopy equivalent to the circle, then ${\cal H}_\gtrless(\Psi)$ is isomorphic to set   $\mathbb Z_+$ with the {\it  division order}, where
$m\succ n$ signifies that $m$ divides $n$. (Thus,  $1$ is the maximal element  here and $0$ as the minimal one.) 

Similarly, one can determine  ${\cal H}_\gtrless(\Psi)$  for general Eilenberg-MacLane spaces $\Psi=K(\Pi, n)$.
This seems transparent 
for Abelian  groups $\Pi$. But
 if  a space $\Psi$, not necessarily a $K(\Pi, 1)$, has a  non-Abelian fundamental group $\Pi=\pi_1(\Psi)$, such as the above space $\Psi_N(X)$ of subsets  $\psi\subset X$ of cardinality $N$, 
  then keeping  track of 
the conjugacy classes of subgroups  $\Pi'\subset \Pi$ and 
maps $\phi: S\to \Psi$  that send $\pi_1(S)$ to these $\Pi'$ becomes more difficult.
\vspace {1.5mm}

If one wishes (simple mindedly?) to remain with  integer valued spectral  values, one has to pass to some numerical  invariant that  that takes values in a quotient of  ${\cal H}_\gtrless$ isomorphic to $\mathbb Z_+$, e.g. as follows.  

{ \it  Homotopy Height. } Define  {\it the homotopy (dimension) height } of a  homotopy class $[\phi]$  of continuous map $\phi:S\to \Psi$  as the minimal integer $n$ such that the $[\phi]$ factors  as  $S \to K\to \Psi$, where $K$ is a cell complex of dimension (at most)  $n$. \vspace {1mm}

{\it "Stratification" of  Homotopy Cohomotopy Spectra by Hight.} \vspace {1mm}
This  "hight"  or  a similar hight-like function defines a  partition of the homotopy spectrum into the subsets, call them $Hei_n\subset \mathbb R$, $n=0,1,2,...$, of the values of the energy $E[\phi]\in \mathbb R$ on  the homotopy classes  $[\phi]$ with homotopy heights $n$, where either  the {\sl supremum} or  the {\sl infimum} of  the numbers $r\in Hei_n$ may serve as the {\sl "$n$-th HH-eigenvalue of $\psi$".}\vspace {1mm}

One  also  may  "stratify" cohomotopy spectra by replacing "{\sl contractibility condition} of  maps  $\psi_{|\Psi_r}:\Psi_r\to T$ by  $\psi_{|\Psi_r}\leq n$.\vspace {1mm}

In the classical case of  $\Psi=P^\infty$  any  such  "stratification"  of  homotopy  "eigenvalues"  lead the usual indexing of the spectrum.  
where, besides the homotopy hight, among other   hight-like invariant invariants   we indicate the following. \vspace {1mm}

{\it Example 1: Total Cell Number.} Define $N_{cell}[\phi]$  as the minimal $N$ such that $[\phi]$ factors  as  $S \to D\to \Psi$, where $D$ is a cell complex   with (at most)  $N$ cells in it. \vspace {0.7mm}

What are, roughly, the total   cell numbers of the  classifying maps from  packing spaces of an $X$ by $N$ balls to the classifying space $B\mathbb S_N$?

\vspace {0.7mm} What are these numbers for the maps between classifying spaces of "classical" finite groups $G$ corresponding to standard  injective homomorphisms $G_1\to G_2$?

\vspace {1mm}

{\it Example 2: Homology Rank.}  Define $rank_{H_\ast}[\phi]$ as the maximum over all fields $\mathbb F$ of the the $\mathbb F$-ranks of the induced homology homomorphisms 

\hspace {-6mm} $[\phi]_\ast : H_\ast (S;\mathbb F)\to H_\ast (\Psi;\mathbb F)$. \vspace {1mm}

{\it On Essentiality of Homology.} There are other  prominent spaces,  $X$, besides  the infinite dimensional projective spaces  $X=P^\infty$, and  energy functions on them,   such as

{\it  spaces $\Psi$ of loops $\psi: S^1\to X$ in simply connected Riemannian  manifolds $X$  

with  length($\psi$) taken for $E(\psi)$}\footnote{This instance of essentiality of the homotopy heights is explained in my article  {\sl   Homotopical Effects  of  Dilatation}, while      the full range of this property among "natural" spaces $\Psi$  of maps  $\psi$ between Riemannian manifolds  and  energies  $E(\psi)$  remains unknown.}, 

\hspace {-6mm} where the cell numbers and the homology ranks spectra  for   $  E(\psi)=lenght (\psi)$ are "essentially determined" by  the homotopy height.
(This is why  the homotopy height  was singled out under the name of    "essential dimension"  in my paper {\sl Dimension, Non-liners Spectra and Width.}) 

However, the homology carries   significantly more information  than the homotopy hight for the $k$-volume function on the  {\it spaces of $k$-cycles of codimensions} $\geq 2$ as it was revealed by Larry Guth in his  paper  {\sl Minimax problems related to cup powers and Steenrod squares}.  \vspace{1mm}

{\it On Height and the Cell Numbers of Cartesian Products}. If the homotopy heights  and/or cell numbers of  maps $\phi_i:S_i\to \Psi_i$, $i=1,...,k $, can be expressed in terms of the corresponding homology homomorphisms over some filed $\mathbb F$ independent of $i$, then, according to {\it K\"unneth formula,} the homotopy hight  of the Cartesian product of maps,
$$\phi_1\times...\times \phi_k:S_1\times...\times S_k\to \Psi_1\times...\times \Psi_k,$$
is additive
$$height[\phi_1\times...\times \phi_k]=height[\phi_1]+...+ height[\phi_k]$$
and the cell number is multiplicative
$$N_{cell}[\phi_1\times...\times \phi_k]=N_{cell}[\phi_1]\times...\times N_{cell}[\phi_k].$$

{\it What are other  cases where these relation remain valid?}

Specifically, we want to know what happens in this regard to the following classes of maps:

(a)   {\it maps between spheres} $\phi_i:S^{m_i+n_i}\to S^{m_i} $,  

(b)  {\it maps  between locally symmetric spaces}, e.g. compact manifolds of {\it constant negative curvatures},

(c) {\it high Cartesian powers}  $\phi^{\times N}: S^{\times N}\to A^{\times N}$  of a single map $\phi: S\to \Psi$. 

When do, for instance,   the limits
$$\lim_{N\to \infty}\frac{height[\phi^{\times N}]}{N} \mbox {  and   } \lim_{N\to \infty}\frac{\log N_{cell}[\phi^{\times N}]}{N}$$ 
{\it not vanish}? (These limits exist, since the the hight and the logarithm of the cell number are sub-additive under Cartesian product of maps.) \vspace {1mm}

Probably, the   general question for  {\it "rational homotopy classes $[...]_{\mathbb Q}$"}   (instead of  "full" homotopy classes" $[...]=[...]_{\mathbb Z}$) of maps into  { \it simply connected} spaces $\Psi_i$  is   easily solvable with   {\it Sullivan's minimal models}. 
 
 Also, the question may be more manageable  for   {\it  homotopy classes $\mod p$.}\vspace {1mm}

{\it  Multidimensional Spectra Revisited.} Let $h=h^T$ be a cohomotopy class of $\Psi$, that is a homotopy class of maps $\Psi\to T$,  and let  $ \upsilon$ be  a function on homotopy classes of maps  $U\to T$ for open subsets $U\subset \Psi$, where the above height-like functions, such as the {\sl homology rank}, are relevant examples of such a  $ \upsilon$.

Then the values of $ \upsilon$ on $h$ restricted to open  subsets $U\subset \Psi$ define a numerical (set) 
function, $U\mapsto   \upsilon(h_{|U})$ 
and  every  continuous map ${\cal E}:\Psi\to Z$  pushes down  this function  to open subsets in  $ X$.

For instance, if $  \upsilon=0,1$ depending on whether a map $U\to T$ is contractible or not  and if $Z=\mathbb R^l$,
then this function on the "negative octants"  $\{x_1<e_1,...,x_k<e_k,...,,x_l <e_l\}$  in $\mathbb  R^l $ carries the same message as $\Sigma_h$ from  section 4.

\section{Graded Ranks, Poincare Polynomials, Ideal Valued Measures  and Spectral $\smile$Inequality.} 
The images as well as kernels  of (co)homology homomorphisms that are induced by continuous maps are {\it graded}  Abelian groups and  their ranks are properly  represented not  by individual numbers but by {\it Poincar\'e polynomials}.   

Thus, sublevel  $\Psi_r=E^{-1}(-\infty,r] \subset \Psi$ of   energy functions  $E(\psi)$ are characterised by the {\it  polynomials}  {\upshape  Poincar\'e$_r(t;\mathbb F)$}  of the the inclusion homomorphisms $\phi_i(r): H_i(\Psi_r; \mathbb F) \to H_i(\Psi; \mathbb F)$,
that are
$$\mbox {{\upshape Poincar\'e}$_r$=Poincar\'e$_r(t;
\mathbb F)=\sum_{i=0,1,2,...} t^i rank_\mathbb F\phi_i(r)$}.$$
Accordingly, the homology spectra, that are the sets of those $r\in \mathbb R$ where the ranks of $\phi_\ast(r)$ change,  are indexed by such   polynomials  with positive  integer coefficient. (The semiring structure  on the   set of such  polynomials coarsely agrees with basic topological/geometric constructions, such as taking $E(\psi) =E(\psi_1)+E(\psi_2)$ on $\Psi=\Psi_1\times \Psi_2$.)\vspace {1mm}

The   set function $D\mapsto$ {\upshape Poincar\'e}$_D$  that assigns these  Poincar\'e polynomials   to  subsets $D\subset \Psi$, (e.g.  
$D=\Psi_r$) has some measure-like properties  that become more pronounced for the set function
$$\Psi\supset D\mapsto\mu(D)=\mu^\ast(D;\Pi)= \mathbf 0^{\setminus \ast}(D; \Pi)\subset H^\ast= H^\ast(\Psi;\Pi),$$
where  $\Pi$ is an Abelian (homology coefficient)  group, e.g. a field $\mathbb F$, and  $\mathbf 0^{\setminus \ast}(D; \Pi)$ is the {\it kernel} of the cohomology restriction homomorphism for the complement  $\Psi\setminus D\subset \Psi$, 
$$H^\ast(\Psi;\Pi)\to H^\ast(\Psi\setminus D;\Pi).$$

Since the  cohomology classes   $h\in\mathbf 0^{\setminus \ast}(D; \Pi)\subset H^\ast=H^\ast(\Psi;\Pi)$ are representable by cochains with the support in $D$,{ \footnote{ This property   suggests
an extension of  $\mu$ to
 multi-sheated {\it domains $D$ over $\Psi$} where   $D$  go to $\Psi$ by    non-injective, e.g. locally homeomorphic  finite to one,  maps  $ D \to A $.}\vspace {1mm}

  {\it the set function  
  $$\mu^\ast: \{subsets \subset  \Psi\}  \to  \{subgroups \subset   H^\ast\}$$ 
  is additive  for the sum-of-subsets  in   $H^\ast$   and  super-multiplicative\footnote{ This, similarly to   {\it Shannon's subadditivity inequality}, implies the existence of {\sl "thermodynamic limits"} of {\it Morse Entropies}, see \cite {bertelson-gromov}. } for the  the $\smile$-product of ideals in the case  where $\Pi$ is a commutative ring}:

 \vspace {1mm}

$$\mu^\ast(D_1\cup D_2)=\mu^\ast(D_i)\bigplus \mu^\ast(D_2)\leqno [\cup\bigplus]$$
for {\it disjoint} open subsets $D_1$ and $D_2$  in $\Psi$, and
$$\mu^\ast(D_1\cap D_2)\supset \mu^\ast(D_1)\smile\mu^\ast(D_2) \leqno [ \cap\smile]$$
for all  open $D_1,D_2\subset \Psi$.\footnote {  See section 4 of  my article {\sl Singularities, Expanders and Topology of Maps. Part 2.} for further properties and  applications of these "measures" .}\vspace {1mm}


 \vspace {1mm}
The relation  $[ \cap\smile] $, applied to  $D_{r,i}= E_i^{-1}(r,\infty)\subset \Psi$ can be equivalently expressed in terms of cohomomoly spectra as follows. \vspace {1mm}

{\it Spectral  $[\min\smile]$-Inequality.}\footnote{ This inequality implies the existence of {\it Hermann Weyl limits} of   energies of cup-powers in infinite dimensional projective (and similar) spaces, see \cite {non-linear}.} Let   
$E_1,...,E_i,..,E_N:\Psi\to  \mathbb R$ be continuous functions/energies and  let $E_{min}:\Psi\to \mathbb R$ be the minimum of these,
$$E_{min}(\psi)=\min_{i=1,...,N} E_i(\psi), \mbox { } \psi\in \Psi.$$
Let $h_i\in H^{k_i} (\Psi;\Pi)$ be cohomology classes, where $\Pi$ is a commutative ring, and let 
$$h_\smile\in H^{\sum_i k_i} (\Psi;\Pi)$$
be the $\smile$-product of these classes,
$$h_\smile=h_1\smile ...\smile h_i\smile...\smile h_N.$$
Then 
$$E^\ast_{min} (h_\smile)\geq \min_{1=1,...,N} E^\ast_i(h_i).\leqno [\min \smile]$$
Consequently, the value of the  {"\sl total energy"} $$E_\Sigma=\sum_{i=1,...,N}E_i:\Psi\to  \mathbb R$$  on this cohomology 
class $h_\smile\in H^\ast(\Psi;\Pi )$ is bounded from below by
$$E^\ast_\Sigma(h_\smile)\geq \sum_{i=1,...,N} E^\ast_i(h_i).$$
(This has been  already used in the homological localisation of the volume energy in section 7.)
kk
\vspace {1mm}

{\it On Multidimensional  Homotopy Spectra.} These spectra, as defined in sections 4, 10,11  represent the values of the pushforward of the "measure" $\mu^\ast$
under maps ${\cal E} :\Psi\to \mathbb R^l$   on special subsets $\Delta\subset \mathbb R^l$;  namely, on complements to  $\bigtimes_{k=1,...,l} (-\infty,e_k]\subset \mathbb R^l$ and the spectral information is encoded by $\mu^\ast({\cal E}^{-1}(\Delta))\subset H^\ast(\Psi).$

On may generalise this by enlarging the domain of $\mu^\ast$, say, by evaluating   $\mu^\ast({\cal E}^{-1}(\Delta))$   for some class of simple subsets    $\Delta$ in  $\mathbb R^l$, e.g. convex sets and/or their complements. 

\vspace {1mm}

{\it On $\wedge$-Product.} The  (obvious) proof of   $[ \cap\smile]$  (and of $[\min \smile]$)  relies
 on locality of the $\smile$-product that, in homotopy theoretic terms, amounts to factorisation of $\smile$ via $\wedge$ that is
the {\it smash product} of (marked) Eilenberg-MacLane  spaces that  represent cohomology, where, recall, 
the {\it smash product}  of  spaces  with marked points, say $T_1=(T_1, t_1)$ and $T_2=(T_2,t_2)$ is 
$$T_1\wedge T_2=T_1\times T_2/T_1\vee T_2$$
where the factorisation 
"$/T_1\vee T_2$" means   {\sl"with  the subset  $(T_1\times t_2)\cup (t_1\times T_2)\subset T_1\times T_2$ shrunk to a point"} (that serves to mark} $T_1\wedge T_2$).

In  fact, general cohomotopy  "measures"  (see 1.9)  and spectra defined with maps $\Psi\to T$  satisfy  
natural (obviously defined)  counterparts/generalizations  of  $[ \cap\smile] $  and   $[\min \smile]$, call them  $[ \cap\wedge] $ and  $[ \min\wedge] $ that are

 \vspace {1mm}

{\it On Grading  Cell Numbers.} 
Denote by   $N_{i\_cell}[\phi]$ the 
 minimal number  $N_i $ such that homotopy class  $[\phi]$ of maps $S\to \Psi$  factors as $S\to K\to \Psi$  where $K$ is a cell complex with (at most) $N_i$ cells {\it of dimension $i$} and observe that the total cell number is bounded by the sum of these,
$$N_{cell}[\phi]\leq\sum_{i=0,1,2,...}N_{i\_cell}[\phi].$$

Under what conditions on $\phi$  does the sum $\sum_{i}N_{i\_cell}[\phi]$ (approximately)  equal $N_{cell}[\phi]$? 

What are relations  between the cell numbers of the  covering  maps $\phi$  between     (arithmetic)  locally symmetric spaces $\Psi$ 
besides $N_{cell}\leq\sum_{i}N_{i\_ cell}$ ?\footnote{ The identity maps  $\phi=id:\Psi\to \Psi$ of locally symmetric spaces $\Psi$   seem  quite   nontrivial in this regard. On the other hand, general locally isometric maps $\phi:\Psi_1\to \Psi_2$ between symmetric spaces  as well as continuous maps $S\to \Psi$ of positive degrees, where $S$ and $\Psi$  are equidimensional manifolds with only $\Psi$ being locally symmetric, are also   interesting.  }

\section{Symmetries, Equivariant Spectra  and Symmetrization.}

. If  the energy function  $E$ on $\Psi$ is invariant under a  continuous  action of a group $G$ on $\Psi$ --  this happens frequently --   then  the relevant category is that of {\it $G$-spaces} $S$, i.e. of topological spaces  $S$ acted upon by $G$, where one works with  {\it $G$-equivariant}   continuous maps $\phi:S\to \Psi$,  {\it equivariant} homotopies,  equaivariinat  (co)homologies, decompositions, etc.\vspace {1mm}

Relevant  examples of this are provided by symmetric  energies $E=E(x_1,...,x_N)$ on Cartesian powers of  spaces, $\Psi=X^{\{1,...,N\}}$, such as 
  our  (ad hoc)  packing energy for  a metric space $X$,
$$ E\{x_1,..., x_i,...,x_N\}= \sup_{i\neq j=1,...,N}dist^{-1} (x_i,x_j)$$
that is invariant under the  symmetric group $Sym_N$.  It is often profitable, as we shall see later on, to exploit the symmetry under certain {\it subgroups} $G\subset Sym_N$.\vspace {1mm}

Besides the group $Sym_N$,   energies $E$ on $X^{\{1,...,N\}}$ are often   invariant under some  groups $H$ acting on   $X$, such as the  isometry group $Is(X)$  in the case of packings. 

If such a group  $H$ is compact, than  its role   is less significant than that of  $Sym_N$, especially for large $N\to \infty$; yet, if $H$ properly acts on 
a {\it non-compact} space $X$,    such as  $X=\mathbb R^m$  that is acted upon by its isometry group, then $H$ and its  action become essential.  

 \vspace {1mm}

{\it MIN-Symmetrized Energy.} An arbitrary function $E$ on a $G$-space $\Psi$ can be rendered $G$-invariant by taking a symmetric function of the numbers  $e_g=E(g(\psi))\in \mathbb R$, $g\in G$. Since we are  mostly concerned  with the order structure in $\mathbb R$,  our  preferred  symmetrisation is
$$E(\psi)\mapsto \inf_{g\in G} E(g(\psi)).$$

{\it Minimization with Partitions.} This  $inf$-symmetrization  does not fully depends on the action of $G$ but rather on the partition of $\Psi$ into orbits of $G$. In fact, given an arbitrary partition   $\alpha$ of $\Psi$ into subsets that we call {\it  $\alpha$-slices}, one  defines the function
$$E_{inf_\alpha}= inf_\alpha E :\Psi\to \mathbb R $$
where $E_{inf_\alpha}(\psi)$ equals the infimum of $E$ on the $\alpha$-slice that contains  $\psi$  for all $\psi\in \Psi$.
Similarly, one defines $E_{sup_\alpha}= sup_\alpha E$ with $E_{min_\alpha}$   and $ E_{max_\alpha}$ understood accordingly.

{\it Example:  Energies on Cartesian Powers.} The energy $E$ on $\Psi$  induces $N$ energies  on  the space  $\Psi^{\{1,...,N\}}$ of  $N$-tuples $\{\psi_1,...,a_i,...,a_N\}$,
that are  $$E_i:\{\psi_1,...,a_i,...,a_N\}\mapsto E(a_i).$$

It is natural, both from a geometric as well as from  a physical prospective, to  symmetrize by taking {\it the total energy} $E_{total}=\sum_i E_i$.  But in what follows we shall resort to 
$E_{min}=\min_i E_i=\min_i E(a_i)$ and use it for bounding  the total energy from below by 
$$E_{total}\geq N\cdot E_{min}.$$

For instance, we shall do it for families of $N$-tuples of balls $U_i$  in a Riemannian manifold $V$, thus bounding the $k$-volumes of $k$-cycles $c$ in the unions $ \cup_iU_i$, where, observe, 
$$vol_k(c)=\sum_ivol_k(c\cap U_i)$$
if the balls $U_i$ {\it do not intersect}.\vspace {1mm}


$$E(c)=\sum_i c\cap U_{x_i}.$$
  This, albeit obvious, leads, as we shall see later on, to  non-vacous  relations between\vspace {1mm}
  
   \hspace {0mm}  {\it    homotopy/homology spectrum of the $vol_k$-energy 
   on the space  $C_k(X;\Pi)$ 
   
    \hspace {0mm}  and

 equivariant  homotopy/homology of the spaces  of   packings of $X$ by $\varepsilon$-balls. }\vspace {1mm}

\section{Equivariant Homotopies of Infinite Dimensional Spaces.}

If we want to understand  homotopy spectra of spaces of "natural energies" on spaces of  infinitely many  particles in 
non-compact manifolds, e.g. in  Euclidean spaces, we need to extend the concept of the homotopy and homology spectra to infinite dimensional spaces $\Psi$, where infinite dimensionality is compensated by an additional structure, e.g. by an action of an infinite group  $\Upsilon$ on $\Psi$.

The simplest instance of this is where $\Upsilon$ is a countable group that we prefer to call $\Gamma$  and $\Psi=B^\Gamma$  be the space of maps $\Gamma\to B$ with the (obvious)  {\it shift action} of  $\Gamma$ on this $\Psi$, motivates the following definition (compare  \cite {bertelson-gromov}).
Let   $H^\ast$ be a graded algebra (over some field)   acted upon by  a   countable  amenable group $\Gamma$. Exhaust $\Gamma$ by finite {\it  F{\o}lner subsets} $\Delta_i\subset \Gamma$, $i=1,2,...$,  and,  given a finite dimensional graded  subalgebra $K=K^\ast \subset H^\ast$, let $P_{i,K}(t)$
denote the Poincare polynomial of the graded  subalgebra in $H^\ast$ generated by the $\gamma$-transforms
 $\gamma^{-1}(H_K^\ast) \subset H^\ast$ for all $\gamma \in \Delta_i$.

Define {\it  polynomial entropy} of the action of $\Gamma$ on  $H^\ast$ as follows.
$$  Poly.ent ( H^\ast: \Gamma   )=\sup_K \lim_{i\to \infty} \frac {1}{card(\Delta_i)}\log P_{i,K}(t).$$

Something of this kind  could be applied to subalgebras  $H^\ast\subset H^\ast(\Psi; \mathbb F)$, such as images and/or kernels of the restriction cohomology homomorphisms for (the energies sublevel) subsets $U\subset \Psi$,   IF   the following  issues are  settled. 
\vspace {1mm} 

\textbf {1. } In our example of moving balls or particles in $\mathbb R^m$  the relevant groups $\Upsilon$,  such as    the group of the orientation preserving  Euclidean  isometries   are connected and act {\it trivially} on the cohomologies of our spaces   $\Psi$.

For instance, let $\Gamma\subset \Upsilon$ be a discrete subgroups and $\Psi$ equal the {\it dynamic  $\Upsilon$-suspension of} $B^\Gamma$, that is $B^\Gamma\times \Upsilon$  divided by the diagonal action of $\Gamma$.
$$\Psi=\left ( B^\Gamma\times \Upsilon\right) /\Gamma.$$

The (ordinary) cohomology  of this space   $\Psi$  are bounded 
by those  of $B$ tensored by the cohomology of $\Upsilon/\Gamma$  that would  give {\it zero} polynomial entropy for finitely generated  cohomology algebras $H^\ast(B)$.

In order to have something more interesting, e.g.  the { \it   mean    Poincar\'e polynomial} equal that of $B^\Gamma$, which is the ordinary  \upshape Poincar\'e$(H^\ast (B))$, 
one needs 
  a definition of  some    { \it   mean (logarithm) of the  Poincar\'e polynomial} that might  be   {\it far from zero}  even if the ordinary cohomology of $\Psi$ vanish.

There are several  candidates for  such  { \it mean Poincar\'e polynomials}, e.g the one is suggested in section1.15 of my article {\it Topological Invariants of Dynamical Systems and Spaces of Holomorphic Maps}.

Another  possibility that is  applicable to the above $\Psi=\left ( B^\Gamma\times \Upsilon\right) /\Gamma$ with  {\it residually finite} 
groups $\Gamma$ is using finite $i$-sheeted  covering $\tilde \Psi_i$  corresponding to subgroups $\Gamma_i\subset \Gamma$ of order $i$ and taking the limit of  
$$\lim_{i\to \infty}\frac {1}{i}\log\mbox{\upshape Poincar\'e$(H^\ast (B))$}.$$,

(Algebraically, in terms of actions of  groups $\Gamma$ on  abstract graded algebras  $H^\ast$, this corresponds to taking 
 the normalised limit of logarithms of {\it $\Gamma_i$-invariant} sub-algebras in $H^\ast$; this brings to one's mind a possibility of a  generalisation of the above polynomial entropies   to {\it sofic groups} (compare \cite {arzhantseva}).

    \textbf {2. } The above numerical definitions of the polynomial  entropy   and of the mean Poincar\'e polynomials beg to be rendered in categorical  terms similarly to the ordinary entropy (see \cite {structure}).

  \textbf {3. } The spaces $\Psi_{\infty}(X)$  of  (discrete) infinite   countable subsets $\psi\subset X$ that  are meant to represent infinite  ensembles of  particles in    non-compact manifolds $X$, such as  $X=\mathbb R^m$, are more complicated than   
 $\Psi=B^\Gamma$, $\Psi=\left ( B^\Gamma\times \Upsilon\right) /\Gamma$ and other "product like" spaces studied  eralier.

These  $\Psi_{\infty}(X)$ may be seen as as limits of finite spaces   $\Psi_{N}(X_N)$ for $N\to 
\infty $ of $N$-tuples of points in compact manifolds $X_N$ where one has to chose suitable approximating  sequences 
$X_N$.

For instance, if $X=\mathbb R^m$ acted upon by some isometry group $\Upsilon$ of $\mathbb R^m$  one may use either the balls $B^m(R_N)\subset \mathbb R^m$ of radii $R_N= const\cdot R^{ N/\beta}$ in $\mathbb R^m$,  $\beta>0$  for $X_N$ or the tori
$\mathbb R^m/\Gamma_N$ with  the  lattices $\Gamma_N=const\cdot M\cdot \mathbb Z^m$ with 
$M=M_N \approx  N^\frac {1}{\beta}$ for some $\beta>0$.\footnote{The  natural value   is $\beta=m$ that make the volumes of $X_N$ proportional to $N$ but smaller   values, that correspond to ensembles of points in $\mathbb  R^m$ of {\it zero densities}, also make sense as we shall see later on.}
 
Defining  such limits and working out  functional    definitions of relevant structures the   limit spaces, collectively callused    $\Psi_{\infty}(X)$  are the  problems we need to solve  where, in particular, we need to

$\bullet $  incorporate  actions of the group  $\Upsilon$ coherently with  (some subgroup) of the infinite  permutations group acting on subsets $\psi\subset X$ of particles in $X$
 that represent points in $\Psi_\infty(X)$  
 
\hspace {-3mm} and 
 
 $\bullet $  define   (stochastic?)  homotopies and (co)homologies in the spaces $\Psi_\infty(X)$, where   these  may be associated to  limits of families of $n$-tuples  $\psi_{\it P_N}\subset   X_N$ parametrised by some $P_N$ where $dim(P_N)$ may tend to infinity for $N\to \infty$.\footnote{We shall  meet  families of dimensions  $dim(P_N)\sim N^\frac {1}{\gamma}$ where ${\gamma}+{\beta}=m$ for the above $\beta$.}

     \textbf {4. } Most natural energies $E$ on infinite particle  spaces $\Psi_\infty(X)$ are everywhere infinite\footnote{In the optical  astronomy, this is called {\sl Olbers'  dark night sky paradox.}}
 and defining "sublevels" of  such $E$ needs  attention.
 
\section {Symmetries, Families and  Operations Acting on Cohomotopy Measures.}

 {\it Cohomotopy "Measures}". Let $T$ be  a space  with a distinguished  {\it marking point} $t_0\in T$,  let 
 $H^\circ(\Psi;T)$ denote the set of homotopy classes of maps $\Psi\to T$ and 
 define  the "$T$-measure" of an open subset $U\subset \Psi$,
 $$\mu^T(U)\subset H^\circ(\Psi;T),$$
as the set of homotopy classes of maps $\Psi\to T$ that send the complement $\Psi\setminus U$ to $t_0.$

For instance,  if $T$ is the Cartesian product of Eilenberg-MacLane spaces  $K(\Pi;n)$, $n=0,1,2,...,$, then $H^\circ(\Psi;T)=  H^\ast(\Psi;\Pi)$ and $\mu^T$  identifies with the  (graded cohomological)  ideal valued "measure"   $U\mapsto \mu^\ast (U;\Pi)\subset H^\ast(\Psi;\Pi)$ from section 1.5.

Next, given a category $\cal T$ of  marked spaces $T$ and homotopy classes of  maps between them,   denote by $\mu^{\cal T}(U)$ the totality of the  sets 
 $\mu^T(U)$,  $T\in \cal T$, where  the category $\cal T$ acts on $\mu^{\cal T}(U)$  via composition  
$\Psi\overset {m}\to T_1\overset {\tau}\to T_2$ for all $m\in  \mu^T(U)$  and $\tau\in \cal T$.

For instance,  if $\cal T$ is   a category of  Eilenberg-MacLane spaces $K(\Pi;n)$,  this amounts to the natural  action of the (unary) cohomology operations (such as Steenrod squares $Sq^i$ in the case  $\Pi=\mathbb Z_2$) on  ideal valued measures. 

\vspace {1mm}
The above definition can be adjusted for spaces  $\Psi$ endowed with additional structures.

 For  example, if   $\Psi$ represents a {\it family of spaces} by being endowed with a partition $\beta$  into closed subsets -- call them  {\it $\beta$-slices} or {\it fibers}  -- then  one  
restricts to  the space of maps $\Psi\to T$ { \it constant on these slices}  (if $T$ is also partitioned, it would be logical to deal with maps sending slices to slices)  defines    $H^\circ_\beta(\Psi;T)$  as the set of the homotopy classes of these slice-preserving maps and 
accordingly  defines 
$\mu_\beta^{\cal T}(U)\subset H^\circ_\beta(\Psi;T). $ 

Another kind of a relevant structure is an action of a group $G$ on $\Psi$. Then one may (or may not) work with categories $\cal T$ of $G$-spaces $T$ (i.e. acted upon by $G$) and perform
 homotopy, including (co)homology, constructions equivariantly. Thus, one defines equivariant $T$-measures  
$\mu^T_G(U)$ for $G$-invariant subsets $U\subset \Psi$.

(A group action on  a space , defines a partition of this space  into orbits, but this is a weaker structure than that of the the action itself.)\vspace {1mm}

{\it  Guth' Vanishing Lemma.} The supermultiplicativity property of the cohomology measures with arbitrary coefficients $\Pi$ (see 1.5) for spaces $\Psi$  acted upon by finite groups $G$ implies that 
$$\mu^\ast\left (\bigcap_{g\in G} g(U;\Pi)\right)\supset \underset {g\in G} \smile\mu^\ast (g(U;\Pi))$$
for all open subset $U\subset \Psi$.

This, in the case $\Pi=\mathbb Z_2$  was generalised by Larry Guth for families of spaces parametrised by spheres $S^j$ 
as follows.

Given a space $\Psi$  endowed with a partition
 $\alpha$,  we say that a subset in $\Psi$ is {\it $\alpha$-saturated} if it equals the  union of some $\alpha$-slices in $\Psi$ and 
define two operations  on subsets  $U\subset \Psi$,
$$U\mapsto\cap_\alpha (U)\subset U\mbox { and   }  U\mapsto\cup_\alpha (U)\supset U, $$ 
where  

$\cap_\alpha (U)$ is  the {\it maximal  $\alpha$-saturated subset that is contained in $U$}  

\hspace {-6mm} and

 $\cup_\alpha (U)$  is the {\it minimal $\alpha$-saturated subset
that  contains $U$.}

Let, as in the case  considered by Guth,  $\Psi= \Psi_0\times S^j$  where $S^j\subset \mathbb R^{j+1}$ is the $j$-dimensional sphere, let  
$\alpha$ be 
 the partition into the orbits of   $\mathbb Z_2$-action on $\Psi$ by $(\psi_0,s)\mapsto (\psi_0,-s)$ (thus, "$\alpha$-saturated" means "$\mathbb Z_2$-invariant") and let $\beta$ be the partition into the fibres of the projection $\Psi\to \Psi_0$  (and "$\beta$-saturated" means "equal the pullback of a subset in $\Psi_0$").

Following Guth, define \vspace {0.5mm}

\hspace {13mm}$Sq_j:H^{\ast\geq j/2}(\Psi;\mathbb Z_2) \to H^\ast(\Psi;\mathbb Z_2)$  by  $Sq^j : H^p\to H^{2p-j}$   \vspace {0.5mm}

\hspace {-6mm} and formulate his "Vanishing Lemma" in $\mu_\beta$-terms as  follows,\footnote{Guth formulates his lemma in terms of the complementary set $V=\Psi\setminus U$:

{\it if a cohomology class $h\in H_\beta^\ast(\Psi;\mathbb Z_2)$ vanished on $V$, then $St_j(h)$ vanishes on 
$\cap_\beta(\cup_\alpha(V))$.}}
$$ \mu_\beta^\ast \left(\cup_\beta(  \cap_\alpha (U));\mathbb Z_2\right)  \supset Sq_j( \mu_\beta^\ast(U;\mathbb Z_2))\subset H_\beta^\ast (\Psi;\mathbb Z_2), \leqno [\cup  \cap]$$
where, according to our notation,  $H_\beta^\ast (\Psi;\mathbb Z_2)\subset H^\ast (\Psi;\mathbb Z_2)$ equals the image of 
$H^\ast (\Psi_0;\mathbb Z_2)$ under the cohomology homomorphism induced by the projection $\Psi\to \Psi_0$.

If $E:\Psi\to \mathbb R$ is an energy function, this lemma  yields the  lower bound on the  {\it maxmin}-energy\footnote{Recall that  $min_\alpha E(\psi)$, $\psi\in \Psi$, denotes the minimum of $E$ on the $\alpha$-slice containing $\psi$ and $max_\beta $ stands for similar maximisation with $\beta$ (see 1.12).}
$$ E_{max_\beta min_\alpha}=  max_\beta min_\alpha E $$  
evaluated at the cohomology class $St_j(h)$, $h\in H^\ast_\beta(\Psi;\mathbb Z_2)$:
$$ E^\ast_{max_\beta min_\alpha}(St_j(h))\geq  E^\ast(h).\leqno{[maxmin]}$$ 
.


{\it Question.}  What are  generalisations of   $[\cup  \cap]$ and $[maxmin]$  to other cohomology and   cohomotopy measures 
on spaces  with partitions $\alpha, \beta, \gamma$,...? 

\section{Pairing Inequality  for Cohomotopy Spectra.}

Let $\Psi_1, \Psi_2$ and $\Theta$  be topological spaces and let 
$$\Psi_1\times \Psi_2\overset { \circledast}\to \Theta$$
be a continuous map where we write
$$\theta=\psi_1\circledast \psi_2\mbox { for }  b=\circledast (\psi_1,\psi_2).$$
For instance, 
 composition  $\psi_1\circ \psi_2:X\to Z$ of morphisms $X\overset {\psi_1} \to Y\overset {\psi_2}\to Z$ in  a topological category defines  such a map
 between sets of morphisms,
$$mor(X\to Y)\times mor(Y\to Z)\overset{\circledast}\to mor(X\to Z). $$

A more relevant example  for us is the following\vspace {2mm}

\hspace {35mm}{\sc Cycles $\bigtimes$ packings.}\vspace {1mm}

Here, $\Psi_1$ is a  space of locally diffeomorphic maps  $U\to X$  between
 manifolds 
 
 $U$ and $X$,  \vspace {0.5mm}

$\Psi_2$ is the space of cycles   in $X$  with some coefficients $\Pi$,  \vspace {0.5mm}

$\Theta$ is the space of cycles in  $U$  with the same coefficients,  \vspace {0.5mm}

  $\circledast $  stands for  {\sl "pullback"}
 $$\theta=\psi_1\circledast \psi_2=_{def}\psi_1^{-1}(\psi_2)\in \Theta.$$

This  $U$ may  equal the disjoint unions of $N$  manifolds $U_i$ that, in the spherical packing problems, would  go to balls  $B_x(r)\subset X$; since we want these balls {\it not to intersect},   we take the space of {\it  injective} maps $U\to X$  for $\Psi_1$.\vspace {1mm}

If the manifold $X$ is {\it parallelizable} and the  balls $B_{x_i}(r)\subset X$ all  have some  radius   radius $r$  smaller then the injectivity radius of $X$, then  corresponding $U_i=B_{x_i}(r)$  can be identified (via the exponential maps) with the $r$-ball $B^n(r)$ in the Euclidean space $\mathbb R^n$, $n=dim(X)$. Therefore,   the space   of $k$-cycles in $U=
\bigsqcup_i U_i=B^n(r)$ equals in this case  the Cartesian $N$-th power of  this space for the $r$-ball:
$$\Theta={\cal C}_k(U;\Pi)=({\cal C}_k(B^n(r);\Pi))^N.$$

{\it Explanatory  Remarks.}  (a)  Our "{\sl cycles}"  are defined as {\it  subsets} in relevant manifolds  $X$ and/or   $U$  with {\it $\Pi$-valued functions} on these subsets.

 (b) In the case of {\it open} manifolds, we speak of  cycles with {\it infinite supports}, that, in the case of compact manifolds with boundaries or of  open subsets  $U\subset X$, are, essentially, {\it cycles modulo the boundaries} $\partial X$.

(c) {\sl "Pullbacks of cycles"} that preserve their codimensions   are defined, following Poincar\'e for  a wide  class of  smooth {\it generic} (not necessarily equividimensional) maps $U\to X$ (see \cite {manifolds}).

(d) It  is easier  to work with {\it cocycles} (rather than with cycles) where contravariant functoriality needs   no extra assumptions  on spaces and maps in question (see \cite {singularities2}).\vspace {2mm}

Let $h^T$ be a  (preferably non-zero)  cohomotopy class in $\Theta$, that is a homotopy class of non-contractible maps $\Theta\to T$ for some space $T$, (where "{\sl cohomotopy"} reads  "{\sl cohomology}" if $T$ is an Eilenberg-MacLane space) and let $$h^\circledast=\circledast\circ h^T:[\Psi_1\times \Psi_2\to T]$$ be the induced class on $\Psi_1\times \Psi_2$, that  is the homotopy class of the composition of the maps $\Psi_1\times \Psi_2\overset{\circledast}\to \Theta\overset{h^T}\to T$. 

(Here and below, we do not always notationally distinguish  {\sl maps} and   {\sl homotopy classes} of  maps.)

Let $h_1$ and $h_2$ be homotopy classes of maps $S_1\to \Psi_1$ and  $S_2\to \Psi_2$ for  some spaces  $S_i$, $ i=1,2$,

 (In the case where  $h^T$ is a {\it cohomology}  class, these   $h_i$ may be replaced by {\it homology} -- rather than homotopy -- classes represented by these maps.) 
 
 Compose the three maps, 
$$S_1\times S_2\overset {h_1\times h_2}\to \Psi_1\times \Psi_2 \overset{\circledast}\to \Theta \overset {h^T} \to T,$$
and denote the homotopy class of the resulting map $ S_1\times S_2\to T$ by
$$[h_1\circledast h_2]_{h^T}=h^\circledast\circ (h_1\times h_2) : [S_1\times S_2\to T]$$

Let $\chi=\chi(e_1,e_2)$  be a   function in two real variables that is monotone unceasing in each variable. Let   $E_i:\Psi_i\to \mathbb R$, $i=1,2$, and $F:\Theta\to \mathbb R$ be (energy) functions  on the spaces $\Psi_1, \Psi_2$ and $ \Theta$, such that the  $\circledast$-pullback of $F$ to
$\Psi\times \Theta$ denoted 
$$F^\circledast=F\circ \circledast:\Psi_1\times \Psi_2\to \mathbb R$$ satisfies 
$$F^\circledast(\psi_1,\psi_2)\leq \chi(E(\psi_1),E(\psi_2)).$$

 In other words, the  $\circledast$-image of the product of the sublevels $$(\Psi_1)_{e_1}=E_1^{-1}(-\infty, e_1)\subset \Psi_1\mbox {  and  }(\Psi_2)_{e_2}=E_2^{-1}(-\infty, e_2)\subset \Psi_2$$ 
is contained in the $f$-sublevel  $B_f=F^{-1}(-\infty, f)\subset \Theta$ for $f=\chi(e_1,e_2)$,
 $$\circledast\left ((\Psi_1)_{e_1}\times (\Psi_2)_{e_2}\right)\subset \Theta_{f=\chi(e_1,e_2)}.$$

\vspace {1mm} 

\hspace {34mm}{\sc $\circledast$-Pairing Inequality.}\vspace {1mm}

Let
$[h_1\circledast h_2]_{h^T}\neq 0,$
that is 
{\it the composed map $$S_1\times S_2\to \Psi_1\times \Psi_2\to \Theta \to T$$
is  non-contractible.}
Then   the values of $E_1$ and $E_2$ on the homotopy classes $h_1$ and $h_2$  are {\it bounded from below} in terms of a lower bound on  $F^\circ[h^T]$  as  follows.

$$\chi( E_{1\circ}[h_1], E_{2\circ}[h_2])\geq  F^\circ[h^T].\leqno {[_{\circ\circ}\geq ^\circ  ]}$$
In other words

$$\big( E_{1_\circ}[h_1]\leq e_1\big) \&\big ( E_{2\circ}[h_2]\leq e_2\big) \Rightarrow \big(F^\circ[h^T] \leq \chi(e_1,e_2)\big) $$
for all real numbers $e_1$ and $e_2$;  thus,\vspace {0.5mm}
 
\hspace {10mm} {\it  \textbf {upper} bound  $E_1^\circ[h_1]\leq e_1$ $\bigplus$  \textbf {lower} bound  $F^\circ[h^T]\geq  \chi(e_1,e_2)$  

\vspace {0.5mm}

\hspace {48mm}yield\vspace {0.5mm}

 \hspace {30mm}  \textbf {upper} bound   $E_2^\circ[h_2]\geq e_2$,}\vspace {1.1mm}
 
 \hspace {-6mm}where, observe, $E_1$ and $E_2$ are interchangeable in this relation. \vspace {1mm}

In fact,  all one needs for   verifying $ [_{\circ\circ}\geq ^\circ  ]$ is  unfolding the definitions.  

Also  $ [_{\circ\circ}\geq ^\circ  ]$ can be  
  visualised without an explicit use of  $\chi$  by looking at 
the {\it $h^\circledast$-spectral line} in the $(e_1,e_2)$-plane
$$\Sigma_{h^\circledast}=\partial \Omega_{h^\circledast}\subset \mathbb R^2$$
 (we met this $\Sigma$  section 1.3) where
$\Omega_{h^\circledast}\subset \mathbb R^2$ consists 
of the pairs $(e_1,e_2)\in \mathbb R^2$ such that the restriction of $h^\circledast$ to 
the Cartesian product of the sublevels $\Psi_{1e_1}=E_1^{-1}(-\infty, e_1)\subset \Psi_1$ and  
$\Psi_{2 e_2}=E_2^{-1}(-\infty, e_2)\subset \Psi_2$ vanishes, 
$$h^\circledast_{| \Psi_{1 e_1}\times \Psi_{2e_2}}=0.$$

\section{Inequalities between Packing Radii,  Waists and Volumes of Cycles   in a  Presence of  Permutation Symmetries.}

The most essential aspect of the homotopy/homology structure in the   space of packings of a space $X$  by  $U_i\subset X$, $i\in I$,  is associated   with  the permutation group $Sim_N=aut(I)$  that acts on these spaces.

 A proper description of this needs  a use of the concept of "homological substantionality for {\sl variable} spaces" as in section 8. This  is  adapted to the present situation  in the definitions below. \vspace {1mm}

Let  \vspace {1mm}

$\bullet $ $\Psi_1$ be   space of $I$-packing $\{U_i\}$  of $X$ by disjoint open subsets $U_i\subset X$, $i\in I$,

 \hspace {1,5mm}for some set $I$ of cardinality $N$, e.g. by  $N$ balls  of radii $r$, $N=card (I)$,   \vspace {0.6mm}

$ \bullet $   ${\cal U}^\cup$ be 
 the space of  pairs $\big(\{U_i\},u\big )_{i\in I}$, $i\in I$,  where $\{U_i\}\in \Psi_1$   and $u\in \bigcup_{i\in I}U_i$, 

$ \bullet $  ${\cal U}^\times $ be 
 the space of  pairs $\big(\{U_i\},\{ u_i\} \big )_{i\in I}$,  where $\{U_i\}\in \Psi_1$   and $u_i\in U_i$,  \vspace {0.6mm} 
 
$ \bullet $ ${\cal C}^\cup_k$ be the space of $k$-cycles with $\Pi$ coefficient  in the unions 
 of $U_i$ for all 
 
\hspace {1,5mm}  packing  $\{U_i\}\in \Psi_1$ of $X$,
$${\cal C}^\cup_k=\bigcup_{\{U_i\}\in \Psi_1}{\cal C}_k(\bigcup_{i\in I} U_i;\Pi),$$

$ \bullet $  ${\cal C}^\times _{Nk}$ be the space of $N\cdot k$-cycles with the $N$-th 
tensorial power 
 coefficients  
 
 \hspace {1,5mm} in the 
 Cartesian products of  $U_i$ for all  $\{U_i\}\in \Psi_1$,
 $${\cal C}^\times _{Nk}=\bigcup_{\{U_i\}\in \Psi_1}{\cal C}^\cup_{Nk}\left (\bigtimes_{i\in I} U_i;\Pi^{\otimes N}\right).$$

\vspace {1mm}

The four spaces  ${\cal U}^\cup$,     ${\cal U}^\times$,  ${\cal C}_k^\cup$ and   ${\cal C}_{NK}^\times$ tautologically  "fiber" over    $\Psi_1$    by the maps denoted 
$\varpi^\cup$,  $\varpi^\times,$   $\varpi_k^\cup$,   and  $\varpi_{Nk}^\times;$ 
 besides, the  Cartesian products of cycles $C_i\subset U_i$  defines  an    embedding}   
$${\cal C}_k^\cup \hookrightarrow  {\cal C}^\times _{Nk}.$$

Observe that 

 { \tiny ${\blacksquare}$ }  the fibers $\bigtimes_iU_i$ of the  map  $\varpi^\times :  {\cal U}^\times  \to \Psi_1$ for $(\big(\{U_i\},\{ u_i\} \big ) \overset  { \varpi^\times} \mapsto \{U_i\}$ are what we call "variable spaces" in section 8 where  
$ {\cal U}^\times$ playing the role of $\cal X$ and $\Psi_1$   that  of $Q$ from  section 8,

 { \tiny ${\blacksquare}$ }  the above   four spaces are naturally/tautologically  acted upon, along with $\Psi_1$,  by   the symmetric group $Sym_N$ of permutations/automorphisms of the index set $I$  and the above four maps  from these spaces to  $\Psi_1$ are {\it $Sym_N$-equivariant,}

 { \tiny ${\blacksquare}$ }  the  $G$-factored  map $ \varpi^\times$, denoted 
 $$\varpi^\times_{/G}:  {\cal U}^\times/G \to \Psi_1/G,$$
has the same fibres as  $\varpi^\times$ , namely, the products $\bigtimes_i U_i$.
However, if $card(G)>1$ and $dim(X)=n>1$,  then   $\varpi^\times_{/G}$ is a {\it non-trivial }fibration,  even  if all  $U_i$ equal  translates of a ball $U=B^n(r)$ in  $X=\mathbb R^n$,  where the  ${\cal U}^\times= 
\Psi_1\times U^I$ and where $\varpi^\times: \Psi_1\times U^I \to \Psi_1$ equals the coordinate projection.



\vspace {1mm}


 
 The   "$\circledast$-pairing"  described in the previous section via  the  intersection of  cycles in $X$ with  $U_i\subset X$ for $\{U_i\}_{i\in I}\in \Psi_1$ followed by taking the Cartesian product of these intersections defines a pairing
 $$\Psi_1\times\Psi_2 \overset {\circledast^\times} \to \Theta^\times={\cal C}^\times _{Nk}=\bigcup_{\{U_i\}\in \Psi_1} {\cal C}_{Nk}\left (\bigtimes_{i\in I} U_i;\Pi^{\otimes N}\right),$$
where, recall,  $\Psi_1$ is the space of packings of $X$ by  $U_i\subset X$, $i\in I$, and 
  $\Psi_2={\cal C}_k(X;\Pi)$ is the space of cycles in $X$.

And since this  map  ${\circledast^\times}$   is equivariant for  the   action of  $Sym_N$ on   $\Theta^\times$ and on the first factor in   $ \Psi_1\times\Psi_2 $, it  descends to 
 $$  {\circledast^\times_{/G}}  : \Psi_1/G\times\Psi_2 \to \Theta^\times/G$$
for all subgroups $G\subset Sym_N$.

\vspace {3mm}

 \hspace {3mm}  {\sc  Detection Of Non-Trivial Families of Cycles and of Packings.}\vspace {2mm}

A family $S_2\subset \Psi_2$ of cycles in $X$ is called  {\it  homologically $G$-detectable} by (a family of) packings of $X$, if 
 there exists a   family  $S_1$ of cycles in  $\Psi_1/G\times\Psi_2$ such that the corresponding  family of  product $N\cdot k$-cycles in  "variable spaces" $\bigtimes_i U_i$, that is a map from $S_1\times S_2$ to    
 $ \Theta^\times/G$ by  ${\circledast^\times_{/G}}$ is  {\it homologically substantial.}

Recall (see 8)   that this   substantiality  is non-ambigously defined if  ${\circledast^\times_{/G}}$ is a fibration with {\it contractible fibres}, which is  the case  in our examples where $U_i\subset X$ are topological $n$-balls and observe that  {\it unavoidably} variable nature
 of $\bigtimes_i U_i$ is due to non-triviality of the fibration    $\varpi^\times_{/G}:  {\cal U}^\times/G \to \Psi_1/G$  for permutation groups    $G\neq \{id\}$ acting  on packings  that is most  essential for  what we do.


\vspace{2mm}

In his "{\sl Minimax-Steenrod}" paper   Guth  shows  that all  homology classes $h_2$ of the space $\Psi_2={\cal C}_k(B^n;
 \mathbb Z_2)$  of relative  $k$-cycles in the $n$-ball are     {\it $G$-detectable by   families $S_1=S_1(h_2)$  of packings} of $B^n$ by sufficiently  small $\delta$-balls  $U_i=B_{x_i}^n(\delta) \subset B^n=B^n(1)$, for some $2$-subgroup $G =G(h_2)  \subset Sym_N$.

  This is established with the help of  "Vanishing Lemma" stated in section 14; where, if understand it correctly, the detective power of  such a  family  $S_1\subset \Psi_1/G$ (of  "moving packings" of $X$ by $N$ balls)   is due to   {\it non-vanishing} of some cohomology class in  $\Psi_1/G$ that  comes  from the classifying space ${\cal B}_{cla}(G)$ via the classifying map 
  $ \Psi_1/G\to {\cal B}_{cla}(G)$.
  
(Probably, a  proper incorporation of the cohomology coming from $X$ would imply similar "detectability" {\it  all} manifolds $X$, where it is quite obvious for {\it parallelizable} $X$. )

Apparently Guth' "Vanishing Lemma" shows  that, every  $\mathbb Z_2$-homology class $h_\ast$  in  $K(\mathbb Z_2, m)$    equals the image of a homology class from   the classifying space  ${\cal B}_{cla}(G)$ of some  finite $2$-groups $G=G(h_\ast)$ under a map  ${\cal B}_{cla}(G) \to K(\mathbb Z_2, m)$. Equivalently, this means that given {\it a non-trivial} $\mathbb Z_2$-cohomology operation $op$ from degree $m$ to $n$,  (that is, necessarily,  a polynomial combination of Steenrod squares) there exists a class $h\in  H^m({\cal B}_{cla}(G);\mathbb Z_2)$ for some finite  $2$-group $G$,  such that $op(h)\neq 0$.

 Possibly,  this  is also true(?)   for other  primes  $p$  (which, I guess,  must be known to people working on cohomology of $p$-groups).  
 


On the other hand it seems unlikely that the (co)homologies of spaces of $I$-packing  for all finite sets $I$ and/or  of  (the classifying spaces of) all $p$-groups  are fully  detectable by the  (co)homologies of  the space ${\cal C}_{k} ( B^n; \mathbb Z_p) $ of cycles in the  $n$-ball,  where, recall,    the space  
${\cal C}_{k} (  B^n; \mathbb Z_p)$ is homotopy equaivalent to Eilenberg-MacLane's $K(\mathbb Z_p, n-k)$.

\vspace {2mm}

\hspace {3mm}{\sc Pairing Inequalities  Between $k$-Volumes and Packing Radii.}\vspace {1mm}

 Let $S_1\subset \Psi_1$ be a $G$-invariant family of $I$-packings of a Riemannian manifold  $X$  by open subsets  $U_i=U_{i,s_1}$, $i\in I$, $s_1\in S$, where $G$ is a subgroup of the group $Sym_N=aut(I)$, and let $S_2$ be a family of $k$-cycles (or more general  "virtually $k$-dimensionl entities")  $Y=Y_{s_2}$   in $X$.
 
 Let  the  coupled  family  that is a map from $S_1/G\times S_2$ to    
 $ \Theta^\times/G$ by  ${\circledast^\times_{/G}}$ be  {\it homologically substantial.}
 Then, by the definition of " the wast of a variable space "  (see section 8) the supremum of the volumes of $Y_{s_2}$
 is bounded from below by 
 $$\sup_{s_2\in S_2} vol_k(Y_{s_2})\geq \sum_{i\in I} waist_k(U_{i,s_1}).$$
 In particular, since small balls of radii $r$ in $X$ have $k$-$waists\sim r^k$
 the above implies
 $$\frac {\inf_{s_1\in S_1}\sum_{i\in I}inrad_k(U_{i,s_1})^k} {\sup_{s_2\in S_2}vol_k(Y_{s_2})}\leq const(X),$$ 
where $inrad(U)$, $U\subset X$ denote the radius of the largest ball contained in $U$.

This inequality, in the case where  $U_i\subset B^n= B^n(1)\subset \mathbb R^n$  are  Euclidean $r$-balls,  is used by Guth in \cite {guth-steenrod},  (as it was  mentioned earlier, for obtaining a (nearly sharp)  {\it lower bound} on the $k$-volume spectrum of the $n$-ball that is on  $\sup_{s_2\in S_2} vol_k(Y_{s_2})$ for families $S_2$ of $k$-cycles $Y$ in  $X=B^n$ such that  a given  cohomology class 
$h\in H^\ast\left ({\cal C}_k(X ;\mathbb Z_2);\mathbb Z_2\right)$ does not vanish on this $S_2\subset {\cal C}_k(X ;\mathbb Z_2)$. This is  achieved  by constructing a $G$-invariant  family $S_1=S_1(h)$  of $I$-packings of $X$ by $r$-balls  for some $2$-group $G\subset Sym_N$, $N=card(I)$, such that the coupled  family  
$$S_1/G\times S_2\overset {\circledast^\times_{/G}}\to    \Theta^\times/G$$  
is  homologically substantial  und such that the radii $r$ of the balls $U_i=B_{x_i}(r)\subset X$ are sufficiently large.

{\it Conversely}, and this is what we emphasise in this paper, one can {\it bound from above} the radii $r$ of the balls in a $G$-invariant  family $S_1$ of $I$-packing of $X$ by $r$-balls  in terms of a cohomology class $h'$ in the space  $\Psi_1/G$ (that is of  packings/$G$)  such that $h'$ does not vanish on $S_1$, where one can use for this purpose suitable   families $S_2=S_2(h')$ of $k$-cycles in $X$ with possibly small volumes, e.g. those used by Guth in \cite {guth-steenrod} for his {\it upper bound} on the $vol_k$-spectra.

 Thus, in the spirit of Guth' paper,  one  gets  such bounds for $N=card(I)\to \infty$ and 
some  $h'\in H^\ast( \Psi_1/G;\mathbb Z_2)$ for certain $2$-groups $G\subset Sym_N$,  where        $ord(G)\to \infty $ and $deg(h')\to \infty$ along with $N\to \infty$.
 
 And even though  the packing spaces $\Psi_1$ looks quite innocuous being the complements to the $2r$-neighbourhoods to the unions of the  diagonals in $X^I$  there is no(?) parent  alternative method to obtain such bounds.\footnote {A natural candidate for such method would be the Morse theory for the distance function to the union of the diagonals in $X^I$, see \cite {baryshnikov}, but this does not(?) seem to yield such bounds.}

 \vspace {2mm}
 {\it Multidimensional  Rendition of the Pairing Inequalities.}  
  Let $\Psi$ be a space of pairs $\psi=(\{U_i\}_{i\in I}, Y)$ where $U_i\subset X$ are disjoint open subsets and $Y\subset X$ is a $k$-cycle, say with $\mathbb Z$ or $\mathbb Z_2$ coefficients.

Let   $inv(U) $  be some  geometric invariant of open subsets  $U\subset X$  that is monotone under inclusions between subsets,  e.g. the $n$-volume, inradius,
some kind of waist of $U$, etc.

Define   ${\cal E}=(E_0, E_1, ..., E_N):\Psi\to \mathbb R^{N+1}$ by  
$$E_0(\{U_i\}_{i\in I}, Y)=vol_k(Y)\mbox {  and } E_i(\{U_i\}_{i\in I}, Y)=inv(U_i)^{-1}$$
assume that $\Psi$ is invariant under the action of some subgroup $G\subset Sym_N$ that permutes $U_i$, and observe that all of the above can be  seen in  the   kernels  of the cohomology homomorphisms from $H^\ast(\Psi/G)$ to $H^\ast(\Psi_{e_0,e_1,...,e_N}/G)$ for the  subsets $ \Psi{e_0,e_1,...,e_N} \subset \Psi$  defined by the inequalities
$E_i(\psi) < e_i $, $\psi=(\{U_i\},Y)$, $i=0,....N$,  as in the definition of the  multidimentional
 (co)homology spectra in  sections 4,10,11. 

For example, besides the pairing inequalities for packing by balls, this also allows an encoding of similar inequalities for {\it convex partitions} from \cite {waists}.

 \section { Sup$_\vartheta$-Spectra, Symplectic Waists and  Spaces of Symplectic  Packings.} 
  
 Let $\Theta $  be a set of metrics  $\vartheta$ on   a topological space $X$ and define {\it $sup_\vartheta$-invariants } of $(X,\Theta)$   
as     {\it the suprema}  of the corresponding invariants $(X,\vartheta)$  over all $\vartheta\in \Theta.$ (In many cases, this definition  makes sense for more general classes  $\Theta $ of metrics spaces that do not have to be homeomorphic to a fixed $X$.)  \vspace {1mm}

{\it Problem.}   Find general criteria for finiteness of  $sup_\vartheta$-invariants.

 \vspace {1mm}

 {\it Two Classical Examples: Systoles and Laplacians.} (1)  Let $\Theta $ be the space of Riemannian metrics on the $2$-torus $X$  with 
 $$sup_{\theta\in \Theta} area_\theta(X)\leq 1.$$
 Then  the   {\it $\sup_\vartheta$-systole$_1$ of $(X,\vartheta)$  is $<\infty$.}
 
  In fact,  
 $$sup_\vartheta\mbox {-}systole_1(X,\Theta)=\sqrt{\frac {2}{\sqrt {3}}}$$
by {\it Lowener's torus inequality}  of 1949.
 
This means that  all  toric surfaces $(X,\vartheta)$   of {\it unit areas} admit closed {\it non-contractible} curves of lengths $\leq \sqrt{\frac {2}{\sqrt {3}}}$, where, observe,  
 the equality  $systole_1= \sqrt{\frac {2}{\sqrt {3}}}$ holds for $\mathbb R^2$ divided by the {\it hexagonal lattice}. (See  Wikipedia article on systolic geometry and references therein  for further information.) 

(2) Let $\Theta$ be the space of Riemanian metrics $\vartheta$ on the $2$-sphere $X$  with
$area_\vartheta(X)\geq 4\pi$ (that is the area if the unit sphere).  Then 

\hspace {5mm}{\it the first    $\sup_\vartheta$-eigenvalue of the Laplace operator on $(X,\Theta)$ is $<\infty$. }

In fact $sup_\vartheta$-$\lambda_1(X,\Theta)=  2$, that is the first eigenvalue of the Lapalce operator on the unit sphere,  by { \it the Hersch inequality} of 1970.

\vspace {1mm}

{\it Symplectic Area Spectra and Waists.} Let $X$ be a smooth manifold of even dimension $n=2m$ and let $\omega=\omega (x)$ be a differential 2-form on $X$. 
 
 A   {\it Riemannian metric}  $\vartheta$ on $X$,  is called {\it  adapted to} or {\it  compatible with} $\omega$, if \vspace {0.7mm}
 
$$area_\vartheta(Y)\geq |\int_Y\omega | \leqno {\bullet_{\geq \omega }}$$ for all  smooth oriented surfaces $Y\subset X$; 
 \vspace {0.7mm}
 
 2. the  $n$-volume  $d_\vartheta x$ element, satisfies 
  $$d_\vartheta (x)\leq |\omega^m| \mbox { at all ponts $x\in X$} .\leqno{\bullet_{\leq \omega^m }}, $$
 that is $vol_\vartheta(U)\leq \int_U\omega^m$ for all open subsets $U\subset X$.\footnote {The {\it inequalities} $area_\vartheta(Y)\geq |\int_Y\omega |$ and $d_\vartheta (x)\leq  |\omega^m|$ imply the {\it equality} $ d_\vartheta x= |\omega^m|$ and if $\omega^m(x)\neq 0$, then  there is a $\mathbb R$-linear isomorphism of the  tangent space $T_x$ to $ \mathbb C^m$, such that    $(\sqrt{-1}\omega,\vartheta)_x$ go to the imaginary  and the  real parts of the diagonal  Hermitian form on $ \mathbb C^m$.  But $\bullet_{\leq \omega^m }$ is better adapted for   generalisations than $\bullet_{= \omega^m }$.}\vspace {1mm}

 {\it Question.}  Which part of  the (suitably factorized/coarsend)  homotopy/homology  area  spectra of $(X,\vartheta)$ remains finite after taking suprema  over  $\vartheta\in \Theta(\omega)$?  \vspace {1mm}
 
{\it Partial  Answer Provided by the Symplectic Geometry.} The form $\omega $ is called {\it symplectic } if it is {\it closed}, i.e $d\omega=0$,  and  $\omega^m=\omega^m(x)$ does not vanish on $X$. In this case $X=(X,\omega)$ is called {\it a symplectic  manifold.}

  {\it The  symplectic  $k$-waist} of    $X$ may be defined as   the {\it supremum of the  $k$-waists}\footnote{ This may be any kind of a waist  defined with families $S$  of   "virtually $k$-dimensional entities" of a particular kind and with a given type of  homological substantiality required from $S$, where  the most relevant for the symplectic geometry are $\mathbb Z$-waists $waist_k(X;\mathbb Z)$ defined  with   families of $\mathbb Z$-cycles that are assumed as regular as one wishes.}  of the Riemannian manifolds $(X,\vartheta)$ for all metrics $\vartheta$ {\it compatible}  with $\omega$. 
  
 It is easy to see that 
  
 \hspace {9mm}  the space $\Theta=\Theta(\omega)$  of metric $\vartheta$ compatible with $\omega$ is  {\it contractible,}\vspace {0.7mm}

  and that 
  
  \hspace {9mm}  $\omega$-compatible  metrics are {\it extendable} from open subsets $X_0\subset X$ 
  
   \hspace {9mm} to all of $X$ with the usual  precautions at the boundaries of $X_0$.\vspace {0.7mm}
  
 It follows,  that the symplectic waists are 
  monotone under equividimensional symplectic embeddings:   
  
   $$sympl\mbox{-}waist_k (X_0)\leq     
   sympl\mbox{-}waist_k(X)$$ for all open subsets $X_0\subset X$.
   
 If  $k=2$, then upper bounds on symplectic waists are obtained by proving  homological stability of certain families of {\it $\vartheta$-psedoholomorphic curves}  in  $ X$ under deformation of compatible  metrics $\vartheta$, where  "psedoholomorphic curves" are oriented   surfaces $Y\subset X=(X,\omega, \vartheta)$,  such that $area_\vartheta(Y) =\int_Y\omega$.
 
 This, along with the  symplectomorphism of the ball $B=B^{2m}(1)\subset \mathbb C^m$    
onto  the complement $ \mathbb CP^m\setminus  \mathbb CP^{m-1}$, where  the symplectic form  in  $ \mathbb CP^m$  is normalised in order to have the area   of the projective line  $\mathbb CP^1\subset  \mathbb CP^m$  equal  $area(B^2(1))=\pi$, implies that  
   $$  sympl\mbox {-}waist_2(\mathbb CP^m;\mathbb Z) = sympl\mbox {-}waist_2(B;\mathbb Z)=waist_2(B; \mathbb Z)=\pi,$$
where   $waist_2(...; \mathbb Z)$  stands for the  $\mathbb Z$-waist   that is defined with    homologically substantial families of $\mathbb Z$-cycles.
\vspace {1mm}

 \vspace {1mm}

On the other hand, it is probable that     $sympl$-$waist_k(X;\mathbb Z)=\infty$ for all $X$, unless $k= 0, 2, n=dim(X)$ and this is also what one   regretfully expects  to happen to {\it the  symplectic  $k$-systoles}   $syst_k(X,\omega)=_{def}\sup_ \vartheta syst_k(X, \vartheta)$,     for $k\neq 0,2,n$. 

For instance, the complex  projective space $\mathbb CP^m$  may(?) carry Riemannian metrics $\vartheta_s$  for all $s>0$ compatible with the standard symplectic form on $\mathbb CP^m$ such that the $k$-systoles, i.e the minimal $k$-volumes   of all non-homologous zero $k$-cycles in $(\mathbb CP^m,\vartheta_s)$ of all dimensions $k$ except for $0,2,2m$ are $>s$.
 
Yet, some  geometric (topological?)  invariants  of the functions  $\vartheta\mapsto  syst_k(X,\vartheta)$   and $\vartheta\mapsto  waist_k(X,\vartheta)$ on the space of  metrics $\vartheta$ compatible with $\omega $  may shed some light  on the symplectic geometry of $(X, \omega)$, where possible invariants of such a  function $F( \vartheta)$ may be  the asymptotic rate of  some kind of "minimal complexity" of the Riemannian manifolds $(X,\vartheta)$  (e.g.  some integral curvature or  something like  the minimal number of contractible metric  balls needed to cover $(X,\vartheta)$) for which $F( \vartheta)\geq s$, $s\to \infty$.
 
 \vspace {2mm}

Let us generalise the above in the spirit of  "multidimensional spectra"  by  introducing the space $\Psi =\Psi(X,I)$  of triples $\psi=(\{U_i\}_{i\in I}, Y, \vartheta)$, where $U_i$ are disjoint open subsets, $Y\subset X$ is an integer $2$-cycle and $\vartheta$ is a Riemannian metric compatible with $\omega$.

Let ${\cal E}: \Psi\to \mathbb R^{N+1}$ be the map defined by 
$${\cal E}(\{U_i\}_{i\in I}, Y, \vartheta)= (sympl\mbox{-}waist_2 (U_i), area_ \vartheta(Y))_{i=1,2,...,N},$$
where, as an alternative to  $sympl${-}$waist_2(U)$, one may use  $inrad_\omega (U)$, $U\subset (X,\omega)$, that is  the supremum of the radii $r$  of the balls  $B^{2m}(r)\subset \mathbb C^m$ that admit symplectic embedding into $U$.
Let $f: \mathbb R^{N+1}\to \mathbb R$  be a positive function that is symmetric   and monotone decreasing in the first $N$-variables and monotone increasing in the remaining variable (corresponding to $area (Y)$), where the simplest instance of this is $-\sum_ {i=1,...N}z_i+z_{N+1}$

Let $G\subset Sym_N$ be a permutation  group that observe, naturally  acts on $\Psi$,  let $S$ be a topological  space with a free action of $G$ and let   $[\varphi/G]$ be a homotopy class of maps $\varphi :S/G\to \Psi/G$.

 Let $ E_f(\psi)=f({\cal E}(\psi))$ and  define {\it the $\sup\vartheta$-homotopy spectrum} of $E_f$ (compare section 4)  
 $$E_f[\varphi/G]_{\sup_ \vartheta}=\sup_ \vartheta \inf_{\varphi\in [\varphi]} \sup_{s\in S}  E_f(\varphi(s)). $$

\vspace {1mm}

 {\it Playing  $\inf_{\phi}$  Against $ \sup_ {\vartheta}$.} Our major concern here is the  possibility of 
$ E_f[\varphi/G]_{\sup_ \vartheta}=\infty$ that can be outweighed  by enlarging the homotopy  class 
 $[\phi]$ to a  homology class or to a set of such classes. With this in mind, given a cohomology class $h$ in $H^\ast(\Psi/G,\Pi)$ with some coefficients $\Pi$,    one defines  
 $$E_f[h]_{\sup_ \vartheta}=\sup_ \vartheta \inf_{\varphi^\ast(h)\neq 0} \sup_{s\in S}  E_f(\varphi(s)). $$
 where the infimum is taken over all $G$-spaces $S$ and all maps  $\varphi :S/G\to \Psi/G$.
  such that $\varphi^\ast(h)\neq 0$.
 
 Another possible measure against $ E_f[\varphi/G]_{\sup_ \vartheta}=\infty$ is taking $\mathbb Z_2$-cycles $Y$ instead of $\mathbb Z$-cycles, but this is less likely to tip the  balance in our favour. 

On the other hand, one may enlarge/refine  the outcome of minimisation over $\varphi$, yet, still keeping the final result finite,  by restricting the topology of $Y$, e.g.  by allowing only $Y$ represented by surfaces of genera bounded by a given number. Also one may incorporate the integral $\int_Y\omega$ into $E_f$.

 But  usefulness of all these variations is limited by the means we have at our disposal for  proving  finiteness
  of the $\sup_\vartheta$-spectra  that are limited to  the homological (sometimes homotopical) stability of families of psedoholomorphic curves. \vspace {1mm} 
   
{\it Packing Inequalities.} If $X=(X,\omega)$   admits a nontrivial stable family of psedoholomorphis curves $Y\subset X$, then, there are  non-trivial constrains on the topology of the space of packings of $X$   by  $U_i$ with $inrad_\omega (U_i)\geq r$.

Namely there are  connected   $G$-invariant  subsets $S_\circ$ in the space of $I$-tuples of disjoint   topological balls  in $ X$  for all finite sets $I\ni i$  of sufficiently large cardinalities $N$ and some groups $G\subset Sym_N=aut(I)$, such that

{\it every  family  $S$ of  $I$-tuples of disjoint subsets  $U_{i,s}\subset X$, $s\in S$, that is  $G$-equaivariantly homotopic, or just homologous,  to  $S_\circ$ satisfies,
$$\inf_{s\in S}\sum_{i\in I}   sympl\mbox{-}waist_2 (U_{i,s})\leq \int_Y\omega.$$
Consequantly, 
$$\inf_{s\in S}\sum_{i\in I} \leq \pi\cdot inrad_{\omega}(U_{i,s})\leq \int_Y\omega.$$}

 And   effective and rather precise unequalties of this kind   estimates are possible for 
 particular manifolds, say for the projective $\mathbb C P^m$ and domains in it where psedoholomorphic curves are abundant.   

But if $X$ is a closed symplectic manifold $X=(X,\omega)$ {\it with no psedoholomorphic curves in it }, e.g. the $2m$-torus,   one does not know whether  there are  non-trivial constrains on the homotopy   types  of symplectic packing spaces.

Also it is unclear if such an $X$ must have $syst_2(X,\omega)=\infty$ and/or   $sympl\mbox{-}waist_2(X)=\infty$.   

   (See \cite {anjos}  \cite {biran-packings}  \cite  {biran-connectedness}  \cite{buhovski}   \cite {cieliebeck-hofer}   \cite {guth-polydiscs}  \cite {entov}  \cite {mcduff}  \cite {mcduff-polterovich}  \cite {schlenk} for what is known  concerning  {\it individual} symplectic packings of $X$ and  {\it spaces} of embeddings of {\it a single} ball into $X$.)
   
   \vspace{1mm}
   
Conclude this section by observing that   spaces of  certain  symplectic packings can be described entirely in terms of the set   $\Theta=\Theta(\omega)$ of $\omega$-compatible metrics $\vartheta$ on $X$ with the following definition applicable to general classes  $\Theta$ of metric spaces.    \vspace{1mm}
   
  {\it  $\Theta$-Packings by Balls.}    A packings of $(X,\Theta)$ by $I$-tuples of  $r$-balls $B^n(r)\subset \mathbb R^n$ is  a pair $(\vartheta, \{f_i\})$ where $\vartheta \in  \Theta$ and an  $\{f_i\}$, $ i\in I$, is an   $I$-tuple of {\it expanding maps} $f_i:B^n(r)\to (X,\vartheta)$  with disjoint images. \vspace{1mm}
  
 {\it Problem. } Find  further (non-symplectic)    classes $\Theta$ with "interesting" properties of the corresponding spaces of $\Theta$-Packings.
    
    For example,  a Riemannain metric $\vartheta$ may be regarded  as compatible with a {\it pseudo-Riemannian} (i.e. non-degenerate  indefinite) $h$ on a compact manifold $X$ if  the $h$-lengths of all curves are bounded in absolute values  by their $\vartheta$-length and if the $\vartheta$-volume of $X$ equals the $h$-volume.  
    
    Are there instances of $(X,h)$  where some $\Theta$-packings and/or  $sup_\vartheta$-invariants carry a non-trivial information about $h$?

    Are there such examples of other  {\it $G$-structures} on certain  $X$ for  groups  $G\subset GL(n)$, $ n=dim (X)$, besides the symplectic and the orthogonal ones,  where  metrics $\vartheta$ serve  for   reductions of   groups $G$ to their maximal compact subgroups?
    
    \section {Packing Manifolds  by $k$-Cycles and $k$-Volume Spectra   of Spaces of Packings.} 
 
 Define {\it an $I$-packing of the space  ${\cal C}_\ast(X;\Pi)=\bigoplus_{k=0,1,...} {\cal C}_k(X;\Pi)$   of cycles with $\Pi$-coefficients} in a Riemannian manifold $X$   as  {\it  an $I$-tuple $\{V_i\}$  of cycles in $X$} with a given {\it lower bound} on {\it some distances}, denoted $"dist"$, between these cycles, say 
 $$"dist"(V_i,V_j)\geq d,$$
or, more generally,
 $$"dist"(V_i,V_j)\geq d_{ij},\mbox  { } i,i \in I.$$

\vspace {1mm}

{\it  Leading Example: $dist_X$-Packings of the Space  ${\cal C}_\ast(X;\Pi)$.}  A significant  instance of "distance between cycles in $X$" is 
 $$dist_X(V,V')=_{def}\inf_{v\in V,v'\in V'}dist_X(v,v'),$$ where we use here the same notation for cycles $V$ and their supports in $X$,   both denoted $ V_i\subset X$.
 
{\it  Question 1.} What are the homotopy/homology properties  of the spaces $$\Psi_{d_{ij}}={\cal C}_\ast(X;\Pi)^I_{<d_{ij}}\subset {\cal C}_\ast(X;\Pi)^I$$   of $k$-tuples of cycles  $V_i\subset X$  that satisfy the inequalities $dist_X(V_i,V_j)\geq d_{ij}$.?
 
 (Here, as at the other  similar packing  occasions, one should think in $Sym_N$-equivaruinat terms that makes sense since the {\it totality} of the spaces ${\cal C}_\ast(X;\Pi)^I_{<d_{ij}}$ for all $d_{ij}>0$ is $Sym_N$-invariant.)

 \vspace {1mm}
 
  {\it Volume Sectra of  $X^I_{> d_{ij}}$.}  One can approach this question from an opposite angle by  looking  at the  $k_N$-volume spectra  of  (the spaces of cycles in the)  Riemannian manifolds  
 $X^I_{> d_{ij}}$ for various $k_N\leq N\cdot dim(X)$, where   $X^I_{> d_{ij}}\subset X^I$ are  defined by the inequalities  $dist_X(x_i,x_j)>d_{ij} $, since mutual distances between   cycles $V_i$ in $X$ can be seen in terms of locations of their Cartesian products in $X^I$, 
as follows:
$$ dist_X(V_i,V_j)\geq d_{ij}\Leftrightarrow V^\times = \bigtimes_{i\in I} V_i \subset X^I_{> d_{ij}}.$$
 
Recall at this point  
   that the  spaces $X^I_{> d_{ij}=d}$ represent packings of $X$ by balls of radii $d/2$.  
 and 
 observe  that  $V^\times \in {\cal C}_\ast(X^I;\Pi^{N \otimes })$ are rather special, namely split, cycles in $X^I$.

Now the above Question 1 comes with the following companion.\vspace {1mm}

 {\it Question 2.} What are  the volume spectra of the spaces   $X^I_{> d_{ij}}$ and how do they depend on $d_{ij}$?

  \vspace {1mm}

Recall the following relation between volume of cycles $W$ and $W'$  in the Euclidean space and distances between them.  \vspace {1mm}

 {\it Gehring's Linking  Volume Inequality}.   Let $W\subset \mathbb R^n$  be a  $k$-dimensional sub(psedo)manifolds of dimension $k$.\vspace {1mm}
 
  Suppose, \vspace {0.5mm}
  
   {\it $W$ is non-homologous to zero in its open $d$-neighbourhood 
 $U_d(W)\subset \mathbb R^n$,}\vspace {0.5mm}

 or, equivalently, there exists
 an  $ (n-k-1)$-dimensional 
 subpseodomanifold 
 
 $W'\subset \mathbb R^n$ 
 that has\vspace {0.5mm}
 
  {\it a  non-zero 
   linking number with $W$ 
   and  such that $dist_{\mathbb R^n}(W,W')\geq d$}.\vspace {0.5mm}

("Linking" 
is understood mod $2$ if $W$
is non-orientable.) \vspace {1mm}

   \hspace {-5mm} Then,  according to the {\it Federer Fleming  Inequality}, (see next section) 
 $$vol_k(W)\geq \varepsilon_n d^{k}, \mbox { } \varepsilon_n>0,$$
 where, moreover, 
 $$ \varepsilon_n= \varepsilon_k=   vol_{k} (S^{k})$$ 
  ($S^{k}=S^k(1)\subset \mathbb R^{k+1}$ denotes the unit sphere.)
 by the Bombieri-Simon  solution of Gehring's Linking problem (see next section).

  \vspace {1mm}

{\it Proof.} Map the Cartesian product  $W\times W'$ to the sphere  $S^{n-1}(d)\subset \mathbb R^n$  of radius $d$ by 
 $$f:(w,w')\mapsto\frac  {d(w-w')}{dist_{\mathbb R^m}(w,w')}$$ 
 and observe that
 
 $\bullet$ the family of the $f$-images  of the  "slices" $W\times w'\subset W\times W'$, $w'\in W'$, in $S^{n-1}(d)$   is {\it homologically substantial} (in the sense of section 6) since the degree of the map $f:W\times W'\to S^{n-1}(d)$ equals the linking number between $W$ and $W'$;
 
 $\bullet$ the map $f$ is {\it distance  decreasing, hence, $k$-volume decreasing} on the 
 "slices" $W\times w'$ for all $w'\in W'$, since  $dist(w,w')\geq d$ for all $w\in W$ and $w'\in W'$.

Therefore, by the definition of {\sl waist} (see section 6) 
$$vol_k (W)\geq wast_k(S^{n-1}(d))$$
where  $wast_k(S^{n-1}(d))
= vol_k(S^k(d))=d^kvol_k(S^k(1))$ by  the sharp  spherical waist inequality (see sections 6 and 19).
QED.
\vspace {1mm}

{\it Linking Waist Inequality.} The above argument  also shows  that whenever a  $k$-cycle  $W\subset \mathbb R^n $ is {\it non-trivially  homologically  linked} with some $W'\subset \mathbb R^N$, and  $dist(W,W')\geq d$, then 
$$waist_l(W)\geq waist_l (S^n(d))=d^lvol_l(S^l)\mbox { for all $l\leq k$.}$$

This  provides  non-trivial constrains on the spaces of packings of $\mathbb R^n$ by $l$-cycles, since $(k-l)$-cycles in the space of $l$-cycles $Y\subset \mathbb R^n$ make $l$-cycles   $W\subset \mathbb R^n$.\vspace {1mm}

Another generalisation of Gehring inequality concerns several cycles, say $Y_i$ linked to some $k$-cycle  $W$. In this case the intersections of $Y_i$ with any   chain implemented
 by a subvariety 
 $V=V^{k+1}\subset \mathbb R^n$ that fills-in $W$, i.e. has $W$ as its boundary, $\partial V=W$, make  packings of this $V$ by $0$-cycles.

 This, applied to the {\it minimal} $V$ filling-in (spanning) $W$, suggests a (sharp?)    lower bound on distances between such $Y_i$ in terms of what happens to the ordinary packings of the ball $B^{k+1}(r)$ which has  $vol_k(\partial B^k)=vol_k(W)$.
 
 \vspace {1mm}

{\it Question 3.} Let $Y_i$, $i\in I$, be $m$-cycles in $ \mathbb R^{2m+1}$, such that  

\hspace {24mm} $vol_m(Y_i)\leq c$ and  $dist(Y_i,Y_j)\geq d$.

What are, roughly, possibilities for the linking matrices $L_{ij}= \#_{link}(Y_i,Y_j)$ of such $Y_i$ depending on $c$ and $d$?

What are the  homotopies/homologies  of spaces of such $I$-tuples of cycles depending on $L_{ij}$?





\section{Appendix:  Volumes, Fillings, Linkings, Systoles and    Waists.}
 
Let us formulate certain mutually interrelated   fillings, linkings and    waists inequalities  extending those presented  above and  earlier in section 6.
\vspace {1mm}

1. {\it Federer-Fleming Filling-by-Mapping (Isoperometric)  Inequality.} (\cite {federer}) Let $Y\subset \mathbb R^n$ be a closed subset with finite $k$-dimensional Hausdorff measure for an integer $k\leq n$.  Then there exisists a  continuous map $f:Y\to \mathbb R^n$ with the following properties.

   $\bullet_{k-1}$  The image $f(Y)\subset \mathbb R^n$ is {\it at most $(k-1)$-dimensional. }Moreover, $f(Y)$  is contained in {\it a piecewise linear subset $\Sigma^{k-1}= \Sigma^{k-1}(X)\subset \mathbb R^n$ of dimension} ${k-1}$.
 
  $\bullet_{disp}$ {\it The displacement of $Y$ by $f$ is bounded in terms of the Hausdorff measure of $Y$} by 
  $$ \sup_{y\in Y}dist_{\mathbb R^n}(f(y), y)\leq const_nHaumes_k(Y)^{\frac {1}{k}}.$$
 
 $\bullet_{vol}$ The $(k+1)$-dimensional measure of {\it the cylinder  $C_f\subset \mathbb R^n$  of the map $f$}    that is the union of the straight segments 
 $[y, f(y)]\subset \mathbb R^n, y\in Y$,   satisfies
   $$Haumes_{k+1}(C_f)\leq const'_n Haumes_k(Y)^{\frac {k+1}{k}}.$$
 
(A possible choice of  $\Sigma^{k-1}$ is the $(k-1)$-skeleton of  a standard decomposition of 
  the Euclidean $n$-space $\mathbb R^n$ into $R$-cubes for 
 $R=C_n Haumes_k(Y)^\frac {1}{k} $  and  a sufficiently large constant $C_n$, where the map $f:Y\to \Sigma^{k-1}$ is obtained by consecutive  radial projections  from  $Y$ intersected with the  $m$-cubes $\square^m$ to the boundaries  $\partial \square^m$ from certain  points in $\square^m$  starting from $m=n$  and up-to $m=k$.)

 It remains unknown if this  holds with $const_k$ instead of   $const_n$ but the following    inequality  with $const_k$ is  available.  \vspace {1mm}
   
2. {\it Contraction Inequality.}    Let the  $Y\subset \mathbb R^n$ be a $k$-dimensional polyhedral subset.  Then there exists a  continuous map $f:Y\to \mathbb R^n$ with the following properties.

   $\bullet_{k-1}$  The image $f(Y)\subset \mathbb R^n$ is   contained in {\it a piecewise linear subset of dimension} ${k-1}$.
 
  $\bullet_{dist}$ The image of $f$ lies within a {\it controlled distance} from $Y$. Namely, 
  $$ \sup_{y\in Y}dist_{\mathbb R^n}( f(y),  Y)\leq const_kvol_k(Y)^{\frac {1}{k}}.$$
 
 $\bullet_{hmt}$ There exists {\it a homotopy  between the idenetity map and $f$}, say $F:Y\times [0,1]\to \mathbb R^n$ with $F_0=id$ and $F_1=f$, such that  the image of $F$ satisfies 
 $$vol_{k+1}(F(Y\times [0,1]))\leq const'_k vol_k(Y)^{\frac {k+1}{k}}.$$

This is proven   appendix 2 in \cite {filling}  for more general  spaces $X$ in the place of $\mathbb R^n$, including all Banach spaces $X$. 
   \vspace {1mm}

  3. {\it Almgen's Sharp Filling (Isoperimetric) Inequality.} Almgen proved in 1986 \cite {almgren-isoperimetric}
the following sharpening of the (non-mapping aspect of)   Federer-Fleming inequity.
 
{\it the volume minimising $(k+1)$-chains $Z$ in Euclidean spaces   satisfies}:
 $$\frac  {vol_{k+1} (Z)}{vol_k(\partial Z)^\frac {k+1}{k}} \leq \frac {vol_{k+1}(B^{k+1}(1))}{ vol_k(S^k(1)) ^\frac {k+1}{k}}=   (k+1)^{-\frac {k+1}{k}} vol_{k+1}(B^{k+1}(1))^{-\frac {1}{k}}, $$
 where $B^{k+1}(1) \subset\mathbb R^{k+1}$  is the unit Euclidean ball and $S^k(1)=B^{k+1}(1)$ is the unit sphere.

 In fact, {\it  Almgren's local-to global variational principle} \cite {almgren-isoperimetric}, \cite {singularities2}   reduces filling bounds in  Romanian manifolds $X$  to lower bounds of the suprema of  mean curvatures of subvarieties  $Y\in X$ in terms of $k$-volumes of $Y$, where such a sharp(!) bound  for $Y\subset \mathbb R^n$ is obtained by Almgren  by reducing it to that  to that for {\it the Gaussian curvature} of the boundary of the convex hull of $Y$.  
 
 \vspace{1mm}

 4.   {\it Divergence  Inequality.} Recall that the    {\it the $k$-divergence} of a vector field   $\delta =\delta_x$ on a Riemannian manifold $X$ is the  function on the tangent $k$-planes in $X$ that equals {\it the $\delta$-derivative of the $k$-volumes} of these $k$-planes. Thus, \vspace {1mm} 
 
 { \it the $\delta$-derivative of the $k$-volume of each $k$-dimensional submanifold  $V\subset X$ moving by the  $\delta$-flow equals $\int_V div_k(\delta)(\tau_v)dv$, for $\tau_v$ denoting the tangent $k$-plane  to $V\subset X$ at $v$.}\vspace {1mm} 

For instance, the $k$-divergence of {\it the standard   Euclidean radial field focused at zero} $\delta_x=\overrightarrow  x=grad ( \frac {1}{2}||x||^2)$ on $\mathbb R^n$, where $\overrightarrow  x$ denotes this very $x$ seen as the tangent vector
 parallely  transported from $0\in \mathbb R^n$  to $x$,   equals the norm  $||x||$. 
 (The $\delta$-flow  here  equals  the homothety $x\mapsto e^tx$, $t\in \mathbb R$.)
 

 \vspace {1mm}
 
  Let  $Z$ be a compact {\it minimal/stationary} $(k+1)$-dimensional  subvariety with boundary $Y=\partial Z$ and let $\nu=\nu_y$ be the unit vector field tangent to $Z$, normal to $Y$ and facing outside $Z$.
 
 Then 
$$ \int _Zdiv_k(\delta)(\tau_z)dz= \int_Y \langle\delta_y,\nu_y\rangle.\leqno {\star_=}$$

If  $Z$ and $Y$ are non-singular, this equality, that follows from the definition of
 "stationary" and the Gauss-Stokes formula, goes back to 19th century,
 while the singular case  is, probably, due to Federer-Fleming (Reifenberger?  Almgren?). \vspace {1mm} 

5. {\it  Isoperimetric Corollary.} If $div_k(\delta)>0$, then every on subset $Z_0\subset Z$ satisfies:
$$  vol_k(Y) \geq vol_{k+1}(Z_0)   \cdot\frac { \inf_{z\in Z_0} div_{k+1}\delta}{ \sup_{y\in Y}  ||\delta_y|| }. \leqno {\star_{\partial \geq }}$$

This inequality may be    applied to a {\it radial field} $\delta $ focused in $X\supset Z$  at some (say non-singular) point $z_0\in Z$ (such a filed is tangent to the geodesics issuing from $z_0$)  in conjunction  with {\it the coarea  inequality} for the intersections $Y_r$ of $Z$
with spheres  in $X$ around $z_0$ of radii $r$,
$$ \int_0^d vol_k(Y_r)dr \geq vol_{k+1}(Z_d),$$ 
  $Z_r\subset Z$ denoting the intersection of $Z$ with the $d$-ball in $X$ around in $z_0$.

For instance, if $X=\mathbb R^n$ and $\delta$ is the    standard radial field focused at some point $z_0\in Z$, then one arrives this way at   the  classical (known since 1960s, 50s, 40s?)\vspace {1mm}

6.  {\it Monotonicity (Isoperimetric) Inequality.}
 $$\frac  {vol_{k+1} (Z_r)}{vol_k(Y_r)^\frac {k+1}{k}} \leq \frac {vol_{k+1}(B^{k+1}(1))}{ vol_k(S^k(1)) ^\frac {k+1}{k}}.$$
 
{\it Remark.}  If $Z$ is volume minimising, rather than being being only "stationary",  this follows from  the above Almgren sharp filling inequality. 
  \vspace{1mm}
  
  {\it Corollary.} If  the boundary $\partial Z$ lies $d$-far from $z_0$, ten the volume of $Z_d$ is bounded from below by that of the Euclidean ball $B^{k+1}(d)=B_{\mathbf { 0}}(\mathbb R^{k+1}, d) \subset\mathbb R^{k+1}$  of radios $d$:
   $$vol_{k+1}(Z_d)\geq vol_{k+1}( B^{k+1}(d)),$$
 provided the boundary $\partial Z$ lies $d$-far from $z_0$.\vspace {1mm}

 7. { \it Application to  Linking.}  The boundary $Y=\partial Z$ (as well as $Y_d=\partial Z_d$)  satisfies
   $$vol_{k}(Y)\geq vol_{k} (S^k(d)), \leqno {\bullet_{\partial \geq}}$$
where $S^k(d)=\partial   B^{k+1}(d)$ is the Euclidean sphere of radios $d$ and where, we keep assuming that  $Y$ lies {\it outside the $d$-ball} in $\mathbb R^n$ around $z_0$.
 
  Indeed, this follows by applying  ${\star_\geq }$ to the radial field that is   focused at $z_0$, that equals the standard one (i.e. $ \overrightarrow {x-z_0}$) {\it inside the ball} $B^n_{z_0}(d)\subset \mathbb R^n$ and that  has  {\it norm (length) equal  $d$  everywhere outside} this ball.   \vspace {1mm}

  Consequently, every (mildly regular) $k$-cycle $Y\subset  \mathbb R^n$ {\it bounds a $(k+1)$-chain $Z$ in its $d$-neigbourhood for $d$ equal the radius of $k$-sphere with volume equal  $vol_k(Y)$.} 
 
 Namely, the solution  $Z$  of {\it the Plateau problem} with boundary $Y$ does the job.
 (This is how Bombieri and Simon solve the Gehring linking problem,  see \cite {bombieri}.)

\vspace {1mm}
   
8. {\it Systolic Inequality.} A simple adjustment of the above  "monotonicity argument"  yields  the following short-cut in the proof of {\it  the systolic unequalty}\footnote {This was proven in \cite {filling} by a reduction to the contraction (filling) inequality  generalized to Banach spaces, see  \cite {intersystolic}, \cite {katz}.}

Let $K=K(\Pi,1)$, $\Pi=\pi_1(K)$, be {\it an aspherical space} and let ${\cal X}={\cal X} (h,R)$, where $h\in H_n(K)$, $R>0$, be the class of all $n$-dimensional pseudo manifolds $X$ with piecewise Riemannian metrics on them along with maps $f:X\to K$, such  that \vspace {0.8mm}

{\it $\bullet_h$  the fundamental  homology class $[X]_n\in H_n(X)$ of $X$ (defined mod 2 if $X$ is non-orientable) goes 
to a given non-zero homology class   $h\in H_n(K)$;

$\bullet_R$ the restrictions of  the map $f$ to the $R$-balls in $X$, that are maps   $B_x(R)\to K$, are contractible for all $x\in X$.

Then 
 $$vol(X)\geq \alpha_n R^n \mbox   { for }
\alpha_n=\frac  {(2\sigma_{n-1})^n} {n^n \sigma_{n}^{n-1} },$$
where 
$\sigma_n= \frac{2\pi^\frac{n+1}{2}}{\Gamma(\frac{n+1}{2} )}$
is the volume of the unit sphere $S^n$; thus,  $\alpha_n\sim \frac {(2\sqrt e)^n}{n^n}$, that is    $\alpha_n$ equals $\frac {(2\sqrt e)^n}{n^n}$ plus a subexponential term.} (The expected value of the constant in the systolic inequality is  $\sim \frac {c^n}{n^\frac {n}{2}}$.)

\vspace {1mm}
 
 {\it Proof.} Assume that $X$ is {\it volume minimising}  in $\cal X$, i.e. under conditions $\bullet_h$ and $\bullet_R$.\footnote{This is assumption is justifiable    according  to \cite  {wenger}, or, approximately, that is sufficient for the present purpose, by  the formal (and trivial) argument from section 6  in \cite {filling}.} 
  Then   volume of each $r$-ball $B_x(r)\subset  X$ with $r<R$ is bounded by the volume of {\it the spherical $r$-cone } $B_{\circ r}(S)$ over the $r$-sphere $S=S_x(r)= \partial B_x(r)\in X$.\footnote{
 The spherical  $r$-cone  over a piecewise Riemannian manifold $S$ can be seen  by isometric imbedding from $S$ to the equatorial  sphere $S^{N-1}(r)\subset S^N(r)\subset \mathbb R^{N+1}$ and taking the geodesic cone over $S\subset S^{N-1}$ from a pole in  $S^N\supset S^{N-1}$ for  $B_{\circ r}(S)$.}
 
 In fact, if you cut  $B_x(r)$ from $X$ and attach   $B_{\circ r}(S)$  to $X\setminus  B_x(r)$
by the boundary  $\partial B_{\circ r}(S)=S= \partial B_x(r)$, the resulting space $X'$ will admit map $f': X'\to K$ with the conditions  $\bullet_h$ and $\bullet_R$ satisfied.
 
Therefore,
 $$vol_{n-1}(S_x(r))= \frac {d}{dr}vol(B_x(r)) \geq \beta_n vol(B_x(r))^\frac {n-1}{n}
  \mbox { for }\beta_n=\frac {\sigma_{n-1}}{(\frac {1}{2} \sigma_n)^\frac {n-1}{n}},$$
which implies, by integration over $r\in [0,R]$, the bound $vol(B_x(R))\geq  \alpha_n R^n$. for $\alpha_n=(\beta_n/n)^n. $ QED.

\vspace {1mm}

{\it Remark.} Earlier,  Guth, \cite {guth-systolic} suggested a short proof of a  somewhat improved  systolic inequality, that  says, in particular, that  the  $R$-ball $B_x(R)\subset X$ at some point $x\in X$  has  $vol_n(B_x(R))\geq (4n)^{-n} R^n$, provided  the fundamental class  $[X]_n\in H_n(X)$ equals the product of one dimensional  classes.

 His argument, based on    minimal hypersurfaces and induction on dimension, 
generalises to minimal hypersurfaces with boundaries   as in \cite {gromov-lawson} and yields a similar systolic  inequality 
 for   spaces  $X$  with "sufficiently large" fundamental groups  $\Pi=\pi_1(X)$.
 
  But the proof  of the bound 
  $vol_n(B_x(R))\geq \varepsilon_n R^n$ for general groups  $\Pi$, also  due to Guth (see \cite {guth-balls}  where a more general inequality is proven), is rather complicated and gives smaller  constant  $\varepsilon _n$.\footnote{ Most (all?) {\it known} examples of fundamental groups of closed aspherical manifolds, e.g. those with non-postive curvatures, are "sufficiently large".  But, conjecturally, there are many "non-large" examples.}

\vspace {1mm}

9.  {\it Negative Curvature and  Infinite Dimensions.} The divergence inequality and its corollaries applies to many non-Euclidean spaces $X$, such as    {\it $CAT(\kappa)$-spaces}  with $\kappa\leq 0$, that are 
 
 {\it   complete simply connected, possibly   infinite dimensional spaces, e.g. Riemanian/Hilbertian  manifolds,  with non-positive sectional curvatures $ \leq\kappa$  in the sense of Alexandrov }
 
Albeit   vector fields are not, strictly speaking,  defined in singular spaces,  {\it radial  semigroups of  transformations $X\to X$  with controllably  positive $k$-divergence}
are available in $CAT$-spaces. 
This along with  a solution of Plateau's problem,  (let it be only approximate one) shows that

 {\it  every (mildly regular) $k$-cycle $Y$ in a   {\it $CAT(\kappa)$-space}  $X$  is homologous to zero in its $d$-neighbourhood  $U_d(Y)\subset X$ for $d$ equal the radius of  the corresponding $k$-sphere  in the hyperbolic space with curvature $\kappa$,  
 that is    $S^k(\kappa, d)\subset H^{k+1}(\kappa)$, such that  $vol_k( S^k(d))=vol_k(Y)$.}\vspace {1mm}

 {\it Remark.}  If $X$ is {\it finite dimensional as well as non-singular,}  this also follows from the spherical waist inequality that, in fact, does not need $k$-volume contracting (or expanding) radial fields but rather (controllably) $k$-volume contracting maps from $X$ to the unit tangent spheres  $S_x^{n-1}(1)\subset T_x(X)$, $n=dim(X),$ at all 
 $x\in X$, see Appendix 2  in \cite {filling}.
 But proper setting for such an inequality   for infinite dimensional spheres and in a presence of singularities remains problematic.

\vspace {1mm}

.  

10. {\it Almgren's Inequality for  Curvature $\geq 0$.} Let $X$ be a  complete $n$-dimensioanl  Riemannian manifold $X$ {\it  with non-negative sectional curvatures}  and with {\it strictly positive    volume density at infinity}
$$dens_\infty (X)=_{def}\limsup_{R\to \infty} \frac  {vol_n(B_{x_0}(X, R))}{ vol_n(  B_{\mathbf {0}}( \mathbb R^n,R))}>0.$$

Observe, that since $curv(  X) \geq 0$,  this  $X$  admits a {\it  unique  tangent cone $T_\infty(X)$  at infinity} and  the volume  density  of $X$ at infinity  equals the volume of the init ball in this cone centered at   the apex $o\in T_\infty(X)$, where, recall, the tangent cone  $T_\infty(X)$,  of a metric space $X$ at infinity  is {\it the pointed Hausdorff limit} of the  metric spaces obtained by scaling $X$ by $\varepsilon\to 0$:
 $$T_\infty(X)=\lim_{\varepsilon \to 0}( \varepsilon X=_{def}(X,\varepsilon\cdot  dist_X)).$$ 
  
  Let $Y\subset X$ be a compact $k$-dimensional subvariety such that the distance minimizing segments $[x,y]\subset X$  almost all points  $x\in X$ and $Y$
 have their   $Y$-endpoints  $y$ contained in the $C^{1,Lipshitz}$-regular locus of $Y$ where moreover they are normal to $Y$. (This is automatic for  closed smooth  submanifolds $Y\subset X$ but  we need it  for more general $Y$.)

{\it If the norms of the  mean curvatures of $Y$ at almost all of these points $y\in Y$ are bounded by a constant $M$, then   
  $$dens_\infty(X)\leq \frac {vol_k(Y)} {vol_k(S^k(k/M))}$$
 for the $k$-sphere of radius $k/M$ in  $\mathbb R^{k+1}$ that has  its mean curvature equal $M$.}\vspace {1mm}

{\it Proof.} The volumes of the $R$-tubes  $ U_R(Y)\subset X$ around $Y$ are bounded  by those of  $S^k(k/M)\subset \mathbb R^n$ by {\it the  Hermann  Weyl tube formula} extended   as a (volume comparison) inequality to Riemannian  manifolds  $X$ with $curv\geq \kappa$  by Bujalo  and Heintze-Karcher.  

\vspace {1mm}

  Then it  follows by Almgren's variational local-to-global principle mentioned in the above 3, that  
 \vspace {1mm}

 {\it the volume minimising $(k+1)$-chains $Z$ in $X$ satisfy:
$$\frac  {vol_{k+1} (Z)} {vol_k(\partial Z)^\frac{k+1}{k}}  \leq   (k+1)^{-\frac {k+1}{k}} ( dens_\infty(X)\cdot vol_{k+1}(B^{k+1}(1))^{-\frac {1}{k}}. \leqno {[curv\geq 0]_{almg}}$$}

{\it Shareness of ${[curv\geq 0]_{almg}}$.} This inequality  is {\it sharp}, besides $X=\mathbb R^n$, for certain conical (singular if $\delta<1$)  spaces $X$, in particular, for 
$X=(\mathbb R^{k+1}/\Gamma) \times \mathbb R^{n-k-1}$ for finite isometry group $\Gamma$ acting on  $\mathbb R^{k+1}$.

{\it Linking Corollary.}  The    inequality ${[curv\geq 0]_{almg}}$, combined with the above coaria  inequality,  yields the following  generalisation of the Bombieri-Simon (Gehring)  linking volume inequality. \vspace {1mm}

{\it The $(k+1)$-volume minimising chains $Z$ in  a complete Riemannian manifold $X$ of non-negative curvature   that fill-in a  $k$-cycle  $Y$ in $X$ are contained in the $d$-neighbourhood of $Y\subset X$
for 
$d$   equal the radius of the ball $B_o(T_\infty(X),d)\subset T_\infty(X)$, such that 
$vol_k( \partial  B_o(T_\infty(X),d))=vol_k (Y)$.}

 \vspace {1mm}

11. {\it Convex Functions and   Monotonocity Inequality for $curv\geq 0$.}
Let $X$ be a metric space, let  $x_0\in X$ be  a preferred point in $X$ and $\mu_\bullet$, also written as $d\mu_{x_\bullet}$, be a probability measure on $X$.

Let   $h_{x_0}(x,x_\bullet)$, $x,x_\bullet\in X$,  be defined as 
  $$h_{x_0}(x_\bullet,x)=  \max (0, (-dist(x,x_\bullet)+dist (x_0,x_\bullet))$$
and 
$$H_{x_0, \mu_\bullet}(x)=\int_Xh^2_{x_0}(x,x_\bullet)d\mu_{x_\bullet}.$$

If $X$ is a complete manifold with $curv(X)\geq 0$ and $\mu_{\bullet,i}$ is a sequence of measures with supports tending to infinity (thus weakly convergent to $0$), then the limit  
$H(x)=H_{x_0}(x)$ of the functions $H_{x_0, \mu_{\bullet,i}} $ for $i\to\infty$,  assuming this limit exists, is  a {\it convex} function on $X$  by the old Gromoll-Meyer lemma  on convexity of {\it the Busemann functions}. 

Moreover, \vspace {1mm}

if  $ \mu_{\bullet,i}= \rho_i(x_\bullet) dx$ for  {\it radial  functions} $\rho_i$,  i.e $\rho_i(x_\bullet)= \phi_i(dist (x_0,x_\bullet))$,  for  some real functions $\phi_i$, and {\it if $dens_\infty(X)>0$}, then the corresponding  limit  functions $H(x)$ are {\it strictly convex}.  

 In fact this strictness is  controlled by the density $\delta=dens_\infty(X)>0$ as follows. \vspace{1mm}

 {\it The second derivatives of  these $H$ along all geodesics   in $X$ are 
 bounded from below   by $\varepsilon=\varepsilon_n(\delta)>0$.}\footnote{ This, which is obvious once it has been stated, was pointed out to me in slightly different terms  by Grisha Perelman  along with the following  similar observation that constitutes the geometric core of  {\it the Grove-Petersen finiteness theorem.}\vspace {0.5mm}

 {\it If $curv(X)\geq -1$ and if the volume of the unit ball $B_{x_0}(1) \subset  X$ around  an $x_0\in X$ has volume $\geq \delta$, then  $x_0$ admits a  convex neighbourhood 
 $U_\varepsilon \subset  B_{x_0}(1)$  that contains the  ball $B_{x_0}(\varepsilon )$ for $\varepsilon=\varepsilon_n(\delta)>0$.}}
 
Now, the  strict convexity of the (smoothed if necessary)  function  $H(x)$  can be seen as
a lower bound on the $k$-divergence of the gradient of $H$ and, as in the above 6, we arrive at 
 
 {\it the monotonicity inequality}  of intersections of  stationary $(k+1)$-dimensional
  subvarieties $Z\subset X$  with the  balls $B_{x_0}(r)\subset X$, $x_0\in Z$,
 $$\frac  {vol_{k+1} (Z_r)}{vol_k(\partial Z_r)^\frac {k+1}{k}} \leq const =const_n(\delta)\mbox { for } n=dim(X)\mbox { and } \delta=dens_\infty(X).$$

{\it Remark.} A specific evaluation of this   $const$ depends on   a lower bound on the  scalar products $\langle s_0, s \rangle$ averaged over subsets $U\subset S^{n-1}$ with spherical  measures $\geq \delta\cdot vol_{n-1}(S^{n-1})$ in the tangent spheres $S^{n-1}=S^{n-1}_x\subset \mathbb R^n_x=T_x(X)$. But even the best bound
$$\inf_{s_0\in S^{n-1}}\frac {1}{vol(U)}\int_U \langle s_0, s \rangle ds\geq c(\delta)$$
with sharp $c(\delta)$ does not seem(?) to deliver the sharp constant 
$$const_{almgren}=(k+1)^{-\frac {k+1}{k}} ( dens_\infty(X)\cdot vol_{k+1}(B^{k+1}(1))^{-\frac {1}{k}}.$$

 \vspace {1mm}

11.  {\it On Singularities  and Infinite Dimensions with $curv\geq \kappa$.} Since   singular Alexandrov spaces $X$ with $curv\geq \kappa$  admit  strictly contracting  minus "gradient fields" of  strictly convex functions allow  {\it conical filling} of cycles in balls,   $Y\subset B_{x_0}(R)$, such that $vol_{k+1}(Z)\leq cost\cdot  R\cdot vol_k(Y)$ and, according to  \cite {filling}, \cite {wenger} 
{\it the cone inequality} in a space $X$  implies a non-sharp filling inequality in this $X$.

Probably,  these singular  spaces enjoy  the  filling and waist inequalities  with similarly {\it sharp constants} as their non-singular counterparts, (This is easy for    {\it equidimensional Hausdorff limits} of  non-singular spaces, since waists and filling constants  are Hausdorff continuous in the presence of lower volume bounds.) 

 In fact, one expects a full fledged contravariantly  Hausdorff continuous (i.e. for   collapsing $X_i\to X=X_\infty$)    theory of volume minimizing as well as of {\it quasi-stationary} (similar to quasi-geodesics of Milka-Perelman-Petrunin) subvarieties in $X$.

\vspace {1mm} 

Anther avenue of possible generalisations is that of infinite dimensional (singular if needed)  spaces $X$ with positive curvatures.  Here one is encouraged by {\it the stability of $dens_\infty$}: 
$$ dens_\infty(X)=  dens_\infty(X\times \mathbb R^N)$$
 that suggests a class of infinite dimensional spaces $X$  with "small positive"  curvatures  where  the differentials of various  exponential maps  are isometric up-to  small (trace class or smaller) errors. In this case,  one may try to   defined $dens_\infty$ that would allow one  to   formulate and prove an infinite dimensional counterpart  of  $[curv\geq 0]$.

 
 

     \vspace {1mm}

12. {\it Almgren's Morse Theory for Regular Waists}.  Recall  (see section 6) that $k$-waists  of Riemannian manifolds $X$ are defined via classes $\cal D$ of  diagrams  $D_X=\{X\overset {\chi} \leftarrow \Sigma\overset {\varsigma}\to S\}$ where   $S$ and  $\Sigma$ are  pseudomanifolds with $dim(\Sigma)-dim(S)=k$  that represent {\it homologically substantial}  $S$-families of $k$-cycles $Y_s$ in $X$, that are the $\chi$-images of the pullbacks  $\varsigma^{-1}(s), s\in S$  and where homological substantiality may be understood as non-vanishing of the image of the
fundamental homology  class $h\in H_n(Z)$ under the homomorphism   ${\chi}_\ast:H_\ast(\Sigma)\to H_\ast(X)$. Namely 
$waist_k(X)$ is defined as $$waist_k(X)=\inf_{D_X\in \cal D} \sup_svol_k(Y_s)\mbox {  for $Y_s= \chi( \varsigma^{-1}(s))$}.$$

 If  a class $\cal D$  consists of  diagrams  ${\cal D}_X=\{X\overset {\chi} \leftarrow \Sigma\overset {\varsigma}\to S\}$  with { \it sufficiently regular} maps $\varsigma$ and $\chi$, e.g.   {\it piecewise  real analytic} ones, then the resulting waists, call them {\it regular},   admit {\it rather rough} lower bounds in terms of {\it filling} and of  {\it local contractibility} properties of $X$, where the latter referees to  the range of pairs of  numbers $(r,R)$ such that every $r$-ball in $X$ is contractible in the concentric $R$-ball.

On the other hand,   {\it the sharp} lower  bound on  a regular  $k$-waist  of  the unit  $n$-sphere $S^n$, that (trivially) implies the equality 
   $$reg\mbox{-}waist_k(S^n)=vol_k(S^k)$$   
  can be derived from 
    the  Almgren-Morse theory in  {\it  spaces of rectifiable cycles with flat topologies}, \cite {pits}, \cite {guth-waists}. 
     This theory    implies that $$ reg\mbox {-}waist_k (X)\geq \inf_{M^k\in {\cal MIN}_k} vol_k(M^k)$$ 
for  "$\inf$" taken over all {\it minimal/stationary}  $k$-subvarieties $  M^k \subset X$ and the inequality 
 $$reg\mbox{-}waist_k(S^n)\geq vol_k(S^k)$$   follows from the  lower volume bound  for minimal/stationary subvarities in the  spheres $S^n$:
 $$vol_k(M_k) \geq  vol_k(S^k).  \leqno { [volmin_k(S^n)]}$$

    Almgren's theory, that  preceded the  homological localisation method from \cite {waists},  
 albeit limited to the regular case,  
    has an  advantage over the lower bounds on the   $\mathbb Z_2$-waists indicated in section 6  of being applicable  to the integer and to $\mathbb Z_p$-cycles, that allows minimisation over the  maps $\chi:\Sigma\to X$ of non-zero {\it integer degree} in the case of oriented $X$, $S$ and $\Sigma$.\footnote {I am not certain in all implications of Almgren's theory, since I have not studied    the technicalities of this theory  in  detail.}     
   
   Besides, 
Almgren's theory  plus Weyl-Buyalo-Hentze-Karcher volume tube  bound yield the following\vspace {1mm} 

\hspace {-4mm}{\it sharp waist inequality for  closed  Riemannian manifolds $X$ with $curv(X)\geq 1$}, (see   section 3.5 in \cite {singularities2} and  \cite {memarian}).
     $$ waist_k (X)\geq vol_k(S^k)\frac{vol_n(X)}{vol_n(S^n)},\mbox { } n=dim(X).\leqno{[waist]_{curv\geq 1}}$$
     
     \vspace {1mm}

If  $X$ is a manifold with a   {\it convex boundary}, this inequality may be applied to the (smoothed) double of $X$. but the resulting low bound on $waist(X)$ is non-sharp, unlike those obtained by the convex partitioning argument from \cite {waists}. 

Nevertheless, Almgern's Morse theory seems better suitable  for proving the counterpart of $[waist]_{curv\geq 1}$ for  general singular Alexandrov spaces with curvatures $\geq 1$
(where  even a properly formulated volume tube bound is still unavailablele.)     \vspace {1mm}

(13) {\it Regularization Conjecture.} Probably,   general homologically substantial  families of virtually $k$-dimensional subsets $Y_s\subset X$, should admit    an $\varepsilon$-approximation, for all  $\varepsilon>0$,   by {\it regular} nomologically substantial families $Y_{s,\varepsilon}\subset X$ with $$Haumes_k(Y_{s,\varepsilon})\leq Haumes_k(Y_{s})+\varepsilon.$$

This would imply that   all kinds of $k$-waists of $S^n$ are equal to $vol_k(S^k)$ with  no  regularity assumption, that is for  all homologically substantial diagrams with continuous $\chi$ and $\varsigma$ and  with $k$-volumes understood as Hausdorff measures.
(This remains unknown except for $k+1$ and $k=n-1$.)
 \vspace {1mm}

14. {\it Waists at Infinity.} Let $X$ be a complete manifold with $curv(X)\geq 0$. Do the waists of the complements to the  balls $B_{x_0}(R)\subset X$
satisfy 

\hspace {20mm} $waist_k(X\setminus B_{x_0}(R)\geq R^k  dens_\infty(X)vol_k(S^k)$?

 \vspace {0.5mm}

Are there some non-trivial (concavity, monotonicity) inequalities between these waists for different values of $R$?

Does  the above  lower  bound hold for waists of the subsets   $X\setminus B_{x}(R)$  seen as variable ones (in the sense of section 8) parametrised by $X\ni x$?\vspace {1mm}

If "yes" this  would imply that every $m$-dimensional subvariety   $W\subset X$ that is {\it not homologous to zero in its $R$-neighbourhood}, satisfies (compare section  18)

\hspace {25mm}$waist_k(W)\geq R^k  dens_\infty(X)vol_k(S^k)$.\vspace {1mm}

Another class of spaces  $X$ where evaluation of waists (and volume spectra in general) at infinity may be instructive is that of {\it symmetric spaces with non-positive curvatures}.

(One may  define "waists at infinity" via integrals of  positive radial functions $\phi(dist(x,x_0))$   over $k$-cycles in $X$.  For instance, the function  $\phi (d)= \exp\lambda d$, $  \lambda<0$, may serve better than $+\infty$ on the ball $B_{x_0}$ and 1 outside this ball which depicts $X\setminus B(R)$.)\vspace {1mm}

\vspace {1mm} 

(15) {\it Waists of Products and Fibrations.} Let $X$ be a Riemannian product,  $X=\underline X\times X_\varepsilon$, or more generally, let  $X=X(\varepsilon)$ be fibered over $\underline X $.

If  $X_\varepsilon$ is sufficiently small compared to $\underline X$, e.g. $$X_\varepsilon=\varepsilon\cdot  X_0 =_{def} (X_0, dist_\varepsilon=\varepsilon\cdot dist_0)\mbox  { for a small $\varepsilon>0$,}$$
then, conjecturally, 

$$wast_k(X)= wast_k(\underline X)$$
and something similar is expected for  fibrations $X(\varepsilon) \to \underline X$ with {\it small}   fibres $X_{\underline x}\subset X$ that also for this purpose  must vary {\it slowly} as functions of $\underline  x$.

Here, conjecturally, the Hausdorff  $k$-waist of $X(\varepsilon)$   equals
the maximum of  the $k$-volumes of the fibres of the normal projection  $X(\varepsilon)\supset \underline X$ plus a lower order term.\footnote{ This would imply the sharp lower  bound on the Hausdorff  waist of spheres by the homological localisation argument in \cite {waists}. }

This, {\it the regular case,} follows from the Almgren-Morse theory whenever "small"  minimal subvarieties  $M^k\subset X(\varepsilon)$ are "vertical", namely, if  those  with $vol_k(M^k)\leq  wast_k(X(\varepsilon))$, go  to points under  $X(\varepsilon) \to \underline X$, where a sufficient condition to such  "verticality"  is a presence of {\it contracting} vector fields, in sufficiently large balls in $\underline X$ and where such filds often coms as gradients of convex functions. But it is unclear how to sharply  bound from below the  Hausdorff waists in the non-regular case,
compare   7.4 in \cite {waists} and  1.5(B) in \cite {plateau-stein}.

\vspace {1mm}

16. {\it Waists of Thin  Convex Sets.}  A particular case of interest is where $\underline X$ is a compact  convex $m$-dimensional subset in an $n$-dimensional  space, call it  $Z$, of constant curvature,  and $X(\varepsilon)\supset \underline X$ are  $n$-dimensional convex subsets in $Z$ that are $\varepsilon$-close to this $\underline X$.

These $X(\varepsilon)$ {\it normally} project to $\underline X$ with convex $\varepsilon$ small $k$-dimensional fibres $X_{\underline x}$ over  the  interior  points $\underline x\in int(\underline X)= \underline X\setminus \partial \underline X$  for $k=n-m$\footnote {The fibres over the  boundary points  $\underline x\in\partial \underline X$ have 
$dim( X_{\underline x})=n$.} and, conjecturally,
$$\frac { wast_{k}(X(\varepsilon))}{\max_{\underline x\in int(\underline X)}vol_{k}(X_{\underline x})}\to 1 \mbox { for } \varepsilon \to 0.$$

This is known for $\delta$-$Mink_k$-waists,   and  this implies, in conjunction with homological localisation, the   
 sharp lower bounds on the Minkowski waists  spheres, see \cite {waists}. But the case of Hausdorff waistes with {\it no regularity} assumption remains open. \vspace {1mm}

17. {\it Waists of Solids.} Does the regular  $k$-waist of the rectangular  solid 
$$[0,l_1]\times [0,l_2]\times ...\times [0,l_n], \mbox { } 0<l_1\leq l_2\leq...\leq l_i\leq... \leq l_n,$$
equal $l_1\cdot l_2\cdot ...\cdot l_k$?

This is not hard to show  for {\it fast growing} sequences    $l_1<< l_2<<...<< l_i<<... << l_n$
but  the case of  roughly equal $l_i$, especially   that of  the cube $[0,l]^n$, remains problematic.
\vspace {1mm}

18.{\it  Waist of the  Infinite Dimensional  Hilbertian Sphere.}  "Homologically substantial" families of $k$-cycles in $S^\infty$ may be defined with either Fredholm maps $S^\infty\mathbb R^\infty$ of index $k$,  (i.e. with virtually $k$-dimensilnal fibres)  or,   with Fredholm maps $Y\times S^\infty\to S^\infty$ of Fredholm degree one, or by bringing the two diagrammatically  together  as in section 6.

But it is unclear if these waists are  not equal to zero.

 \vspace {1mm}

19. {\it Linking Inequalities with   $\delta$-$Mink_k$.} Such inequalities provide lower bounds on the volumes of $\delta$-neighbourhoods $U_\delta(W)\subset X$ of $k$-dimensional subvarities $W\subset X$   that are not homologous to zero in their $R$-neighbourhoods for some $R>\delta$. For instantce,

If $X=\mathbb R^n$, then 
$$vol_n( U_\delta(W))\geq vol_n( U_\delta(S^k(R)))$$
for the $R$-sphere $S^k(R)\subset \mathbb R^{k+1}\subset \mathbb R^n$.\vspace {1mm}

{\it Proof } The argument  from section18  (also see section 8 in \cite {filling}) with the radial projections from $W$ to the $R$-spheres with centers $x\in \mathbb R^n\setminus W$ applies here, since:

$\bullet$  the spherical waist inequality  holds true with  $\delta$-$Mink_k$ defined  with $\delta$-neighbourhoods of subsets $Y\subset S^{n-1}\subset \mathbb R^n$ taken in 
$\mathbb R^n \supset  S^{n-1}$ (see \cite {waists});

$\bullet$ these projections $W\to S^{n-1}_x(R)$ are not only distance decreasing, but they   are obtained from the identity map by a distance decreasing homotopy. Therefore it diminishes the volume of $\delta$-neighbourhoods by {\it Csik\'os' theorem}. {See \cite{csikos}  and \cite {connely} and references therein.)


 \vspace {1mm}


20. {\it Geometric Linking Inequalities.}  Let two closed subsets $W, W'\subset \mathbb R^n$ of dimensions $k$ and $n-k-1$  {\it can not be unlinked}, i.e.  moved   apart without mutual intersection on the way  by a certain class $\cal M$ of geometric motions. 

Can one bound from below  the Hausdorff measures of these subsets in terms of $dist(W,W')$?

  For instance -- this, according to Eremenko  \cite {eremenko},  was observed by I. Syutric in 1976 or in 1977 --   if $\cal M$ consists of homotheties $h_t: x\mapsto x_t= x_0 +(1-t) (x-x_0)$, $t\in [0,1]$, for some point $x_0\in W$,  then $Haumes_1(W)\geq \pi R$ for connected sets $W$ and 
  $hausmes_1(W)\geq 2\pi R$ for closed curves $W$.

 Indeed,  the image of this homothety, that is  the unit cone over $W$ from $x_0$,  intersects $W'$, say at  $x'\in W'\cap h_{t'}(W) $ for some $t'\in [0,1]$;  thus    the radial projection from $W$ to the $R$-sphere $S^{n-1}_{x'}(R)\subset \mathbb R^n$, that is distance decreasing,  contains two diametrically opposite points, namely the images of $x_0$ and  $x_1= h_{t'}^{-1}(x')\in W$. QED. 
 \vspace {1mm }
 
 {\it A bound on $vol(W\times W')$}. If $W$  and $W'$ can not be unlinked by {\it parallel translations}, then, obviously,  the map $W'\times W\to S_0^{n-1}(R)$  for $(w,w')\mapsto  R\frac {w-w'}{dist( w,w')} $ is onto. Hence, 
 $$Hausmes_k(W)\cdot  Hausmes_{n-k-1}(W)\geq const_nR^{n-1}.$$
  
Thus, either $W$ or $W'$ must have large Hausdorff measure, but one of them may arbitrarily small.  This is seen by taking (large)  $W$ that $\varepsilon$-approximates the $k$-skeleton of a (large) sphere $S^{n-1}(R+r)\subset \mathbb R^n$ and where $W'$ much finer, say with   $\varepsilon'= 0.1 \varepsilon$  approximates   the $(n-k-1)$-skeleton of a  (tiny)  $r'$-ball  $B^n_0(r')\subset   B^n_0(R+r) \subset \mathbb R^n$ for $r'\leq r/2$.

 If $r'\geq 2 \varepsilon$ then $W$  "cages" $W'$ inside $S^{n-1}(R+r)$, provided $r'\geq 10\varepsilon$:
 
  $W'$ can not be moved outside of $S^{n-1}(R+r)$ by an isometric motion without meeting $W$ on the way.

  \vspace {1mm}
  
  {\it Non-accessible  Articles.} There is a dozen or so other papers on Gehring linking problem but,  since they are  not openly accessible, one can not tell what is written in there.


\end {document}